\input ./style/arxiv-general.cfg
\documentclass[aop,MSNbibl,seceqn,dvips]{arximspdf}
\makeatletter
   \@ifpackageloaded{graphicx}{}{\usepackage{graphicx}}
\makeatother
\doi{10.1214/15-AOP1006}% Updated by VTEXPTS2LaTeX.exe, 23.03.2015 13:59
\volume{44}
\issue{2}
\pubyear{2016}
\firstpage{1535}
\lastpage{1598}
\docsubty{FLA}

\makeatletter
\renewcommand{\mid}{|}

\newcommand{\rrvert}{\vert}
\newcommand{\rrVert}{\Vert}
\newcommand{\llvert}{\vert}
\newcommand{\llVert}{\Vert}
\renewcommand{\mid}{|}
\newtheorem{theorem}{Theorem}[section]
\newtheorem{corollary}[theorem]{Corollary}
\newtheorem{lemma}[theorem]{Lemma}
\newtheorem{proposition}[theorem]{Proposition}
\newproclaim{remark}[theorem]{Remark}
\def\R{\mathbb{R}}

\def\Z{\mathbb{Z}}

\def\E{\mathbb{E}}
\def\P{\mathbb{P}}
\def\0{\mathbf{0}}
\def\1{\mathbf{1}}
\def\Var{\mathop{\operatorname{Var}}}
\def\Cov{\mathop{\operatorname{Cov}}}
\makeatother

\begin{document}
\begin{frontmatter}

\title{Spatial asymptotics for the parabolic Anderson models with
generalized time--space Gaussian noise\thanksref{T1}}
\runtitle{Parabolic Anderson models}

\begin{aug}
% Corresponding author: Xia Chen - xchen@math.utk.edu% Updated by
%VTEXPTS2LaTeX.exe, 23.03.2015 13:59
\author[A]{\fnms{Xia}~\snm{Chen}\corref{}\ead[label=e1]{xchen@math.utk.edu}}%,
%\author[]{\fnms{}~\snm{}\ead[label=]{}}
% \and
%\author[]{\fnms{}~\snm{}\ead[label=]{}}
\runauthor{X. Chen}
\affiliation{University of Tennessee}
%\dedicated{}
\address[A]{Department of Mathematics\\
University of Tennessee\\
Knoxville, Tennessee 37996\\
USA\\
\printead{e1}}
\end{aug}
\thankstext{T1}{Supported in part by the Simons Foundation \#244767.}

% HISTORY:
%
\received{\smonth{1} \syear{2014}}% Updated by VTEXPTS2LaTeX.exe,
%23.03.2015 13:59
%
\revised{\smonth{1} \syear{2015}}% Updated by VTEXPTS2LaTeX.exe, 27.03.2015 09:30
%23.03.2015 13:59

% ABSTRACT
%
\begin{abstract}
Partially motivated by the recent papers of
Conus, Joseph and Khoshnevisan
[\textit{Ann. Probab.} \textbf{41} (2013) 2225--2260]
and Conus et al.
[\textit{Probab. Theory Related Fields} \textbf{156} (2013) 483--533],
this work is concerned with the precise spatial asymptotic behavior for
the parabolic Anderson equation
\[
\cases{
\displaystyle {\partial u\over\partial t}(t,x)={1\over2}\Delta
u(t,x)+V(t,x)u(t,x),
\vspace*{3pt}\cr
u(0, x)=u_0(x),}
\]
where the homogeneous
generalized Gaussian noise $V(t,x)$ is, among other forms,
white or fractional white in time and space.
Associated
with the Cole--Hopf solution to the KPZ equation,
in particular, the precise asymptotic form
\[
\lim_{R\to\infty}(\log R)^{-2/3}\log\max_{|x|  \le R}u(t,x) ={3\over4}\root 3\of
{2t\over3}\qquad\mbox{a.s.} %
\]
is obtained for the parabolic Anderson model $\partial_t u
={1\over2}\partial_{xx}^2u+\dot{W}u$
with the $(1+1)$-white noise $\dot{W}(t,x)$. In addition,
some links between
time and space asymptotics for the parabolic Anderson equation
are also pursued.
\end{abstract}

% KEYWORDS
\end{frontmatter}

%s1 #&#
\section{Introduction}\label{intro}

This work is devoted to the analysis of the spatial asymptotics
for the parabolic Anderson model
%
%
%e1.1 #&#
\begin{equation}
\label{intro-1} \cases{ \displaystyle {\partial u\over\partial t}(t,x)=
{1\over2}\Delta u(t,x) +\theta V(t,x)u(t,x),
\vspace*{3pt}\cr
u(0,
x)=u_0(x),}
\end{equation}
where $V(t,x)$ is a centered generalized homogeneous
Gaussian field with the covariance function formally given as
%
%
%e1.2 #&#
\begin{equation}
\label{intro-2}
\qquad \Cov \bigl(V(s, x), V(t,y) \bigr)=\gamma_0(s-t)
\gamma(x-y),\qquad s, t\in\R^+, x, y\in\R^d,
\end{equation}
and $\theta>0$ is a constant playing a role as coefficient.
Some remarkable progress in this direction has been made in recent
papers by Conus, Joseph and Khoshnevisan \cite{CJK} and Conus et al.
\cite{CJKS}
in the case when the time is white, that is, when
$\gamma_0(\cdot)=\delta_0(\cdot)$ (Dirac function) and $\gamma(x)$ takes
a variety of forms. Here we specifically mention the case when
$d=1$, $\gamma_0(u)=\delta_0(u)$
and $\gamma(x)=\delta_0(x)$ in which (\ref{intro-1}) is
formally written as
%
%
%e1.3 #&#
\begin{equation}
\label{intro-7} \cases{ \displaystyle {\partial u\over\partial t}(t,x)=
{1\over2}\partial_{xx}^2u(t,x)+\theta
\dot{W}(t,x) u(t,x),
\vspace*{3pt}\cr
u(0, x)=u_0(x)}
\end{equation}
with $V=\dot{W}$ being a space--time white noise, where $\{W(t,x); t\in\R^+,x\in\R\}$ is
a time--space Brownian sheet. Under the bounded initial condition [given
in (\ref{intro-12}) below],
Conus, Joseph and Khoshnevisan prove (Theorem 1.3, \cite{CJK}) in this case
that
\begin{eqnarray}
\label{intro-4} C^{-1}&\le&\liminf_{R\to\infty}(\log
R)^{-2/3}\log\max_{\llvert  x\rrvert  \le
R}u(t,x)
\nonumber\\[-8pt]\\[-8pt]\nonumber
&\le& \limsup_{R\to\infty}(\log R)^{-2/3}\log\max
_{\llvert  x\rrvert  \le R}u(t,x) \le C\qquad\mbox{a.s.}
\nonumber
\end{eqnarray}
The importance of this result partially lies in the connection (see
\cite{H}) between
(\ref{intro-7}) and
the Kardar--Parisi--Zhang (KPZ) equation
(see \cite{karda-1} and \cite{karda-2} for its
background in the study of interface)
%
%
%e1.4 #&#
%e1.5 #&#
\begin{eqnarray}
\label{intro-8}
{\partial h\over\partial t}(t,x)={1\over2}
{\partial^2 h\over\partial x^2} (t,x)+{1\over2} \biggl({\partial h\over\partial x}(t,x)
\biggr)^2 +\theta\dot{W}(t,x),
\nonumber\\[-8pt]\\[-8pt]
\eqntext{(t,x)\in\R^+\times\R,}
\end{eqnarray}
through the Hopf--Cole transform
%
%
%e1.6 #&#
\begin{equation}
\label{intro-9} u(t,x)=\exp \bigl\{ h(t,x) \bigr\}.
\end{equation}
In particular, (\ref{intro-4}) leads to that
$\max_{\llvert  x\rrvert  \le R}h(t,x)\asymp(\log R)^{2/3}$
($R\to\infty$).

The objectives of this work are set up as follows:

First, we shall install the limits for the asymptotics given in
(\ref{intro-4}) and in some other cases considered in \cite{CJK}
and \cite{CJKS}. Further, we shall identify or compute
the values of these limits.

Second, we shall consider a wider class of Gaussian potentials
where $V(t,x)$ can be white or colored in time. Our first theorem
(Theorem \ref{th-1}) considers
the case
of a general $\gamma_0(\cdot)$ matching with a ``nice'' $\gamma(\cdot)$.
In this paper, however, we are mainly interested in the cases listed in Table~\ref{table1}
where $V(t,x)$ is
fractional (colored)/white in time and space. In the setting labeled
$(1)\times\mathrm{(II)}$, the Gaussian field $V(t,x)$ is formally given
as
%
%
%e1.7 #&#
\begin{equation}
\label{intro-10} V(t,x)=c{\partial^{d+1}W^H(t,x)\over\partial t\,\partial x_1\cdots
\partial x_d},\qquad (t, x_1,
\ldots, x_d)\in\R^+\times\R^d
\end{equation}
and is known as the fractional noise, where $W^H(t,x)$ is a $(d+1)$-parameter
fractional Brownian sheet with the Hurst parameter $H=(H_0,H_1,\ldots, H_d)$.
The settings $(1)\times\mathrm{(III)}$, $(2)\times\mathrm{(II)}$ and
$(2)\times\mathrm{(III)}$ are also interpreted by (\ref{intro-10})
with $H=(H_0, 1/2)$, $H=(1/2, H_1,\ldots, H_d)$ and $H=(1/2, 1/2)$,
respectively. The case of Riesz potential $\gamma(x)=\llvert  x\rrvert
^{-\alpha}$
can be interpreted as a fractional noise with a radially symmetric fractional
spatial component and has close ties to
some classical laws in physics, such as Newton's
gravity law and Coulomb's electrostatics law.

%
%t1 #&#
\begin{table}
\tabcolsep=0pt
\caption{Fractional/white potentials considered in this paper}\label{table1}
\begin{tabular*}{\tablewidth}{@{\extracolsep{\fill}}@{}lccc@{}}
\hline
\textbf{Time/space} & \textbf{(I)} $\bolds{\gamma(x)=\llvert  x\rrvert  ^{-\alpha}}$ & \textbf{(II)} $\bolds{\gamma(x)=\prod_{j=1}^d\llvert  x_j\rrvert  ^{2H_j-2}}$
& \textbf{(III)} $\bolds{\gamma(x)=\delta_0(x)}$\\
\hline
(1) $\gamma_0(\cdot)=\llvert  \cdot \rrvert  ^{-\alpha_0}$ & $\alpha_0\ge0$, & $1/2<H_0\le1$, & $d=1$
\\
\phantom{(1)} $(\alpha_0=2-2H_0)$ & $0<\alpha<d$, & $1/2<H_j<1$ $(1\le j\le d)$,  &
\\
& $2\alpha_0+\alpha<2$ & $2H_0+\sum_{j=1}^dH_j>d+1$ &
\\[3pt]
(2) $\gamma_0(\cdot)=\delta_0(\cdot)$ & $0<\alpha<2\wedge d$ & $1/2 <H_j<1$ $(j=1,\ldots, d)$ & $d=1$\\
\hline
\end{tabular*}
\end{table}

There are some major differences between regime (1) and regime (2)
that lead to difference in treatment between these two regimes.
In regime (1) the solutions $u(t,x)$ have a Feynman--Kac
representation [see (\ref{intro-13}) below] and the solutions in regime (2)
do not. On the other hand, the solutions in regime (1) do not possess
the martingale structure that is related to the mild representation given
in (\ref{intro-34}) below.

Similar to \cite{CJK} and \cite{CJKS},
we assume in (\ref{intro-1}) that $u_0(\cdot)$
is deterministic with
%
%
%e1.8 #&#
\begin{equation}
\label{intro-12} 0<\inf_{x\in\R^d}u_0(x)\le\sup
_{x\in\R^d}u_0(x)<\infty.
\end{equation}
The major development of this paper involves two independent random
systems: one
is a $d$-dimensional Brownian motion $B(t)$ and the other is
a centered generalized homogeneous
Gaussian field $V(t,x)$. Throughout the paper,
by $\E_x$ and $\P_x$, we mean that, respectively,
the expectation and probability law with respect to the Brownian
motion with $B(0)=x$. $\E$ and $\P$ are introduced
for the expectation and probability law with respect to the Gaussian
field.

%s1.1 #&#
\subsection{Results under Feynman--Kac representation}

Solving equation (\ref{intro-1}) may mean different things under different
definitions of stochastic integrals.
The cases considered in this subsection yield
the Feynman--Kac representation
%
%
%e1.9 #&#
\begin{equation}
\label{intro-13}
\qquad u(t,x)=\E_x \biggl[\exp \biggl\{\theta\int
_0^tV \bigl(t-s, B(s) \bigr)\,ds \biggr\}
u_0 \bigl(B(t) \bigr) \biggr],\qquad x\in\R^d
\end{equation}
for the solution of equation (\ref{intro-1}). When $V(t,x)$
has sufficiently nice trajectories, (\ref{intro-13}) is a well-known
fact. This is not a trivial
matter in our context, as the Gaussian field $V(t,x)$ is not even
(necessarily) point-wise defined.
Mathematically speaking, a generalized
centered
time--space Gaussian field $V$ can be defined as a
random linear operator on a Schwartz
space ${\mathcal S}(\R^+\times\R^d)$ of rapidly decreasing (at $\infty
$) and
infinitely smooth functions $\varphi(t, x)$
on $\R^+\times\R^d$ with $
\lim_{t\to0^+}\partial_t^{(n)}\varphi(t, x)=0$
for all $n=0, 1,\ldots$ such that for each
$\varphi\in{\mathcal S}(\R^+\times\R^d)$, $\langle V, \varphi\rangle$
is a centered normal random variable.
In the settings considered in this paper, there are (probably
generalized) functions $\gamma_0(\cdot)$
on $\R$ and $\gamma(\cdot)$ on $\R^d$ such that for
any $\varphi, \psi\in{\mathcal S}(\R^+\times\R^d)$
%
%
%e1.10 #&#
\begin{eqnarray}
\label{intro-14}
&& \Cov \bigl(\langle V,\varphi\rangle,\langle V,\psi
\rangle \bigr)
\nonumber\\[-8pt]\\[-8pt]\nonumber
&&\qquad =\int_{(\R^+\times\R^d)^2}\gamma_0(s-t)
\gamma(x-y)\varphi(s, x)\psi(t,y) \,ds\,dt\,dx\,dy.
\end{eqnarray}
This relation is formally written as in the form given in (\ref{intro-2}).
Given a probability density $h\in{\mathcal S}(\R^+\times\R^d)$, write
$h_\varepsilon(s,x)=\varepsilon^{-(d+1)}h(\varepsilon^{-1}s, \varepsilon^{-1}x)$.
Notice that $V_\varepsilon(t,x)\equiv
\langle V, h_\varepsilon(t-\cdot, x-\cdot)\rangle$ is a point-wise
defined Gaussian field on $\R^+\times\R^d$ and has a sufficient regularity
if $h(t,x)$ is smooth enough. The time integral in~(\ref{intro-13}) is defined as the ${\mathcal L}^2$-limit
%
%
%e1.11 #&#
\begin{eqnarray}
\label{intro-15}
&& \int_0^tV \bigl(t-s, B(s)
\bigr)\,ds
\nonumber\\[-8pt]\\[-8pt]\nonumber
&&\qquad \buildrel \operatorname{def}\over= \lim_{\varepsilon\to0^+}\int
_0^tV_\varepsilon \bigl(t-s, B(s) \bigr)\,ds
-{\mathcal L}^2(\Omega, {\mathcal A}, \P\otimes\P_x),
\end{eqnarray}
provided that the right-hand side converges in
${\mathcal L}^2(\Omega, {\mathcal A}, \P\otimes\P_x)$.

Conditioning on the Brownian motion, this integral is a centered
Gaussian process
(in $t$) with the conditional variance
%
%
%e1.12 #&#
\begin{equation}
\label{intro-16} \int_0^t\!\!\int
_0^t\gamma_0(r-s)\gamma
\bigl(B(r)-B(s) \bigr)\,dr\,ds\qquad (t\ge0).
\end{equation}
The exponential integrability required by the construction of Feynman--Kac
representation in (\ref{intro-13}) can be established by the conditional
Gaussian property, given the exponential integrability of the
Hamiltonian in (\ref{intro-16}). An interested reader is referred to
\cite{HNS} for details.

Under the usual conditions (satisfied by the theorems in this
subsection), the
Feynman--Kac representation given in
(\ref{intro-13}) is a weak solution (Theorem 4.3, \cite{HNS}) to
(\ref{intro-1})
in the sense that
\begin{eqnarray*}
\int_{\R^d}u(t,x)\varphi(x)\,dx&=&\int_{\R^d}u_0(x)
\varphi(x)\,dx+{1\over2} \int_0^t\!\!\int_{\R^d}u(s,x)\Delta\varphi(x)\,dx\,ds
\\
&&{}+\int_0^t\!\!\int_{\R^d}u(s,x)V(s,x)
\varphi(x)\,dx\,ds
\end{eqnarray*}
for any $C^\infty$ and compactly supported function $\varphi(x)$ on $\R^d$,
where the last term is a Stratonovich stochastic integral (Definition
4.1, \cite{HNS}).

One such case is when $\gamma_0(\cdot)$ satisfies some local integrability
and $\gamma(\cdot)$ satisfies
%
%
%e1.13 #&#
\begin{equation}
\label{intro-18} \int_{\R^d} \bigl(1+\llvert \lambda \rrvert
^\delta \bigr) \hat{\gamma}(\lambda) \,d\lambda<
\infty \qquad\mbox{for some $\delta>0$},
\end{equation}
where $\hat{\gamma}$ represents the Fourier transform
(which is non-negative, and exists possibly in the distributional sense).
%
%
%e1.14 #&#
\begin{equation}
\label{intro-33} \hat{\gamma}(\lambda) =\int_{\R^d}
\gamma(x)e^{i\lambda\cdot x}\,dx,\qquad \lambda\in\R^d.
\end{equation}
By Fourier inversion
\[
\gamma(x)-\gamma(y)=(2\pi)^{-d}\int_{\R^d}
\bigl(e^{-i\lambda\cdot
x}-e^{-i\lambda\cdot y} \bigr) \hat{\gamma}(\lambda)\,d\lambda.
\]
Therefore, under (\ref{intro-18}) $\gamma(\cdot)$ is H\"older continuous
with the exponent $\delta$ given in~(\ref{intro-18}).

%
%th1.1 #&#
\begin{theorem}\label{th-1} Assume that $\gamma(\cdot)$
satisfies (\ref{intro-18}) and
$\gamma_0(\cdot)\ge0$ satisfies
%
%
%e1.15 #&#
\begin{equation}
\label{intro-17} \int_0^t\!\!\int
_0^t\gamma_0(r-s)\,dr\,ds<\infty
\qquad (t>0).
\end{equation}
Then for any $t>0$
the weak solution in (\ref{intro-13}) obeys the asymptotic law
%
%
%e1.16 #&#
\begin{eqnarray}
\label{intro-19}
&& \lim_{R\to\infty}(\log R)^{-1/2}\log\max
_{\llvert  x\rrvert  \le R}u(t,x)
\nonumber\\[-8pt]\\[-8pt]\nonumber
&&\qquad =\theta \biggl(2d\gamma(0)\int_0^t\!\!\int_0^t
\gamma_0(r-s)\,dr\,ds \biggr)^{1/2}\qquad\mbox{a.s.}
\end{eqnarray}
\end{theorem}

Theorem \ref{th-1} here is comparable to Theorem 2.5 in \cite{CJKS}
under a different assumption that appears to be
not so comparable to (\ref{intro-18}).

The other cases are those labeled as (1) in Table~\ref{table1} where the Gaussian
potential is fractional in time.
The legitimacy of the Feynman--Kac representation (\ref{intro-13})
is secured for $(1)\times\mathrm{(II)}$ by Theorem 4.3, \cite{HNS},
for $(1)\times\mathrm{(III)}$
by Theorem 6.2, \cite{HNS}
and for $(1)\times\mathrm{(I)}$
by an obvious
modification of the approach used in \cite{HNS}.\vspace*{1pt}

Let $W^{1,2}(\R^d)$ be the Sobolev space of all functions $g$ on
$\R^d$ such that $g, \nabla g\in\mathcal{L}^2(\R^d)$. Denote
%
%
%e1.17 #&#
\begin{eqnarray}\label{intro-22}
{\mathcal A}_d&=& \biggl\{g(s, x); g(s,\cdot)\in
W^{1,2} \bigl(\R^d \bigr), \int
_{\R^d}g^2(s,x)\,dx=1\hspace*{-25pt}\nonumber
\\
&&\hspace*{7pt} \forall0\le s\le1\mbox{ and } \int_0^1\!\!\int_{\R^d}\bigl\llvert \nabla_x g(s, x) \bigr
\rrvert ^2\,dx\,ds<\infty \biggr\},\hspace*{-25pt}
\nonumber\\[-8pt]\\[-8pt]\nonumber
{\mathcal E}(\alpha_0, d,\gamma)&=&\sup
_{g\in{\mathcal A}_d} \biggl\{\int_0^1\!\!\int_0^1\!\int_{\R^d\times\R^d}
{\gamma(x-y)\over\llvert  r-s\rrvert  ^{\alpha_0}} g^2(s, x)g^2(r, y)\,dx\,dy\,dr\,ds\hspace*{-25pt}
\\
&&\hspace*{123pt}{}-{1\over2}\int_0^1\!\!\int
_{\R^d}\bigl\llvert \nabla_x g(s, x)\bigr\rrvert
^2\,dx\,ds \biggr\}.\hspace*{-25pt}\nonumber
\end{eqnarray}
By Lemma 7.2, \cite{CHSX}, ${\mathcal E}(\alpha_0, d,\gamma)$ is finite
under the assumptions in any of the cases listed in Table~\ref{table1} with
the label (1).

Consistently with the parameter $\alpha$ in
setting $\mathrm{(I)}$, we set
%
%
%e1.18 #&#
\begin{eqnarray}
\label{intro-22prime} \alpha=2d-2\sum_{j=1}^dH_j
\end{eqnarray}
in the setting labeled $\mathrm{(II)}$.
A common property shared by $\mathrm{(I)}$ and $\mathrm{(II)}$
is the spatial scaling
%
%
%e1.19 #&#
\begin{eqnarray}
\label{intro-23prime} \gamma(cx)=c^{-\alpha}\gamma(x),\qquad x\in
\R^d, c>0,
\end{eqnarray}
which plays a major role in the formulation of the following theorem.

%
%th1.2 #&#
\begin{theorem}\label{th-2} In settings
$(1)\times\mathrm{(I)}$ and $(1)\times\mathrm{(II)}$ listed
in Table~\ref{table1},
we have that for any $t>0$, the weak solution in
(\ref{intro-13}) satisfies
%
%
%e1.20 #&#
\begin{eqnarray}
\label{intro-23}
&& \lim_{R\to\infty}(\log R)^{-{2/(4-\alpha)}}\log\max
_{\llvert  x\rrvert  \le R}u(t,x)\nonumber
\\
&&\qquad  ={4-\alpha\over4} \biggl(
{4{\mathcal E}(\alpha_0, d, \gamma)\over2-\alpha} \biggr)^{(2-\alpha)/(4-\alpha)}
\\
&&\quad\qquad{}\times \theta^{4/(4-\alpha)}d^{2/(4-\alpha)}
t^{(4-\alpha-2\alpha_0)/(4-\alpha)} \qquad\mbox{a.s.}\nonumber
\end{eqnarray}
\end{theorem}

Relation (\ref{intro-23prime}) remains valid for the setting of
$\gamma(\cdot)=\delta_0(\cdot)$ and $d=1$ with $\alpha=1$.
Consistently with (\ref{intro-22}), set
\begin{eqnarray*}
{\mathcal E}(\alpha_0, 1,\delta_0)&=&\sup
_{g\in{\mathcal A}_1} \biggl\{\int_0^1\!\!\int_0^1\!\int_{-\infty}^\infty
{g^2(s, x)g^2(r, x)
\over\llvert  r-s\rrvert  ^{\alpha_0}}\,dx\,dr\,ds
\\
&&\hspace*{52pt}{}  -{1\over2}\int_0^1\!\!\int_{-\infty}^\infty\bigl\llvert \nabla_xg(s,
x) \bigr\rrvert ^2\,dx\,ds \biggr\}. %
\end{eqnarray*}

%
%th1.3 #&#
\begin{theorem}\label{th-2prime} In setting
$(1)\times\mathrm{(III)}$, listed
in Table~\ref{table1},
we have that for any $t>0$, the weak solution in
(\ref{intro-13}) satisfies
%
%
%e1.21 #&#
\begin{eqnarray}
\label{intro-24}
&& \lim_{R\to\infty}(\log R)^{-2/3}\log\max
_{\llvert  x\rrvert  \le R}u(t,x)
\nonumber\\[-8pt]\\[-8pt]\nonumber
&&\qquad  ={3\over4}\theta^{4/3}
t^{(3-2\alpha_0)/3}\root3\of{{4{\mathcal E}(\alpha_0, 1,
\delta_0)}} \qquad\mbox{a.s.}
\end{eqnarray}
\end{theorem}

For any $t>0$, let ${\mathcal S}([0,t]\times\R^d)$ be the sub-class of
${\mathcal S}(\R^+\times\R^d)$ consisting of $\varphi$ supported on
$[0, t]$ such that
$ \lim_{s\to t^-}\partial_s^{(n)}\varphi(s,x)=0$ for
$n=0,1,\ldots.$ By comparing the covariance functions one
can see that
\[
\bigl\{ \bigl\langle V, \varphi(t-\cdot, \cdot) \bigr\rangle;
\varphi \in{\mathcal S} \bigl([0,t]\times\R^d \bigr) \bigr\}\buildrel
d\over= \bigl\{\langle V, \varphi\rangle;  \varphi\in{\mathcal S}
\bigl([0,t] \times\R^d \bigr) \bigr\}.
\]
Therefore,
\begin{eqnarray}
\label{intro-24prime}
&&\biggl\{\E_x\exp \biggl\{ \theta\int
_0^tV \bigl(s, B(s) \bigr)\,ds \biggr\}; x\in\R^d \biggr\}
\nonumber\\[-8pt]\\[-8pt]\nonumber
&&\qquad \buildrel d\over= \biggl\{ \E_x\exp \biggl\{\theta\int
_0^tV \bigl(t-s, B(s) \bigr)\,ds \biggr\};
x\in\R^d \biggr\}.
\nonumber
\end{eqnarray}
Together, (\ref{intro-13}), Theorems \ref{th-1}~and~\ref{th-2}
(with $u_0\equiv1$)
lead to the following spatial asymptotics for the models of directed
polymers.

%
%co1.4 #&#
\begin{corollary}\label{co-1} Under the assumption of Theorem \ref{th-1},
\begin{eqnarray}
\label{intro-25} &&\lim_{R\to\infty}(\log R)^{-1/2}\log\max
_{\llvert  x\rrvert  \le R} \E_x\exp \biggl\{\theta\int
_0^tV \bigl(s, B(s) \bigr)\,ds \biggr\}
\nonumber\\[-8pt]\\[-8pt]\nonumber
&&\qquad =\theta \biggl(2d\gamma(0)\int_0^t\!\!\int
_0^t\gamma_0(r-s)\,dr\,ds
\biggr)^{1/2}\qquad\mbox{a.s.}
\nonumber
\end{eqnarray}
Under the assumption of Theorem \ref{th-2},
%
%
%e1.22 #&#
\begin{eqnarray}
\label{intro-26} &&\lim_{R\to\infty}(\log R)^{-{2/(4-\alpha)}}\log\max
_{\llvert  x\rrvert
\le R} \E_x\exp \biggl\{\theta\int
_0^tV \bigl(s, B(s) \bigr)\,ds \biggr\}\nonumber
\\
&&\qquad ={4-\alpha\over4} \biggl({4{\mathcal E}(\alpha_0, d, \gamma)\over
2-\alpha} \biggr)^{(2-\alpha)/(4-\alpha)}
\\
&&\quad\qquad{}\times
\theta^{4/(4-\alpha)}d^{2/(4-\alpha)} t^{(4-\alpha-2\alpha_0)/(4-\alpha)} \qquad\mbox{a.s.}
\nonumber
\end{eqnarray}
Under the assumption of Theorem \ref{th-2prime},
%
%
%e1.23 #&#
\begin{eqnarray}\label{intro-26prime}
&& \lim_{R\to\infty}(\log R)^{-2/3}\log\max
_{\llvert  x\rrvert  \le R} \E_x\exp \biggl\{\theta\int
_0^tV \bigl(s, B(s) \bigr)\,ds \biggr\}
\nonumber\\[-8pt]\\[-8pt]\nonumber
&&\qquad ={3\over4}\theta^{4/3} t^{(3-2\alpha_0)/3}\root3\of{{4{
\mathcal E}(\alpha_0, 1, \delta_0)}} \qquad\mbox{a.s.}
\end{eqnarray}
\end{corollary}

We now consider the special case when $\alpha_0=0$ (equivalently,
$H_0=1$) in Theorems \ref{th-2}~and~\ref{th-2prime}.
The Gaussian potential $V$ is time-independent.
Corresponding to $(1)\times\mathrm{(II)}$, for example,
\[
V(x)=c{\partial^d W^H\over\partial x_1\cdots\partial x_d}(x_1,\ldots, x_d),
\qquad x=(x_1,\ldots, x_d)\in\R^d,
\]
where $W^H(x_1,\ldots, x_d)$
is a spatial fractional Brownian sheet
with the Hurst index $H=(H_1,\ldots, H_d)$ satisfying
$1/2<H_1,\ldots, H_d<1$. As for $(1)\times\mathrm{(III)}$ with $\alpha_0=0$
in Table~\ref{table1}, $V(x)=\dot{W}(x)$, a spatial white noise on $\R$.

Write
\[
{\mathcal F}_d= \biggl\{g\in W^{1,2} \bigl(\R^d
\bigr); \int_{\R
^d}g^2(x)\,dx=1 \biggr\}.
\]
By Lemma \ref{va-3} in the \hyperref[appendix]{Appendix},
when $\alpha_0=0$, ${\mathcal E}(0, d, \gamma)$ becomes
%
%
%e1.24 #&#
\begin{eqnarray}
\label{intro-27} {\mathcal E}(d, \gamma)&\equiv&\sup_{g\in{\mathcal F}_d} \biggl
\{\int_{\R^d\times\R^d}\gamma(x-y) g^2(x)g^2(y)\,dx\,dy
\nonumber\\[-8pt]\\[-8pt]\nonumber
&&\hspace*{87pt}{} -{1\over2}\int_{\R^d}\bigl\llvert \nabla g(x)
\bigr\rrvert ^2\,dx \biggr\}.
\end{eqnarray}
In the special case when $d=1$ and $\gamma(\cdot)=\delta_0(\cdot)$,
%
%
%e1.25 #&#
\begin{equation}
\label{intro-28} {\mathcal E}(1, \delta_0)=\sup_{g\in{\mathcal F}_1}
\biggl\{ \int_{-\infty}^\infty g^4(x)\,dx -
{1\over2}\int_{-\infty}^\infty\bigl\llvert
g'(x)\bigr\rrvert ^2\,dx \biggr\}={1\over6}.
\end{equation}
Indeed, the original version (page 291, \cite{CKS}) of the above
identity is
\[
\sup_{g\in{\mathcal F}_1} \biggl\{ 2\int_{-\infty}^\infty
g^4(x)\,dx -{1\over2}\int_{-\infty}^\infty
\bigl\llvert g'(x)\bigr\rrvert ^2\,dx \biggr\}=
{2\over3}. %
\]
Replacing $g(x)$ by $\sqrt{2}g(2x)$ on the left-hand side, we have that
\begin{eqnarray*}
&& \sup_{g\in{\mathcal F}_1} \biggl\{ 2\int_{-\infty}^\infty
g^4(x)\,dx -{1\over2}\int_{-\infty}^\infty
\bigl\llvert g'(x)\bigr\rrvert ^2\,dx \biggr\}
\\
&&\qquad = 4\sup
_{g\in{\mathcal F}_1} \biggl\{ \int_{-\infty}^\infty
g^4(x)\,dx -{1\over2}\int_{-\infty}^\infty
\bigl\llvert g'(x)\bigr\rrvert ^2\,dx \biggr\}.
\end{eqnarray*}
So we have (\ref{intro-28}).

%
%co1.5 #&#
\begin{corollary}\label{co-2}
When $\gamma(\cdot)$ satisfies the assumptions given in
Theorem \ref{th-1},
%
%
%e1.26 #&#
\begin{eqnarray}
\label{intro-29}
&& \lim_{R\to\infty}(\log R)^{-1/2}\log\max
_{\llvert  x\rrvert  \le R} \E_x\exp \biggl\{\theta\int
_0^tV \bigl(B(s) \bigr)\,ds \biggr\}
\nonumber\\[-8pt]\\[-8pt]\nonumber
&&\qquad =t\theta
\bigl(2d\gamma(0) \bigr)^{1/2}\qquad\mbox{a.s.}
\end{eqnarray}
When $\gamma(\cdot)$ is given in $\mathrm{(I)}$ or $\mathrm{(II)}$
with $0<\alpha<2\wedge d$,
\begin{eqnarray}
\label{intro-30}
&&\lim_{R\to\infty}(\log R)^{-{2/(4-\alpha)}}\log\max
_{\llvert  x\rrvert
\le R} \E_x\exp \biggl\{\theta\int
_0^tV \bigl(B(s) \bigr)\,ds \biggr\}
\nonumber\\[-8pt]\\[-8pt]\nonumber
&&\qquad ={4-\alpha\over4}t \biggl({4{\mathcal E}(d, \gamma)\over2-\alpha}
\biggr)^{(2-\alpha)/(4-\alpha)}\theta^{4/(4-\alpha)}d^{2/(4-\alpha)} \qquad\mbox{a.s.}
\nonumber
\end{eqnarray}
When $d=1$ and $\gamma(x)=\delta_0(x)$,
%
%
%e1.27 #&#
\begin{eqnarray}\label{intro-31}
&& \lim_{R\to\infty}(\log R)^{-2/3}\log\max
_{\llvert  x\rrvert  \le R} \E_x\exp \biggl\{\theta\int
_0^t\dot{W} \bigl(B(s) \bigr)\,ds \biggr\}
\nonumber\\[-8pt]\\[-8pt]\nonumber
&&\qquad = {3t\over4}\theta^{4/3}\root3\of{2\over3}
\qquad\mbox{a.s.},
\end{eqnarray}
where $\dot{W}(x)$ ($-\infty< x<\infty$)
is an 1-dimensional spatial white noise.
\end{corollary}

%s1.2 #&#
\subsection{Results for mild solutions}

We now consider the cases labeled by (2) in Table~\ref{table1}, in which
the Gaussian noise $V(t,x)$ is
white in time.
The Feynman--Kac
representation (\ref{intro-13}) is no longer available as $\gamma
(0)=\infty$.
Indeed, one can easily see
that the Hamiltonian in (\ref{intro-16}) [given as the conditional variance
of the time integral (\ref{intro-15}) that would be conditionally
Gaussian if defined] diverges in this case. In spite of this, the parabolic
Anderson equation (\ref{intro-1}) can be solvable
in a slightly different
sense which is briefly described below; we refer to \cite{Dalang,HY} and \cite{Walsh} for details.

The spatial covariance functions considered here have the representation
\[
\gamma(x)=\int_{\R^d}K(x-y)K(y)\,dy,\qquad x\in
\R^d,
\]
where $K(x)$ is symmetric and nonnegative; see (\ref{lo-12}) below.
Assume that the spatial covariance function $\gamma(\cdot)$ satisfies
$\gamma(\cdot)\ge 0$ and
the Dalang condition
%
%
%e1.28 #&#
\begin{equation}
\label{intro-32} \int_{\R^d}{\hat{\gamma}(\lambda)\over1+\llvert  \lambda \rrvert  ^2}\,d\lambda <
\infty,
\end{equation}
where $\hat{\gamma}(\cdot)$ is the Fourier transform of $\gamma(\cdot
)$; see
(\ref{intro-33}). Notice that $\hat{\gamma}(\cdot)\ge0$ as $\gamma
(\cdot)$ is
nonnegative definite.

Let $W(t,x)$ be a $(d+1)$-parameter Brownian sheet,
and consider the Gaussian
field
%
%
%e1.29 #&#
\begin{equation}
\label{intro-35}
\qquad M_t(\varphi)=\int_0^t\!\!\int_{\R^d} \biggl[\int_{\R^d}
\varphi(y-x)K(y)\,dy \biggr]W(ds\,dx), \qquad \varphi\in{\mathcal S} \bigl(
\R^d \bigr),
\end{equation}
where ${\mathcal S}(\R^d)$ is the Schwartz space of the infinitely smooth
and rapidly decaying functions on $\R^d$.

By the theory of Walsh (Chapter~2, \cite{Walsh})
and Dalang \cite{Dalang}, this field can be extended into a
martingale measure $M(t,A)=M_t(1_A)$ such that up to the distributional
equivalence
\[
\langle V, \varphi\rangle=\int_{\R^+\times\R^d} \varphi(s, x)M(ds\,dx),
\qquad \varphi\in{\mathcal S} \bigl(\R^+\times\R^d \bigr).
\]
By the Dalang--Walsh theory, (\ref{intro-32}) ensures the existence and
uniqueness (with a.s. equivalence) of the solution to the
parabolic Anderson equation in the sense that
%
%
%e1.30 #&#
\begin{equation}
\label{intro-34} u(t,x)=(p_t\ast u_0)+\theta\int_0^t\!\!\int_{\R
^d}p_{t-s}(y-x)u(s,y)M(ds,
dy),
\end{equation}
where $p_t$ is the density function of the $d$-dimensional
Brownian motion $B(t)$, and the stochastic
integral on the right-hand side
is taken in the sense of It\^o--Skorokhod. We point out that $u(t,x)\ge0$
in regime (2), labeled as $u(t,x)\ge0$ in Table~\ref{table1}. Indeed, our claim
follows from the following facts: (1) the uniqueness of solution, which
implies that $u(t, x)\equiv0$ if $u_0(0)\equiv0$; (2) the monotonicity
in initial condition. In comparison to the zero solution,
we conclude the solution $u(t,x)\ge0$ if $u_0(x)>0$. The
monotonicity in initial condition was established in \cite{Mueller}
and \cite{Shiga} in the setting of $(2)\times\mathrm{(III)}$. See
(\ref{lo-1}) below for its generalization to whole regime~(2).

An alternative but equivalent
view is to interpret the product in (\ref{intro-1}) between $V(t,x)$
and $u(t,x)$ as
the Wick product; see \cite{HY} for an over view
of the Wick product. When $\gamma(\cdot)$ is bounded and continuous,
the solution
has a ``renormalized'' Feynman--Kac representation,
%
%
%e1.31 #&#
\begin{equation}
\label{intro-36}
\qquad\quad
u(t,x)= e^{-(\theta^2t/2)\gamma(0)} \E_x\exp \biggl[\exp
\biggl\{\theta \int_0^tV \bigl(t-s, B(s)
\bigr)\,ds \biggr\}u_0 \bigl(B(t) \bigr) \biggr].
\end{equation}
We refer to the argument used in the proof of Theorem 7.2, \cite{HNS}
for a
proof of~(\ref{intro-36}). This representation is no longer valid whenever
$\gamma(0)=\infty$. In the cases labeled (2) in Table~\ref{table1}, however,
the solution
$u(t,x)$ can be obtained as the \mbox{$L^2$-}limit
$ \lim_{\varepsilon\to0^+} u_\varepsilon(t,x)$ of~$u_\varepsilon(t,x)$, represented in
(\ref{intro-36}), that appears as the
solution of~(\ref{intro-1}), with $V(t,x)$ being replaced by the
Gaussian potential $V_\varepsilon(t,x)$
of the modified spatial covariance; see, for example, \cite{HN} for details.

By comparing (\ref{intro-13}) and (\ref{intro-36}),
we observe some obvious
differences between solutions
in the Stratonovich sense (\ref{intro-13}) and in the Skorokhod sense
(\ref{intro-36}).
On the other hand,
the solutions given in (\ref{intro-13}) and
(\ref{intro-36}) follow the same limiting
behavior as that stated in Theorem \ref{th-1} for the case
$\gamma_0(\cdot)=\delta_0(\cdot)$
in which (\ref{intro-19}) becomes
%
%
%e1.32 #&#
\begin{equation}
\label{intro-37} \lim_{R\to\infty}(\log R)^{-1/2}\log\max
_{\llvert  x\rrvert  \le R}u(t,x) =\theta \bigl(2dt\gamma(0) \bigr)^{1/2}
\qquad\mbox{a.s.},
\end{equation}
which is the precise form of the limit law stated in
Theorem 2.5, \cite{CJKS}.

%
%th1.6 #&#
\begin{theorem}\label{th-3} In settings $(2)\times\mathrm{(I)}$
and $(2)\times\mathrm{(II)}$ listed in Table~\ref{table1},
%
%
%e1.33 #&#
\begin{eqnarray}\label{intro-38}
&& \lim_{R\to\infty}(\log R)^{-{2/(4-\alpha)}}\log\max
_{\llvert  x\rrvert  \le R}u(t,x)
\nonumber\\[-8pt]\\[-8pt]\nonumber
&&\qquad  ={4-\alpha\over4} \biggl(
{4t{\mathcal E}(d, \gamma)\over2-\alpha} \biggr)^{(2-\alpha)/(4-\alpha)} \theta^{4/(4-\alpha)}d^{2/(4-\alpha)}
\qquad\mbox{a.s.},
\end{eqnarray}
where ${\mathcal E}(d, \gamma)$ is the variation given in (\ref{intro-27}).
\end{theorem}

%
%th1.7 #&#
\begin{theorem}\label{th-4} When $d=1$, $\gamma_0(\cdot)=\delta_0(\cdot)$
and $\gamma(\cdot)=\delta_0(\cdot)$ [i.e., $(2)\times\mathrm{(III)}$
in Table~\ref{table1}],
%
%
%e1.34 #&#
\begin{equation}
\label{intro-39} \lim_{R\to\infty}(\log R)^{-2/3}\log\max
_{\llvert  x\rrvert  \le R}u(t,x) ={3\over4}\theta^{4/3}
\root3\of{2t\over3}\qquad\mbox{a.s.}
\end{equation}
\end{theorem}

In the context of Theorem \ref{th-4}, the parabolic Anderson equation
(\ref{intro-1}) becomes (\ref{intro-7}), which connects
the KPZ equation given in (\ref{intro-8}) through the Hopf--Cole transform
(\ref{intro-9}) in some proper sense.

%
%co1.8 #&#
\begin{corollary}\label{co-4}
Under the deterministic initial condition
\[
-\infty<\inf_{x\in\R}h_0(x)\le\sup
_{x\in\R}h_0(x)<\infty, %
\]
the Hopf--Cole solution $h(t,x)$ to the KPZ equation in (\ref{intro-8})
satisfies
%
%
%e1.35 #&#
\begin{equation}
\label{intro-40} \lim_{R\to\infty}(\log R)^{-2/3}\max
_{\llvert  x\rrvert  \le R}h(t,x) ={3\over4}\theta^{4/3}
\root3\of{2t\over3}\qquad\mbox{a.s.}
\end{equation}
\end{corollary}

%s1.3 #&#
\subsection{Discussion and comment}

As expected, the spatial asymptotics given in the main theorems
are mainly determined by the spatial covariance function $\gamma(\cdot
)$, and more specifically, by the scaling rate $\alpha$ of $\gamma(\cdot)$
[see (\ref{intro-23prime})] when it comes to the settings in Table~\ref{table1}.
On the other hand, cases (1) and
(2) (labeled in Table~\ref{table1})
require different approaches, as we shall see
in Sections~\ref{lo}~and~\ref{mo}.
The tail probability asymptotics (Theorems \ref{th-4prime}--\ref{th-7})
that support the theorems listed above have their independent values,
so we treat them as a part of major theorems of this paper and list
them in Section~\ref{tail}.

The case of time-independence and the case of white time are two
extremes: the least singular and the most singular, respectively.
As we have seen, the
former is associated to $\alpha_0=0$. Because the Fourier transform of
$\gamma_0=\llvert  u\rrvert  ^{-\alpha_0}$ is
$\hat{\gamma}_0(\lambda)=c(\alpha_0)\llvert  \lambda \rrvert  ^{-(1-\alpha_0)}$
($\lambda\in\R^d$)
and $\hat{\delta}_0(\lambda)=1$, the function
$\gamma_0(\cdot)=\delta_0(\cdot)$ is naturally classified as
the extension of $\gamma_0(\cdot)=\llvert  \cdot \rrvert  ^{-\alpha_0}$
($0\le\alpha_0<1$)
to $\alpha_0=1$.
A big surprise is that
these two extreme settings share the same variation
${\mathcal E}(d,\gamma)$ while
the cases with $0<\alpha_0<1$ are formulated by the different variation
${\mathcal E}(\alpha_0, d, \gamma)$. In view of the moment representation
(\ref{mo-9}) below, and knowing that
the difference between two independent Brownian motions
is a Brownian motion, one would bet on
\[
\sup_{g\in{\mathcal F}_d} \biggl\{\int_{\R^d}
\gamma(x)g^2(x)\,dx-{1\over2}\int_{\R^d}
\bigl\llvert \nabla g(x)\bigr\rrvert ^2\,dx \biggr\} %
\]
rather than
${\mathcal E}(d,\gamma)$, as the variation relevant to Theorem \ref{th-3}
and Theorem \ref{th-4}. See Remark \ref{remark} for the discussion.

A central piece of the approach that allows us to compute
the limit values is the precise high moment asymptotics
\[
\log\E u(t,0)^m\qquad (m\to\infty) %
\]
given in Section~\ref{mo}. For the cases labeled (1) in Table~\ref{table1},
our treatment starts at the moment representation (\ref{mo-1}).
The problem can be essentially reduced to the long-term asymptotics
for the annealed moment, to which some results and ideas developed
in the recent work \cite{CHSX} apply. Perhaps the hardest part of this paper
is the computation of the
high moment, when $V(t,x)$ is white in time [i.e., (2) in Table~\ref{table1}].
Unlike the cases labeled (1) in Table~\ref{table1}, the high moment asymptotics
in~(2) do not agree with the long-term asymptotics
such as $\log\E u(t,0)^2$ ($t\to\infty$) at the constant level;
see Remark \ref{remark}
below for details. Under a proper time scaling,
the problem becomes a combination of
high moment and large time with the ratio $t_m\asymp m^{2/(2-\alpha)}$.
The package of methodology includes the Feynman--Kac type
large deviations for time-dependent additive functionals,
newly developed in \cite{CHSX}, some ideas along the line of
probability in Banach spaces and smooth approximation
at the exponential scale.

Another novelty of the paper comes from the proof of the lower bounds
requested by the major theorems listed above, even with
the large deviation
estimates given in Theorems \ref{th-4prime}--\ref{th-7} below.
The relevant idea is clear and simple
in principle: when the space points
$x_1,\ldots, x_n$ are sufficiently spread out, the random variables
$u(t,x_1),\ldots, u(t,x_n)$ are close to being independent. In practice,
carrying this idea out is not easy at all, as indicated by the delicate
steps taken in
\cite{CJK} and~\cite{CJKS}. We adopt this idea and the estimate of localization
developed in \cite{CJK} and~\cite{CJKS}, and use them in our proof
(given in
Section~\ref{lo-s}) in the setting of Theorems \ref{th-3} and \ref{th-4}.
This treatment does not apply to Theorems \ref{th-1}, \ref{th-2} and
\ref{th-2prime} due to its heavy dependence
on the martingale structure associated to
equation (\ref{intro-34}) that defines the mild solution.
The proof (given in Section~\ref{lo-f}) of the lower bounds for
Theorems \ref{th-1}, \ref{th-2} and \ref{th-2prime}
relies on Gaussian property in a substantial way and appears to
be new in methodology.

In comparing their estimates of
the high moment to the literature on intermittency, Conus et al.
(Remark 9.2, \cite{CJKS}) raise the issue of the link between the time
and spatial asymptotics.
In this paper, we pursue this link on two fronts: the first is the
connection
between the long term asymptotics for the annealed moments of $u(t,0)$
and the high moment asymptotics for $u(t,0)$. Indeed, the main
development of
our argument in Section~\ref{mo} is to utilize the link between
annealed intermittency and high moment asymptotics.
We observe a ``perfect match''
when it comes to the Feynman--Kac solution given
by (\ref{intro-13}) and a small but interesting gap when $V(t,x)$ is white
in time. We refer the reader to Remark \ref{remark} below. Our second
concern is the connection between the quenched spatial asymptotics
and the quenched time asymptotics. In Section~\ref{time} we demonstrate
our finding in the setting of time-independence.

We now comment on the relation between the current paper and
\cite{CJK} and \cite{CJKS}. Whenever possible, we adopt the results and
ideas in \cite{CJK} and \cite{CJKS} to our setting. The list includes
the localization (Section~\ref{lo-s}) and estimate for modulus
continuity (Lemma \ref{le-4}) in the case when $V(t,x)$ is white in time.
Estimation by the martingale bound, a substantial idea
in \cite{CJK} and \cite{CJKS},
does not apply to the setting labeled~(1) in Table~\ref{table1}.
The use of the moment representations (\ref{mo-1}) and (\ref{mo-9}),
together with some newly developed ideas in large deviations
for time--space Hamiltonians, allow
us to obtain a form of high moment asymptotics sharper
than those achieved in \cite{CJK} and \cite{CJKS}. On the
other hand, the dependence on the moment representations
(\ref{mo-1}) and (\ref{mo-9}) makes
our method unsuitable to the nonlinear stochastic heat equations
labeled (SHE) in \cite{CJK} and \cite{CJKS}.

In view of the assumption
(\ref{intro-12}) on the initial condition
and the Feynman--Kac
representation (\ref{intro-13}), we have
\[
\underline{u}_0\E_x\exp \biggl\{\int
_0^tV \bigl(t-s, B(s) \bigr)\,ds \biggr\} \le
u(t,x)\le\overline{u}_0\E_x\exp \biggl\{\int
_0^tV \bigl(t-s, B(s) \bigr)\,ds \biggr\}
\]
in the context of Theorems \ref{th-1}, \ref{th-2} and \ref{th-2prime}
where $ \underline{u}_0=\inf_{x\in\R^d}u_0(x)$ and
$ \overline{u}_0=\sup_{x\in\R^d}u_0(x)$.
Or,
%
%
%e1.36 #&#
\begin{equation}
\label{lo-1prime} \underline{u}_0 \tilde{u}(t,x)\le u(t,x)\le
\overline{u}_0 \tilde{u}(t,x)\qquad\mbox{a.s.},
\end{equation}
where
$\tilde{u}(t,x)$ is the solution of
%
%
%e1.37 #&#
\begin{equation}
\label{lo-2} \cases{
\displaystyle {\partial u\over\partial t}(t,x)={1\over2}
\Delta u(t,x) +V(t,x)u(t,x),
\vspace*{3pt}\cr
u(0, x)=1.}
\end{equation}

Relation (\ref{lo-1prime}) remains in the setting of Theorems \ref{th-3}~and~\ref{th-4}.
Indeed, the monotonicity of the solution of (\ref{intro-1}) in
the initial value $u_0$ was established by Mueller \cite{Mueller}
in the special setting $\gamma(x)=\delta_0(x)$. This should be
true in a more general setting. More precisely, let $\tilde{u}(t,x)$
be the solution of (\ref{intro-1}) in the sense of~(\ref{intro-34})
with $u_0(x)$ being replaced by $\tilde{u}_0(x)$. We claim that
%
%
%e1.38 #&#
\begin{eqnarray}\label{lo-1}
&& \tilde{u}_0(x)\le u_0(x)\qquad
\bigl(\forall x\in\R^d \bigr)
\nonumber\\[-8pt]\\[-8pt]\nonumber
&&\qquad \Longrightarrow\quad
\tilde{u}(t,x)\le u(t,x)\qquad\mbox{a.s. }\forall (t,x)\in\R^+\times
\R^d.
\end{eqnarray}
In fact, this becomes obvious in the case when $\gamma(x)$ is well bounded,
in view of~(\ref{intro-36}). For the cases labeled (2)
in Table~\ref{table1}, it is well
known \cite{HN} that $u(t, x)$ can be obtained as the $L^2$-limit
\[
u(t,x)\equiv\lim_{\varepsilon\to0^+}u_\varepsilon(t,x),
\]
where $u_\varepsilon(t,x)$ is the solution of (\ref{intro-1}) with the modified
Gaussian potential $V_\varepsilon(t,x)$ replacing $V(t,x)$
that justifies the
renormalized Feynman--Kac representation~(\ref{intro-36}). Consequently,
the monotonicity in $u_0(x)$ stated in (\ref{lo-1}) passes to $u(t,x)$
through the limit.

Let $\underline{u}(t,x)$ and $\overline{u}(t,x)$
be the solutions of (\ref{intro-1}), corresponding to
the initial conditions $u_0(x)=\underline{u}_0$ and $u_0(x)=\overline{u}_0$,
respectively. By (\ref{lo-1}),
$\underline{u}(t,x)\le u(t,x)\le\overline{u}(t,x)$ a.s. for every
$(t,x)\in\R^+\times\R^d$.
By the linearity of (\ref{intro-1}),
$\underline{u}(t,x)=\underline{u}_0 \tilde{u}(t,x)$
and $\overline{u}(t,x)=\overline{u}_0 \tilde{u}(t,x)$.
This leads to (\ref{lo-1prime}).

By (\ref{lo-1prime}), it is sufficient to establish
our theorems for $\tilde{u}(t,x)$. Equivalently, we replace (\ref{intro-1})
by (\ref{lo-2}) in the rest of the paper. This reduction
results in the stationarity
of $u(t,x)$ in $x$ which substantially simplifies our argument.

In the following proof, we often treat Theorems \ref{th-2}~and~\ref
{th-2prime} together; likewise, we treat Theorems \ref{th-3} and \ref{th-4}
together. In view of (\ref{intro-28}), Theorems \ref{th-2prime}~and~\ref{th-4} can be viewed as, respectively,
Theorems~\ref{th-2}~and~\ref{th-3} in the special case when $d=1$
and $\alpha=1$. This agreement will be reinforced throughout our argument
in order to have a more uniform treatment.

The rest of the paper is organized as follows: Section~\ref{lo} is
devoted to
the lower bounds for Theorems \ref{th-1}, \ref{th-2}, \ref{th-2prime}, \ref
{th-3} and \ref{th-4}.
Section~\ref{mo} is concerned with the high moment asymptotics
$\log u(t,0)^m$ as $m\to\infty$, which appears to be most critical to
the main
development of this work. In Section~\ref{con}, the modulus continuity
of $u(t,x)$ in $x$ is established. With the high moment asymptotics
and the modulus continuity, we are able to compute the exact tail
asymptotics for $\log u(t,0)$ and $ \log\max_{x\in D}u(t,x)$
(with bounded $D\subset\R^d$) in Section~\ref{tail}, where the upper bounds
requested by Theorems \ref{th-1}, \ref{th-2}, \ref{th-2prime}, \ref{th-3}
and \ref{th-4} are established
as a direct consequence of these tail estimates. In Section~\ref{time},
we compare
the quenched spatial asymptotics to existing quenched time asymptotics
in the case of time-independent Gaussian potential.
Finally, we prove some auxiliary lemmas needed for this paper in the
\hyperref[appendix]{Appendix}.

%s2 #&#
\section{Lower bounds}\label{lo}

The proof of the lower bound for a limit law usually appears to
be the most revealing side.
In this section we establish the lower bounds requested by
Theorems \ref{th-1}, \ref{th-2}, \ref{th-2prime}, \ref{th-3} and \ref{th-4}.

%s2.1 #&#
\subsection{The setting of Theorems \texorpdfstring{\protect\ref{th-1}}{1.1}, \texorpdfstring{\protect\ref{th-2}}{1.2} and \texorpdfstring{\protect\ref{th-2prime}}{1.3}}\label{lo-f}

Let $m=m(R)\ge1$ be an integer valued function satisfying
%
%
%e2.1 #&#
\begin{eqnarray}
\label{lo-3}
m\gg\cases{
\displaystyle \sqrt{\log R}, &\quad in the context of Theorem \ref{th-1},
\vspace*{3pt}\cr
\displaystyle (\log R)^{(2-\alpha)/(4-\alpha)}, &\quad in the context of Theorems \ref{th-2}~and~\ref{th-2prime}}\hspace*{-30pt}
\end{eqnarray}
as $R\to\infty$.
Let $\{B_k(t)\}_{k\ge1}$ be an i.i.d. sequence of $d$-dimensional
Brownian motions
with $B_k(0)=0$ ($k=1, 2,\ldots$). The notation $\E_0$ is extended to
the expectation with respect to $\{B_k(t)\}_{k\ge1}$. Write $\tau_k$
for the exit time of $B_k(s)$,
\[
\tau_k=\inf \bigl\{s\ge0; \bigl\llvert
B_k(s)\bigr\rrvert \ge1 \bigr\},\qquad k=1,2,\ldots. %
\]
In view of (\ref{intro-13}), for any $x\in\R^d$,
\begin{eqnarray*}
u(t,x)^m&=&\E_0\exp \Biggl\{\theta \sum
_{k=1}^m\int_0^tV
\bigl(t-s, x+B_k(s) \bigr)\,ds \Biggr\}
\\
&\ge&\E_0 \Biggl[\exp \Biggl\{\theta \sum
_{k=1}^m\int_0^tV
\bigl(t-s, x+B_k(s) \bigr)\,ds \Biggr\}; \min
_{k\le m}\tau_k\ge t \Biggr]
\\
&\ge&\E_0 \Bigl[\exp \bigl\{\lambda\theta \sqrt{\log
R}S_m(t) \bigr\};\xi_m(t,x) \ge\lambda\sqrt{
\log R}S_m(t),  \min_{k\le m}
\tau_k\ge t \Bigr],
\end{eqnarray*}
where $\lambda>0$ is a constant less than but close to $\sqrt{2d}$,
\begin{eqnarray*}
\xi_m(t,x) &=&\sum_{k=1}^m\int
_0^tV \bigl(t-s, x+B_k(s)
\bigr)\,ds, %
\\
S_m(t) &=& \Biggl(\sum_{j, k=1}^m
\int_0^t\!\!\int_0^t
\gamma_0(r-s) \gamma \bigl(B_j(r)-B_k(s)
\bigr)\,dr\,ds \Biggr)^{1/2}. %
\end{eqnarray*}

Set ${\mathcal N}_R=N\Z^d\cap B(0, R)$, where
$B(0, R)=\{x\in\R^d;\llvert  x\rrvert  \le R\}$
and $N>0$ is large but fixed. We have
\begin{eqnarray*}
&& \max_{\llvert  x\rrvert  \le R}u(t,x)^m
\\
&&\qquad \ge \max_{z\in{\mathcal N}_R}u(t,z)^m
\ge\# ({\mathcal N}_R )^{-1}\sum
_{z\in{\mathcal N}_R}u(t,z)^m
\\
&&\qquad \ge \# ({\mathcal N}_R )^{-1} \E_0 \Bigl[
\exp \bigl\{\lambda\theta\sqrt{\log R}S_m(t) \bigr\};
\\
&&\hspace*{90pt} \max
_{z\in{\mathcal N}_R}\xi_m(t,z) \ge\lambda\sqrt{\log
R}S_m(t),  \min_{k\le m}
\tau_k\ge t \Bigr]
\\
&&\qquad =\# ({\mathcal N}_R )^{-1}\E_0 \Bigl({
\mathcal Z}_m(R) 1 \Bigl\{\max_{z\in{\mathcal N}_R}
\xi_m(t,z) \ge\lambda\sqrt{\log R}S_m(t) \Bigr\} \Bigr),
\end{eqnarray*}
where
\[
{\mathcal Z}_m(R)=\exp \bigl\{\lambda\theta\sqrt{\log
R}S_m(t) \bigr\} 1_{\{  \min_{k\le m}\tau_k\ge t\}}. %
\]

The big power $m$ is set to undo the price $\# ({\mathcal N}_R )^{-1}$
paid for pushing $\max_z$ into the expectation. Indeed,
%
%
%e2.2 #&#
\begin{eqnarray}
\label{lo-4}
&&  \max_{\llvert  x\rrvert  \le R}u(t,x)
\nonumber\\[-8pt]\\[-8pt]\nonumber
&&\qquad \ge\# ({\mathcal
N}_R )^{-1/m} \Bigl\{\E_0 \Bigl({\mathcal
Z}_m(R) 1 \Bigl\{\max_{z\in{\mathcal N}_R}\xi_m(t,z)
\ge\lambda\sqrt{\log R}S_m(t) \Bigr\} \Bigr) \Bigr\}^{1/m}.\hspace*{-30pt}
\end{eqnarray}
Given that $\# ({\mathcal N}_R )\le CR^d$ for a universal $C>0$,
(\ref{lo-3}) implies the bounds
%
%
%e2.3 #&#
\begin{equation}
\label{lo-5} \# ({\mathcal N}_R )^{-1/m}= \cases{ \exp
\bigl\{-o (\sqrt{\log R} ) \bigr\},
\vspace*{2pt}\cr
\qquad\mbox{in the context of Theorem~\ref{th-1},}
\vspace*{4pt}\cr
\displaystyle \exp \bigl\{-o \bigl((\log R)^{2/(4-\alpha)} \bigr) \bigr\},
\vspace*{2pt}\cr
\qquad\mbox{in the context of Theorems \ref{th-2}~and~\ref{th-2prime},}}
\end{equation}
which are negligible in comparison to the asymptotic order we try to establish.

Write
\[
\eta_R= \bigl(\E_0{\mathcal Z}_m(R)
\bigr)^{-1}\E_0 \Bigl({\mathcal Z}_m(R) 1
\Bigl\{\max_{z\in{\mathcal N}_R}\xi_m(t,z) <\lambda\sqrt{\log
R}S_m(t) \Bigr\} \Bigr). %
\]
We have
%
%
%e2.4 #&#
\begin{equation}
\label{lo-6}
\qquad \E_0 \Bigl({\mathcal Z}_m(R) 1 \Bigl\{
\max_{z\in{\mathcal N}_R}\xi_m(t,z) \ge\lambda\sqrt{\log
R}S_m(t) \Bigr\} \Bigr)= \bigl(\E_0{\mathcal
Z}_m(R) \bigr) (1-\eta_R).
\end{equation}
An important step is to establish
%
%
%e2.5 #&#
\begin{eqnarray}
\label{lo-7} \lim_{n\to\infty}\eta_{2^n}=0\qquad\mbox{a.s.}
\end{eqnarray}

For any $\varepsilon>0$,
%
%
%e2.6 #&#
\begin{eqnarray}
\label{lo-8}
\P\{\eta_R\ge\varepsilon\}&\le&\varepsilon^{-1}
\E \eta_R\hspace*{-10pt}\nonumber
\\
&=& \bigl(\varepsilon\E_0{\mathcal Z}_m(R)
\bigr)^{-1}\E_0\otimes\E \Bigl({\mathcal
Z}_m(R) 1 \Bigl\{\max_{z\in{\mathcal N}_R}\xi_m(t,z)
<\lambda\sqrt{\log R}S_m(t) \Bigr\} \Bigr)\hspace*{-10pt}
\\
&=& \bigl(\varepsilon\E_0{\mathcal Z}_m(R)
\bigr)^{-1} \E_0 \Bigl({\mathcal Z}_m(R) \P
\Bigl\{\max_{z\in{\mathcal N}_R}\xi_m(t,z) <\lambda\sqrt{\log
R}S_m(t) \mid {\mathcal B} \Bigr\} \Bigr),\hspace*{-10pt}\nonumber
\end{eqnarray}
where ${\mathcal B}$ is the $\sigma$-algebra generated by the Brownian
motions $\{B_k(s)\}_{k\ge1}$.

Conditioning on the Brownian motions, $\{\xi_m(t,z);
z\in{\mathcal N}_R\}$ is a mean zero, and identically distributed
Gaussian family
with the common (conditional) variance $S_m^2(t)$. Further, for any
$z, z'\in{\mathcal N}_R$,
\begin{eqnarray*}
&& \Cov \bigl(\xi_m(t,z), \xi_m
\bigl(t,z' \bigr) \mid {\mathcal B} \bigr)
\\
&&\qquad =\sum
_{j, k=1}^m\int_0^t\!\!\int_0^t \gamma_0(r-s)\gamma
\bigl( \bigl(z-z' \bigr)+ \bigl(B_j(r)-B_k(s)
\bigr) \bigr) \,dr\,ds. %
\end{eqnarray*}
We now claim that for any $0<\rho<1$, one can take $N>0$ sufficiently
large
so that on the event $  \{\min_{k\le m}\tau_k\ge t \}$,
\begin{eqnarray}
\label{lo-9} && \gamma \bigl( \bigl(z-z' \bigr)+
\bigl(B_j(r)-B_k(s) \bigr) \bigr)
\le\rho \gamma \bigl(B_j(r)-B_k(s) \bigr),
\nonumber\\[-8pt]\\[-8pt]
\eqntext{z, z'\in{\mathcal N}_R,  z\neq z', j, k=1,\ldots, m.}
\end{eqnarray}

Regarding this claim, the setting associated to (II), labeled in
Table~\ref{table1}, is the most delicate case among all due to
the un-boundedness of $\gamma(\cdot)$ on each coordinate super plane,
so we treat it in detail.
Let the independent 1-dimensional Brownian motions $B_j^1(s),\ldots, B_j^d(s)$
be the\vspace*{2pt} components of $B_j(s)$, and write $z=(z_1,\ldots,z_d)$ for
$z\in{\mathcal N}_R$. Set $\alpha_j=2-2H_j$ ($j=1,\ldots, d$).
By\vspace*{1pt} assumption we have that $\alpha_j>0$ ($j=1,\ldots, d$). Write
\[
J \bigl(z, z' \bigr)= \bigl\{1\le i\le d;
z_i=z_i' \bigr\},\qquad
z,z'\in {\mathcal N}_R. %
\]
For $i\notin J(z',z)$,
\[
\bigl\llvert \bigl(z_i-z_i' \bigr)+
B_j^i(r)-B_k^i(s)\bigr\rrvert
\ge N-2\ge{N-2\over2}\bigl\llvert B_j^i(r)-B^i_k(s)
\bigr\rrvert. %
\]
Consequently,
\begin{eqnarray*}
&& \gamma \bigl( \bigl(z-z' \bigr)+ \bigl(B_j(r)-B_k(s)
\bigr) \bigr)
\\
&&\qquad =\prod_{i=1}^d\bigl\llvert
\bigl(z_i-z_i' \bigr)+ B_j^i(r)-B_k^i(s)
\bigr\rrvert ^{-\alpha_i}
\\
&&\qquad \le \biggl({2\over N-2} \biggr)^{\alpha'}\prod
_{i=1}^d\bigl\llvert B_j^i(r)-B_k^i(s)
\bigr\rrvert ^{-\alpha_i}
\\
&&\qquad \le \biggl({2\over N-2}
\biggr)^{ \min_{1\le i\le d} \alpha_i} \gamma \bigl(B_j(r)-B_k(s)
\bigr)
\end{eqnarray*}
for every pair $z, z'\in{\mathcal N}_d$ with $z\neq z'$, where the
last step
follows from
\[
\alpha'\equiv\sum_{i\notin J(z, z')}
\alpha_i \ge\min_{1\le i\le d}\alpha_i.
\]
Thus, our assertion (\ref{lo-9}) holds in setting (II).

The proof of (\ref{lo-9}) in other cases is similar,
but easier, due to the fact
that $ \lim_{\llvert  x\rrvert  \to\infty}\gamma(x)=0$
which is automatic for the type-(I) and type-(III) $\gamma(\cdot)$
and a consequence of
assumption (\ref{intro-18})
and the Fourier inversion
\[
\gamma(x)=(2\pi)^{-d}\int_{\R^d}\hat{\gamma}(
\lambda)e^{-i\lambda\cdot
x}\,d\lambda %
\]
in the setting of Theorem \ref{th-1}, according to Riemann's lemma.

By (\ref{lo-9}),
%
%
%e2.7 #&#
\begin{equation}
\label{lo-10} \Cov \bigl(\xi_m(t,z), \xi_m
\bigl(t,z' \bigr) \mid {\mathcal B} \bigr) \le\rho
S_m^2(t).
\end{equation}

Recall that $\lambda<\sqrt{2d}$. Take $u, \rho>0$ sufficiently small so
\[
{(1+2\rho)(\lambda+u)^2\over2}<d\quad\mbox{and}\quad {u^2\over4\rho}>d+1.
\]

Recall (Lemma 4.2, \cite{C-1})
that for a mean-zero $n$-dimensional Gaussian vector
$(\xi_1,\ldots, \xi_n)$ with
identically distributed components,
\[
\rho\equiv\max_{i\neq j}\bigl\llvert \Cov (\xi_i,
\xi_j)\bigr\rrvert /\Var(\xi_1) \le \frac{1}{2},
\]
and for any $A, B>0$,
\[
\P \Bigl\{\max_{k\le n}\xi_k\le A \Bigr\} \le
\bigl(\P \bigl\{\xi_1 \le\sqrt{1+ 2\rho}(A +B) \bigr\}
\bigr)^n +\P \bigl\{U\ge B/\sqrt{2 \rho \operatorname{Var}(\xi_1)} \bigr
\},
\]
where $U$ is a standard normal random variable. Applying this inequality
conditionally with $A=\lambda S_m(t)\sqrt{\log R}$ and
$B=u S_m(t)\sqrt{\log R}$ and noticing
$S_m^2(t)=\Var (\xi_m(t, 0) \mid  {\mathcal B} )$,
we have
\begin{eqnarray*}
&&\P \Bigl\{\max_{z\in{\mathcal N}_R}\xi_m(t,z) <\lambda\sqrt{
\log R}S_m(t) \mid {\mathcal B} \Bigr\}
\\
&&\qquad \le \bigl(\P \bigl\{U \le\sqrt{1+ 2\rho}(\lambda+u)\sqrt{\log R} \bigr\}
\bigr)^{\#({\mathcal N}_R)} +\P \bigl\{U\ge(u/\sqrt{2\rho})\sqrt{\log R} \bigr\}
\\
&&\qquad \le \bigl(1-\exp \bigl\{-(d-v)\log R \bigr\} \bigr)^{C^{-1}R^d} +\exp \bigl
\{-(d+1)\log R \bigr\}
\\
&&\qquad =\exp \bigl\{- \bigl(1+o(1) \bigr)C^{-1}R^v \bigr
\}+R^{-(d+1)}\le C R^{-(d+1)}
\end{eqnarray*}
for large $R>0$, where $v>0$ is independent of $R$.
Bringing this to (\ref{lo-8}) we have that
$\P\{\eta_R\ge\varepsilon\}\le C\varepsilon^{-1} R^{-(d+1)}$ for large $R$.
In particular, (\ref{lo-7}) holds by the Borel--Cantelli lemma.

By the monotonicity of $ \max_{\llvert  x\rrvert  \le R}u(t,x)$
in $R$
and by (\ref{lo-4}), (\ref{lo-5}) and (\ref{lo-6}), the limit along
$R=2^n$ established in (\ref{lo-7}) is sufficient for the lower bounds
for Theorems \ref{th-1}, \ref{th-2} and \ref{th-2prime}
if we can show that
[recall that $m=m(R)\to\infty$ as $R\to\infty$]
%
%
%e2.8 #&#
\begin{eqnarray}
\label{lo-11}
&& \liminf_{\lambda\to(\sqrt{2d})^-} \liminf_{R\to\infty}m^{-1}(
\log R)^{-1/2}\log\E_0{\mathcal Z}_m(R)
\nonumber\\[-9pt]\\[-9pt]\nonumber
&&\qquad \ge
\theta \biggl(2d\gamma(0)\int_0^t\!\!\int
_0^t\gamma_0(r-s)\,dr\,ds
\biggr)^{1/2}
\end{eqnarray}
in the context of Theorem \ref{th-1} and that
%
%
%e2.9 #&#
\begin{eqnarray}
\label{lo-12prime} &&\liminf_{\lambda\to(\sqrt{2d})^-} \liminf_{R\to\infty}m^{-1}(
\log R)^{-{2/(4-\alpha)}}\log\E _0{\mathcal Z}_m(R)
\nonumber
\\
&&\qquad \ge {4-\alpha\over4} \biggl({4{\mathcal E}(\alpha_0, d, \gamma)\over2-\alpha}
\biggr)^{(2-\alpha)/(4-\alpha)}
\\
&&\quad\qquad{}\times \theta^{4/(4-\alpha)}d^{2/(4-\alpha)} t^{(4-\alpha-2\alpha_0)/(4-\alpha)}
\nonumber
\end{eqnarray}
in the context of Theorems \ref{th-2} and \ref{th-2prime}.
By now, the problem of the almost sure limits has been reduced to
pursuing the deterministic limits. Unlike setting (2),
the discussion here does not rely on, but contributes to the development
in later sections.

We now prove (\ref{lo-11}). By the continuity of $\gamma(x)$ at $x=0$,
given $\varepsilon>0$ one can take
$0<\delta<1$ sufficiently small so that $\gamma(x)\ge\gamma
(0)-\varepsilon$
as long as
$\llvert  x\rrvert  \le2\delta$. Set
\[
\tau_k(\delta)=\inf \bigl\{s\ge0; \bigl\llvert
B_k(s)\bigr\rrvert \ge\delta \bigr\}. %
\]
By the definition of ${\mathcal Z}_m(R)$,
\begin{eqnarray*}
{\mathcal Z}_m(R) &\ge& \exp \biggl\{\lambda\theta m\sqrt{\log R}
\biggl( \bigl(\gamma(0) -\varepsilon \bigr)\int_0^t\!\!\int_0^t\gamma_0(r-s)\,dr\,ds
\biggr)^{1/2} \biggr\}
\\
&&{}\times  \P_0 \Bigl\{ \min
_{k\le m}\tau_k(\delta)\ge t \Bigr\}. %
\end{eqnarray*}
Therefore, (\ref{lo-11}) follows from the bound given by
the following relation:
\[
\P_0 \Bigl\{\min_{k\le m}\tau_k(
\delta)\ge t \Bigr\} = \Bigl(\P_0 \Bigl\{\max_{s\le t}
\bigl\llvert B(s)\bigr\rrvert \le\delta \Bigr\} \Bigr)^m. %
\]

It remains to prove (\ref{lo-12prime}). First, by the Brownian scaling
\begin{eqnarray*}
\E_0{\mathcal Z}_m(R)&=&\E_0 \Biggl[\exp
\Biggl\{\theta\lambda\sqrt{\log R} \Biggl(\sum_{j, k=1}^m
\int_0^t\!\!\int_0^t
{\gamma (B_j(r)-B_k(s) )
\over\llvert  r-s\rrvert  ^{\alpha_0}}\,dr\,ds \Biggr)^{1/2} \Biggr\};
\\
&&\hspace*{233pt}  \min _{k\le m}\tau_k\ge t \Biggr]
\\
&=&\E_0 \Biggl[\exp \Biggl\{\theta\lambda t_R^{\alpha_0/2}
\Biggl(\sum_{j, k=1}^m\int _0^{t_R}\!\!\int_0^{t_R}
{\gamma
(B_j(r)-B_k(s) )
\over\llvert  r-s\rrvert  ^{\alpha_0}}\,dr\,ds \Biggr)^{1/2} \Biggr\};
\\
&&\hspace*{227pt} \min
_{k\le m}\tilde{\tau}_k\ge t_R \Biggr],
\end{eqnarray*}
where
\[
t_R=t^{(4-\alpha-2\alpha_0)/(4-\alpha)}(\log R)^{2/(4-\alpha)}
\]
and
\[
\tilde{\tau}_k=\inf \bigl\{s\ge0;\bigl\llvert
B_k(s)\bigr\rrvert \ge t^{-{\alpha_0/(4-\alpha)}}(\log R)^{1/(4-\alpha)}
\bigr\}. %
\]

Notice the representations
%
%
%e2.10 #&#
\begin{eqnarray}\label{lo-12}
\llvert u\rrvert ^{-\alpha_0} &=& C_0\int
_{\R}\llvert v\rrvert ^{-{(1+\alpha_0)/2}} \llvert v -u \rrvert
^{-{(1+\alpha_0)/2}}\,dv \quad\mbox{and}
\nonumber\\[-8pt]\\[-8pt]\nonumber
\gamma(x)&=& \int _{\R^d}K(y)K(y-x)\,dy,
\end{eqnarray}
where $C_0>0$ is a constant independent of $u$ and
the function $K(x)\ge0$ is a positive constant multiple of
\[
\llvert x\rrvert ^{-{(d+\alpha)/2}},\qquad \prod_{i=1}^d
\llvert x_i\rrvert ^{-{(1+\alpha_i)/2}}\quad\mbox{and}\quad
\delta_0(x), %
\]
in connection to,
respectively, the space covariances $\gamma(\cdot)$
of type-(I), type-(II) and type-(III) (labeled in Table~\ref{table1}). This leads to
\begin{eqnarray*}
&&t_R^{\alpha_0/2} \Biggl(\sum_{j, k=1}^m
\int _0^{t_R}\!\!\int_0^{t_R}
{\gamma (B_j(r)-B_k(s) )
\over\llvert  r-s\rrvert  ^{\alpha_0}}\,dr\,ds \Biggr)^{1/2}
\\
&&\qquad = \Biggl(\sum_{j, k=1}^m\int _0^{t_R}\!\!\int_0^{t_R}
{\gamma (B_j(r)-B_k(s) )
\over\llvert   t_R^{-1}(r-s)\rrvert  ^{\alpha_0}}\,dr\,ds \Biggr)^{1/2}
\\
&&\qquad = \Biggl(C_0\int_{\R\times\R^d} \Biggl[\sum
_{j=1}^m \int_0^{t_R}
\bigl\llvert u-t_R^{-1}s\bigr\rrvert ^{-{(1+\alpha_0)/2}}K
\bigl(x-B_j(s) \bigr)\,ds \Biggr]^2\,du\,dx
\Biggr)^{1/2}.
\end{eqnarray*}
Let $f(u, x)$ be a bounded, continuous and locally
supported function on $\R\times\R^d$ with $\llVert  f\rrVert  _2= 1$. By the
Cauchy--Schwarz inequality, the right-hand side of the above equation
is no less than
\begin{eqnarray*}
&&\sqrt{C_0}\int_{\R\times\R^d}f(u,x) \Biggl[\sum
_{j=1}^m \int_0^{t_R}
\bigl\llvert u-t_R^{-1}s\bigr\rrvert ^{-{(1+\alpha_0)/2}}K
\bigl(x-B_j(s) \bigr)\,ds \Biggr]\,du\,dx
\\
&&\qquad =\sqrt{C_0}\sum_{j=1}^m
\int_0^{t_R}\bar{f} \biggl({s\over t_R},
B_j(s) \biggr)\,ds,
\end{eqnarray*}
where
\[
\bar{f}(s, x)=\int_{\R\times\R^d}f(u, y) \llvert u-s \rrvert
^{-{(1+\alpha_0)/2}}K(y-x)\,du\,dy. %
\]

Summarizing our computation and by independence we have
%
%
%e2.11 #&#
\begin{equation}
\label{lo-13} \qquad\E_0{\mathcal Z}_m(R)\ge \biggl(
\E_0 \biggl[\exp \biggl\{\theta\lambda\sqrt{C_0} \int
_0^{t_R}\bar{f} \biggl({s\over t_R},
B(s) \biggr)\,ds \biggr\}; \tilde{\tau} \ge t_R \biggr]
\biggr)^m.
\end{equation}

According to Proposition 3.1 and (3.18) in
\cite{CHSX}, for a bounded open domain
$D\subset\R^d$ containing 0, and for a bounded
function $h(s,x)$ defined on $[0, 1]\times\R^d$
that is continuous in $x$ and equicontinuous (over $x\in\R^d$) in $s$,
\begin{eqnarray}
\label{lo-14}
&&\lim_{t\to\infty}{1\over t}\log
\E_0 \biggl[\exp \biggl\{ \int_0^{t}h
\biggl({s\over t}, B(s) \biggr)\,ds \biggr\};
\tau_D \ge t \biggr]
\nonumber\\[-8pt]\\[-8pt]\nonumber
&&\qquad =\sup_{g\in{\mathcal A}_d(D)} \biggl\{\int_0^1\!\!\int_D h(s,x)g^2(s,x)\,dx\,ds -
{1\over2}\int_0^1\!\!\int
_D\bigl\llvert \nabla_x g(x)\bigr\rrvert
^2\,dx\,ds \biggr\},\hspace*{-20pt}
\nonumber
\end{eqnarray}
where $\tau_D=\inf\{s\ge0; B(s)\notin D\}$ and
${\mathcal A}_d(D)$ is the subspace of ${\mathcal A}_d$ consisting of the
functions $g(s,x)$ vanishing for $x\notin D$.
Applying this to (\ref{lo-13}) we can get
%
%
%e2.12 #&#
\begin{eqnarray}
\label{lo-15} &&\liminf_{R\to\infty}{1\over mt_R}\log
\E_0{\mathcal Z}_m(R)\nonumber
\\
&&\qquad \ge\sup_{g\in{\mathcal A}_d} \biggl\{\theta\lambda\sqrt{C_0}
\int_0^1\!\!\int_{\R^d}
\bar{f}(s,x)g^2(s,x)\,dx\,ds
\\
&&\hspace*{86pt}{} -{1\over2}\int_0^1\!\!\int_{\R^d}\bigl
\llvert \nabla_x g(s, x) \bigr\rrvert ^2\,dx\,ds \biggr\}.
\nonumber
\end{eqnarray}
By Fubini's theorem
\begin{eqnarray*}
&& \int_0^1\!\!\int_{\R^d}
\bar{f}(s,x)g^2(s,x)\,dx\,ds
\\
&&\qquad  =\int_{\R\times\R^d}f(u,y) \biggl[
\int_0^1\!\!\int_{\R^d}
\llvert u-s\rrvert ^{-{(1+\alpha_0)/2}}K(y-x)g^2(s,x)\,dx\,ds \biggr]\,du\,dy.
\end{eqnarray*}
We now take the supremum over $f$ on the right-hand side of (\ref{lo-15}).
Notice that the supremums over $g$ and over $f$ are commutative and that
for any dense subset set ${\mathcal S}$ of the unit sphere of
${\mathcal L}^2(\R\times\R^d)$, by the Hahn--Banach theorem,
%
%e2.13 #&#
\begin{eqnarray}
\sup_{f\in{\mathcal S}}\int_{\R\times\R^d}f(u,x)h(u,x)\,dx=
\biggl(\int_{\R
\times\R^d} \bigl\llvert h(u,x)\bigr\rrvert
^2 \,dx\,du \biggr)^{1/2},\nonumber
\\
\eqntext{h\in{\mathcal
L}^2 \bigl(\R \times \R^d \bigr). }
\end{eqnarray}
Therefore, the right-hand side of (\ref{lo-15})
becomes
\begin{eqnarray*}
&&\sup_{g\in{\mathcal A}_d} \biggl\{\theta\lambda \biggl(C_0
\int_{\R\times\R^d} \biggl[\int_0^1\!\!\int_{\R^d} \llvert u-s\rrvert ^{-{(1+\alpha_0)/2}}\nonumber
\\
&&\hspace*{122pt}{}\times K(y-x)g^2(s,x)\,dx\,ds
\biggr]^2\,du\,dy \biggr)^{1/2}\nonumber
\\
&&\hspace*{174pt}{} -{1\over2}\int_0^1\!\!\int_{\R^d}\bigl\llvert \nabla_x g(x) \bigr
\rrvert ^2\,dx\,ds \biggr\}
\\
&&\qquad =\sup_{g\in{\mathcal A}_d} \biggl\{\theta\lambda \biggl( \int_0^1\!\!\int_0^1\!\!\int_{\R^d\times\R^d}{\gamma(x-y)
\over\llvert  r-s\rrvert  ^{\alpha_0}}g^2(r,x)g^2(s,y)\,dx\,dy
\biggr)^{1/2} \nonumber
\\
&&\hspace*{176pt}{}-{1\over2}\int_0^1\!\!\int_{\R^d}\bigl\llvert \nabla_x g(x) \bigr
\rrvert ^2\,dx\,ds \biggr\},\nonumber
\end{eqnarray*}
where the equality follows from Fubini's theorem and the relations in
(\ref{lo-12}). By rescaling the space variable (Lemma 4.1, \cite{CHSX})
this variation is further equal to
\begin{eqnarray*}
&&(\theta\lambda)^{4/(4-\alpha)}\sup_{g\in{\mathcal A}_d} \biggl\{ \biggl(
\int_0^1\!\!\int_0^1\!\int_{\R^d\times\R^d}{\gamma(x-y)
\over\llvert  r-s\rrvert  ^{\alpha_0}}g^2(r,x)g^2(s,y)\,dx\,dy
\biggr)^{1/2}
\\
&&\hspace*{185pt}{} -{1\over2}\int_0^1\!\!\int_{\R^d}\bigl\llvert \nabla_x g(x) \bigr
\rrvert ^2\,dx\,ds \biggr\}
\\
&&\qquad =(\theta\lambda)^{4/(4-\alpha)}{4-\alpha\over4}2^{-{\alpha/(4-\alpha)}} \biggl(
{2{\mathcal E}(\alpha_0, d, \gamma)\over2-\alpha} \biggr)^{(2-\alpha)/(4-\alpha)},
\end{eqnarray*}
where the equality follows from Lemma 7.2, \cite{CHSX}.

Summarizing our computation since (\ref{lo-15}),
%
%
%e2.14 #&#
\begin{eqnarray}
\label{lo-16}
&& \liminf_{R\to\infty} {1\over mt_R} \log
\E_0{\mathcal Z}_m(R)
\nonumber\\[-8pt]\\[-8pt]\nonumber
&&\qquad \ge(\theta\lambda)^{4/(4-\alpha)}
{4-\alpha\over4}2^{-{\alpha/(4-\alpha)}} \biggl({2{\mathcal E}(\alpha_0, d, \gamma)\over2-\alpha}
\biggr)^{(2-\alpha)/(4-\alpha)}.
\end{eqnarray}
This clearly leads to (\ref{lo-12prime}).
%\end{pf}

%s2.2 #&#
\subsection{The setting of Theorems \texorpdfstring{\protect\ref{th-3}}{1.6} and
\texorpdfstring{\protect\ref{th-4}}{1.7}}\label{lo-s}

Our approach is based on the method of localization developed in
\cite{CJK} and \cite{CJKS}. The construction is specifically designed
for the scheme $(2)\times\mathrm{(I)}$
in \cite{CJKS} and for $(2)\times\mathrm{(III)}$
in \cite{CJK}. This method also works
for $(2)\times\mathrm{(II)}$ with minor modification.
In the following we carry out this approach.

Given $\beta>0$ set
\[
l_\beta(x)=\prod_{j=1}^d
\biggl(1-{\llvert  x_j\rrvert  \over\beta} \biggr)^+\quad\mbox{and} \quad
K_\beta(x)=K(x)l_\beta(x),
\]
where $K(x)$ is given in (\ref{lo-12}). Let $M_\beta(t,A)$ be the
martingale measure constructed through (\ref{intro-35}) with $K(x)$
being replaced by $K_\beta(x)$. For each
$n\ge1$, define
$U^{(\beta,n)}(t,x)$ as the $n$th Picard iteration given in the following
integral equation: $U^{(\beta, 0)}=1$ and
%
%e2.15 #&#
\begin{eqnarray}
U^{(\beta, n+1)}(t,x)=1+\theta\int_0^t\!\!\int
_{[x-\beta\sqrt{t},
x+\beta\sqrt{t}]^d} p_{t-s}(y-x)U^{(\beta, n)}(s,y)M_\beta(ds\,dy)\nonumber
\\
\eqntext{\qquad n=1,2,\ldots.}
\end{eqnarray}
In \cite{CJK} and \cite{CJKS}, the process
$U_\beta(t,x)\equiv U^{(\beta, [\log\beta]+1)}(t,x)$ is used to proximate
$u(t,x)$, where $\beta=\beta(R)$ increases in $R$ with a suitable speed.
By \cite{CJKS}, Lemma~9.8, for any $\{x^{(k)}\}\subset\R^d$ with
$\llvert   x^{(j)}-x^{(k)}\rrvert  \ge2\beta([\log\beta]+1)(1+\sqrt{t})$,
$\{U_\beta(t,x^{(k)})\}$ is an i.i.d. sequence. By \cite{CJKS}, Lemma 9.7,
for every $\eta\in(0, 1\wedge\alpha)$ there are finite and positive
constants $l_i=l_i(d,\alpha, \eta)$ ($i=1,2$) such that uniformly
for $\beta>0$ and $m\ge2$
%
%
%e2.16 #&#
\begin{eqnarray}
\label{lo-17} \E \bigl\llvert u(t,0)-U_\beta(t,0) \bigr\rrvert
^m\le \biggl({l_2m\over\beta^\eta
} \biggr)^{m/2} \exp \bigl\{l_1m^{(4-\alpha)/(2-\alpha)} \bigr\}.
\end{eqnarray}

Let $\eta$ be fixed, and set $\beta=\exp \{M(\log R)^{2/(4-\alpha)} \}$
where $M>0$ is large but fixed (will be specified later). Let
$N=\beta([\log\beta]+1)(1+\sqrt{t})$ and ${\mathcal N}_R=N\Z^d\cap
B(0, R-N)$.
By the fact that $\alpha<2$, $\#({\mathcal N}_R)=CR^{d+o(1)}$ ($R\to
\infty$) where
the constant $C>0$ does not depend on $R$.

Given $\varepsilon>0$
\begin{eqnarray*}
&&\P \Bigl\{\log\max_{z\in{\mathcal N}_R} \bigl\llvert u(t,z)
-U_\beta(t,z) \bigr\rrvert \ge \varepsilon(\log R)^{2/(4-\alpha)} \Bigr\}
\\
&&\qquad \le \#({\mathcal N}_R)\P \bigl\{\log \bigl\llvert
u(t,0)-U_\beta(t,0) \bigr\rrvert \ge \varepsilon(\log R)^{2/(4-\alpha)}
\bigr\}.
\end{eqnarray*}
By Chebyshev's inequality and the moment bound given in (\ref{lo-17}),
\begin{eqnarray*}
&&\P \bigl\{\log \bigl\llvert u(t,0)-U_\beta(t,0) \bigr\rrvert \ge
\varepsilon(\log R)^{2/(4-\alpha)} \bigr\}
\\
&&\qquad \le \exp \{-\varepsilon\log R \} \E \bigl\llvert u(t,0)-U_\beta(t,0)
\bigr\rrvert ^{
(\log R)^{(2-\alpha)/(4-\alpha)}}
\\
&&\qquad \le\exp \{-\varepsilon\log R \} \exp \bigl\{-\tfrac{1}{2} \bigl(\eta
M-l_1-o(1) \bigr)\log R \bigr\}
\end{eqnarray*}
when $R$ is large. Make $M$ sufficiently large, and we have
%
%
%e2.17 #&#
\begin{equation}
\label{lo-18} \P \bigl\{\log \bigl\llvert u(t,0)-U_\beta(t,0) \bigr
\rrvert \ge \varepsilon(\log R)^{2/(4-\alpha)} \bigr\}\le R^{-(d+2)}
\end{equation}
for large $R$. Consequently,
%
%
%e2.18 #&#
\begin{equation}
\label{lo-19} \P \Bigl\{\log\max_{z\in{\mathcal N}_R} \bigl\llvert u(t,z)
-U_\beta(t,z) \bigr\rrvert \ge \varepsilon(\log R)^{2/(4-\alpha)} \Bigr\}
\le R^{-2}
\end{equation}
for large $R$. By the Borel--Cantelli lemma
%
%
%e2.19 #&#
\begin{equation}
\label{lo-19prime}
\qquad\limsup_{n\to\infty} \bigl(\log2^n
\bigr)^{-{2/(4-\alpha)}} \log\max_{z\in{\mathcal N}_{2^n}} \bigl\llvert u(t,z)
-U_{\beta(2^n)}(t,z) \bigr\rrvert =0\qquad\mbox{a.s.}
\end{equation}

On the other hand
\[
\max_{z\in{\mathcal N}_R} \bigl\llvert U_\beta(t,z) \bigr\rrvert
\le\max_{z\in{\mathcal N}_R} u(t,z)+ \max_{z\in{\mathcal N}_R} \bigl
\llvert u(t,z) -U_\beta(t,z) \bigr\rrvert. %
\]
Consequently,
%
%
%e2.20 #&#
\begin{eqnarray}
\label{lo-20}
&& \log\max_{z\in{\mathcal N}_R} \bigl\llvert
U_\beta(t,z) \bigr\rrvert
\nonumber\\[-8pt]\\[-8pt]\nonumber
&&\qquad \le \log2+\max \Bigl\{\log\max
_{z\in{\mathcal N}_R} u(t,z), \log\max_{z\in{\mathcal N}_R}
\bigl\llvert u(t,z) -U_\beta(t,z) \bigr\rrvert \Bigr\}.
\end{eqnarray}

For any $\lambda>0$ satisfying
\[
\lambda+\varepsilon<{4-\alpha\over4} \biggl({4t{\mathcal E}(d,\gamma)\over2-\alpha}
\biggr)^{(2-\alpha)/(4-\alpha)} \theta^{4/(4-\alpha)}d^{2/(4-\alpha)}, %
\]
by independence
\begin{eqnarray*}
&&\P \Bigl\{\log\max_{z\in{\mathcal N}_R} \bigl\llvert U_\beta(t,z)
\bigr\rrvert \le\lambda(\log R)^{2/(4-\alpha)} \Bigr\}
\\
&&\qquad = \bigl(1-\P \bigl\{\log\bigl\llvert U_\beta(t,0)\bigr\rrvert >
\lambda(\log R)^{2/(4-\alpha)} \bigr\} \bigr)^{\#({\mathcal N}_R)}.
\end{eqnarray*}
By (\ref{lo-18})
\begin{eqnarray*}
&&\P \bigl\{\log\bigl\llvert U_\beta(t,0)\bigr\rrvert >\lambda(\log
R)^{2/(4-\alpha)} \bigr\}
\\
&&\qquad \ge \P \bigl\{\log u(t,0) >(\lambda+\varepsilon) (\log R)^{2/(4-\alpha)} \bigr
\}-R^{-(d+2)}.
\end{eqnarray*}
By the large deviation result given in (\ref{up-5}), Theorem \ref{th-6} below,
\[
\P \bigl\{\log u(t,0) >(\lambda+\varepsilon) (\log R)^{2/(4-\alpha)} \bigr\} \ge
\exp \bigl\{-(d-\delta)\log R \bigr\} %
\]
for sufficiently large $R$. Thus we have established the bound
\[
\P \Bigl\{\log\max_{z\in{\mathcal N}_R} \bigl\llvert U_\beta(t,z)
\bigr\rrvert \le\lambda(\log R)^{2/(4-\alpha)} \Bigr\} \le\exp \bigl
\{-R^v \bigr\} %
\]
for some $v>0$. By the Borel--Cantelli lemma,
%
%
%e2.21 #&#
\begin{equation}
\label{lo-21} \liminf_{n\to\infty} \bigl(\log2^n
\bigr)^{-{2/(4-\alpha)}} \log\max_{z\in{\mathcal N}_{2^n}} \bigl\llvert
U_{\beta(2^n)}(t,z) \bigr\rrvert \ge\lambda\qquad\mbox{a.s.}
\end{equation}

Combining (\ref{lo-19prime}), (\ref{lo-20}) and (\ref{lo-21}),
\[
\liminf_{n\to\infty} \bigl(\log2^n \bigr)^{-{2/(4-\alpha)}}
\log\max_{z\in{\mathcal N}_{2^n}} u(t,z)\ge\lambda\qquad\mbox{a.s.} %
\]
By the fact that
\[
\max_{z\in{\mathcal N}_R} u(t,z)\le\max_{\llvert  x\rrvert  \le R} u(t,x)
\]
and by the monotonicity of $ \max_{\llvert  x\rrvert  \le R}
u(t,x)$ in $R$,
\[
\liminf_{R\to\infty}(\log R)^{-{2/(4-\alpha)}}\log\max
_{\llvert  x\rrvert
\le R} u(t,x) \ge\lambda\qquad\mbox{a.s.} %
\]
This leads to the lower bound for (\ref{intro-38}) as $\lambda$ can be made
arbitrarily close to the limit value appearing in the right-hand side of
(\ref{intro-38}).

According to the agreement made at the end of Section~\ref{intro},
the lower bound
requested by (\ref{intro-39}) can be viewed as the special case of the
lower bound for (\ref{intro-38}) under the identification $\alpha=1$ and
$d=1$. %\end{pf}

%s3 #&#
\section{High moment asymptotics}\label{mo}

Associated to the main theorems are the tail behaviors of
$\log u(t,0)$, which are relevant to the high moment asymptotics
for $\E u(t,0)^m$ as $m\to\infty$, in light of the G\"artner--Ellis theorem.
The objective of this section is to find the exact high moment
asymptotics required
by our main theorems.

%s3.1 #&#
\subsection{The setting of Theorems \texorpdfstring{\protect\ref{th-1}}{1.1},
\texorpdfstring{\protect\ref{th-2}}{1.2} and \texorpdfstring{\protect\ref{th-2prime}}{1.3}}\label{mo-f}

Recall our extra assumption $u_0(x)\equiv1$.
We begin with the moment representations (Corollary 4.5 and Remark~4.6,
\cite{CHSX})
%
%
%e3.1 #&#
\begin{equation}
\label{mo-1}
\E u(t,x)^m
=\E_0\exp \Biggl\{
{1\over2}\theta^2 \sum_{j, k=1}^m
\int_0^t\!\!\int_0^t
\gamma_0(r-s)\gamma \bigl(B_j(r)-B_k(s)
\bigr)\,dr\,ds \Biggr\}\hspace*{-20pt}
\end{equation}
for each integer $m\ge1$, where $\{B_k(s)\}_{k\ge1}$ is an i.i.d.
sequence of
Brownian motions.

%
%pr3.1 #&#
\begin{proposition}\label{prop-1}
Under the assumption of Theorem \ref{th-1},
%
%
%e3.2 #&#
\begin{equation}
\label{mo-2} \lim_{m\to\infty}m^{-2}\log\E
u(t,0)^m={1\over2}\theta^2 \gamma(0)\int_0^t\!\!\int_0^t
\gamma_0(r-s)\,dr\,ds.
\end{equation}
\end{proposition}

\begin{pf} We first notice that $\gamma(x)$ reaches its maximum at $x=0$.
Indeed, for a infinitely smooth and rapidly decreasing (at $\infty$) function
$\varphi_0(\cdot)\ge0$ on $\R^+$, $\varepsilon>0$ and $x\in\R^d$,
\begin{eqnarray*}
&&\Cov \bigl(\langle V, \varphi_0 p_\varepsilon\rangle,
 \bigl\langle V, \varphi_0 p_\varepsilon(\cdot-x)
\bigr\rangle \bigr)
\\
&&\qquad = \biggl(\int_{\R^+\times\R^+}\gamma_0(r-s)
\varphi_0(r)\varphi _0(s)\,dr\,ds \biggr)
\\
&&\quad\qquad{}\times  \int
_{\R^d\times\R^d}\gamma(y-z)p_\varepsilon(y)p_\varepsilon(z-x)\,dy\,dz.
\end{eqnarray*}
Here we recall our notation $p_s(x)$ for $d$-dimensional Brownian density.

On the other hand, by homogeneity
\begin{eqnarray*}
&& \Cov \bigl(\langle V, \varphi_0 p_\varepsilon\rangle,
\bigl\langle V, \varphi_0 p_\varepsilon(\cdot-x)
\bigr\rangle \bigr)
\\
&&\qquad  \le\Var \bigl(\langle V, \varphi_0
p_\varepsilon\rangle \bigr)
\\
&&\qquad = \biggl(\int_{\R^+\times\R^+}\gamma_0(r-s)
\varphi_0(r)\varphi _0(s)\,dr\,ds \biggr) \int
_{\R^d\times\R^d}\gamma(y-z)p_\varepsilon(y)p_\varepsilon(z)\,dy\,dz.
\end{eqnarray*}
Consequently
\begin{eqnarray*}
&& \int_{\R^d\times\R^d}\gamma(y-z)p_\varepsilon(y)p_\varepsilon(z-x)\,dy\,dz
\\
&&\qquad
\le\int_{\R^d\times\R^d}\gamma(y-z)p_\varepsilon(y)p_\varepsilon(z)\,dy\,dz.
\end{eqnarray*}
Letting $\varepsilon\to0^+$, by continuity of $\gamma(\cdot)$ we have
$\gamma(x)\le\gamma(0)$.

Therefore, the requested upper bound follows from (\ref{mo-1}).

As for the lower bound, we essentially
follow the strategy used in the previous section:
by
continuity, for any $\varepsilon>0$ there is $\delta>0$ such that
$\gamma(x)\ge\gamma(0)-\varepsilon$ as long as $\llvert  x\rrvert  \le2\delta$.
Thus
\begin{eqnarray*}
\hspace*{-3pt}&&\E_0\exp \Biggl\{{1\over2}\sum
_{j, k=1}^m\int_0^t\!\!\int_0^t \gamma_0(r-s)\gamma
\bigl(B_j(r)-B_k(s) \bigr)\,dr\,ds \Biggr\}
\\
\hspace*{-3pt}&&\qquad \ge\E_0 \Biggl[\exp \Biggl\{{1\over2}\sum
_{j, k=1}^m\int_0^t\!\!\int_0^t \gamma_0(r-s)\gamma
\bigl(B_j(r)-B_k(s) \bigr)\,dr\,ds \Biggr\};\min
_{k\le m} \tau_k(\delta)\ge t \Biggr]
\\
\hspace*{-3pt}&&\qquad \ge\exp \biggl\{{m^2\over2} \bigl(\gamma(0)-\varepsilon \bigr)\int_0^t\!\!\int_0^t
\gamma_0(r-s)\,dr\,ds \biggr\}\P_0 \Bigl\{\min
_{k\le m}\tau_k(\delta)\ge t \Bigr\},
\end{eqnarray*}
where $\tau_k(\delta)$ is the time for $B_k(s)$ to exit from the $\delta$-ball.
Therefore, the requested lower bound follows from the facts that
$\varepsilon>0$ can be arbitrarily small and that the probability
\[
\P_0 \Bigl\{\min_{k\le m}\tau_k(
\delta)\ge t \Bigr\}= \Bigl(\P \Bigl\{\max_{s\le t}\bigl\llvert
B(s)\bigr\rrvert \le\delta \Bigr\} \Bigr)^m %
\]
decays at a speed no faster than exponential rate.
\end{pf}

%
%pr3.2 #&#
\begin{proposition}\label{prop-2}
Under the assumptions of Theorem \ref{th-2}
%
%
%e3.3 #&#
\begin{eqnarray}
\label{mo-3}
&& \lim_{m\to\infty}m^{-{(4-\alpha)/(2-\alpha)}} \log\E
u(t,0)^m
\nonumber\\[-8pt]\\[-8pt]\nonumber
&&\qquad = \biggl({\theta^2\over2} \biggr)^{2/(2-\alpha)}
t^{(4-\alpha-2\alpha_0)/(2-\alpha)}{\mathcal E}(\alpha_0, d, \gamma).
\end{eqnarray}
Under the assumptions of Theorem \ref{th-2prime}
%
%
%e3.4 #&#
\begin{eqnarray}
\label{mo-3prime} \lim_{m\to\infty}m^{-3/2} \log\E
u(t,0)^m={\theta^4\over4} t^{3-2\alpha_0}{\mathcal E}(
\alpha_0, 1, \delta_0).
\end{eqnarray}
\end{proposition}

\begin{pf} We need only to prove (\ref{mo-3}) as (\ref{mo-3prime})
can be viewed as a special case under the identification $d=\alpha=1$.
Recall that $\gamma_0(u)=\llvert  u\rrvert  ^{-\alpha_0}$
in this setting. For any $1\le j, k\le m$, by (\ref{lo-12})
%
%
%e3.5 #&#
\begin{eqnarray}
\label{mo-4} &&\int_0^t\!\!\int
_0^t {\gamma (B_j(r)-B_k(s) )\over\llvert  r-s\rrvert  ^{\alpha_0}}\,dr\,ds\nonumber
\\
&&\qquad =C_0\int_{\R\times\R^d} \biggl[\int
_0^t\llvert u-s\rrvert ^{-{(1+\alpha_0)/2}} K
\bigl(x-B_j(s) \bigr)\,ds \biggr]
\\
&&\hspace*{74pt}{}\times \biggl[\int_0^t
\llvert u-s\rrvert ^{-{(1+\alpha _0)/2}} K \bigl(x-B_k(s) \bigr)\,ds
\biggr]\,du\,dx.\nonumber
\end{eqnarray}
Applying the Cauchy--Schwarz inequality
\begin{eqnarray*}
\hspace*{-2pt}&&\int_0^t\!\!\int_0^t
{\gamma (B_j(r)-B_k(s) )\over\llvert  r-s\rrvert  ^{\alpha_0}}\,dr\,ds
\\
\hspace*{-2pt}&&\qquad \le \biggl(\int_0^t\!\!\int
_0^t {\gamma (B_j(r)-B_j(s) )\over\llvert  r-s\rrvert  ^{\alpha_0}}\,dr\,ds
\biggr)^{1/2}
\biggl(\int_0^t\!\!\int
_0^t {\gamma (B_k(r)-B_k(s) )\over\llvert  r-s\rrvert  ^{\alpha_0}}\,dr\,ds
\biggr)^{1/2}.
\end{eqnarray*}
Consequently,
\begin{eqnarray*}
&&\sum_{j, k=1}^m\int_0^t\!\!\int_0^t {\gamma (B_j(r)-B_k(s) )\over\llvert  r-s\rrvert  ^{\alpha_0}}\,dr\,ds
\\
&&\qquad \le \Biggl\{\sum_{k=1}^m \biggl(\int_0^t\!\!\int_0^t
{\gamma (B_k(r)-B_k(s) )\over\llvert  r-s\rrvert  ^{\alpha_0}}\,dr\,ds \biggr)^{1/2} \Biggr\}^2
\\
&&\qquad \le m\sum_{k=1}^m\int_0^t\!\!\int_0^t
{\gamma (B_k(r)-B_k(s) )\over\llvert  r-s\rrvert  ^{\alpha_0}}\,dr\,ds.
\end{eqnarray*}

Write
\[
X_m(t)=\sum_{j, k=1}^m \int
_0^{t} \int_0^{t}
{\gamma (B_j(r)-B_k(s) )
\over\llvert  r-s\rrvert  ^{\alpha_0}}\,dr\,ds. %
\]
By independence, for any $\beta>0$
\begin{eqnarray*}
\E_0\exp \bigl\{\beta X_m(t) \bigr\} &\le& \biggl(
\E_0\exp \biggl\{m\beta\int_0^t\!\!\int_0^t {\gamma (B(r)-B(s) )\over\llvert  r-s\rrvert  ^{\alpha_0}}\,dr\,ds \biggr\}
\biggr)^m
\\
&=& \biggl(\E_0\exp \biggl\{\beta\int_0^{t_m}\!\!\int_0^{t_m} {\gamma (B(r)-B(s) )\over\llvert  r-s\rrvert  ^{\alpha_0}}\,dr\,ds \biggr
\} \biggr)^m,
\end{eqnarray*}
where $  t_m=tm^{2/(4-\alpha-2\alpha_0)}$, and the equality
follows from the Brownian scaling.
Recall (Theorem 1.1, \cite{CHSX}) that
\begin{eqnarray*}
&& \lim_{m\to\infty}t_m^{-{(4-\alpha-2\alpha_0)/(2-\alpha)}}\log
\E_0\exp \biggl\{\beta\int_0^{t_m}\!\!\int_0^{t_m} {\gamma (B(r)-B(s) )\over\llvert  r-s\rrvert  ^{\alpha_0}}\,dr\,ds \biggr\}
\\
&&\qquad  =\beta^{2/(2-\alpha)}{\mathcal E}(\alpha_0, d, \gamma).
\end{eqnarray*}
We conclude that
%
%
%e3.6 #&#
\begin{eqnarray}
\label{mo-5}
&& \limsup_{m\to\infty}m^{-{(4-\alpha)/(2-\alpha)}} \log
\E_0\exp \bigl\{\beta X_m(t) \bigr\}
\nonumber\\[-8pt]\\[-8pt]\nonumber
&&\qquad \le
\beta^{2/(2-\alpha)} t^{(4-\alpha-2\alpha_0)/(2-\alpha)}{\mathcal E}(\alpha_0, d,
\gamma)
\end{eqnarray}
for any $\beta>0$.

On the other hand, let $\tilde{t}_m=tm^{2/(2-\alpha)}$.
An obvious modification of the argument for (\ref{lo-16}) [with
$ (\log R)^{2/(2-\alpha)}$ being replaced by $\tilde{t}_m$]
shows that for any $\lambda>0$,
\begin{eqnarray*}
&& \liminf_{m\to\infty} {1\over m\tilde{t}_m} \log
\E_0\exp \bigl\{\lambda\tilde{t}_m^{\alpha_0/2}X_m(
\tilde {t}_m)^{1/2} \bigr\}
\\
&&\qquad \ge\lambda^{4/(4-\alpha)}
{4-\alpha\over4}2^{-{\alpha/(4-\alpha)}} \biggl({2{\mathcal E}(\alpha_0, d, \gamma)\over2-\alpha}
\biggr)^{(2-\alpha)/(4-\alpha)}.
\end{eqnarray*}

Let $\lambda=\beta t^{-\alpha_0/2}$.
By Brownian scaling, the above limiting bound can be re-written as
%
%
%e3.7 #&#
%e3.8 #&#
\begin{eqnarray}
\label{mo-6}
\qquad &&\liminf_{m\to\infty}m^{-{(4-\alpha)/(2-\alpha)}} \log
\E_0\exp \bigl\{\beta m^{(4-\alpha)/(2(2-\alpha))}X_m(t)^{1/2}
\bigr\}\nonumber
\\
&&\qquad \ge t^{(4-\alpha-2\alpha_0)/(4-\alpha)}\beta^{4/(4-\alpha)} {4-\alpha\over4}2^{-{\alpha/(4-\alpha)}}
\biggl({2{\mathcal E}(\alpha_0, d, \gamma)\over2-\alpha} \biggr)^{(2-\alpha)/(4-\alpha)}
\\
\eqntext{(\beta>0).}
\end{eqnarray}
By the first half of Lemma \ref{A-1} in the
\hyperref[appendix]{Appendix}, (\ref{mo-5})
and (\ref{mo-6}),
\[
\lim_{m\to\infty}m^{-{(4-\alpha)/(2-\alpha)}} \log\E_0\exp \bigl\{
\beta X_m(t) \bigr\}= \beta^{2/(2-\alpha)} t^{(4-\alpha-2\alpha_0)/(2-\alpha)}{\mathcal
E}(\alpha_0, d, \gamma). %
\]
Let $\beta=\theta^2/2$.
Proposition \ref{prop-2} follows from (\ref{mo-1}).
\end{pf}

%s3.2 #&#
\subsection{The setting of Theorems \texorpdfstring{\protect\ref{th-3}}{1.6}~and~\texorpdfstring{\protect\ref{th-4}}{1.7}}\label{mo-s}

The goal here is to establish the following:

%
%pr3.3 #&#
\begin{proposition}\label{prop-3} Under the assumptions given in
Theorem \ref{th-3},
%
%
%e3.9 #&#
\begin{eqnarray}
\label{mo-7} \lim_{m\to\infty}m^{-{(4-\alpha)/(2-\alpha)}} \log\E
u(t,0)^m=t \biggl({\theta^2\over2} \biggr)^{2/(2-\alpha)} {
\mathcal E}(d, \gamma).
\end{eqnarray}
Under the assumptions given in Theorem \ref{th-4},
%
%
%e3.10 #&#
\begin{equation}
\label{mo-8} \lim_{m\to\infty}m^{-3} \log\E
u(t,0)^m=t{\theta^4\over24}.
\end{equation}
\end{proposition}

In view of (\ref{intro-28}), we need only to
prove (\ref{mo-7}), as (\ref{mo-8}) can be viewed as a special case
under a proper identification.
Our starting point is the following moment representation (see Theorem 5.3
in \cite{HN}
and Theorem 3.1 in \cite{Conus}):
%
%
%e3.11 #&#
\begin{eqnarray}
\label{mo-9} \E u(t,0)^m=\E_0\exp \biggl\{
\theta^2\sum_{1\le j<k\le m}\int
_0^t \gamma \bigl(B_j(s)-B_k(s)
\bigr)\,ds \biggr\}.
\end{eqnarray}

The approach here is much more delicate
due to the absence of the diagonal terms in the
$(j, k)$-summation in (\ref{mo-9}) and the fact that
the missing diagonal terms blow up.
The proof consists of several steps.

Let $t_m=tm^{2/(2-\alpha)}$ and $\varepsilon>0$ be small but fixed. Set
\[
\gamma_\varepsilon(x)=\int_{\R^d}p_{2\varepsilon}(x-y)
\gamma(y)\,dy,\qquad x\in\R^d. %
\]
Here we recall that $p_t(x)$ represents the density function of a
$d$-dimensional Brownian motion $B(t)$ starting at 0.
Let the kernel $K(x)$ be defined
in (\ref{lo-12}).
Clearly,
%
%
%e3.12 #&#
\begin{equation}
\label{mo-11} \gamma_\varepsilon(x)=\int_{\R^d}K_\varepsilon(y)K_\varepsilon(y-x)\,dy,
\qquad x\in\R^d,
\end{equation}
where
%
%
%e3.13 #&#
\begin{equation}
\label{mo-11prime} K_\varepsilon(x)=\int_{\R^d}p_{\varepsilon}(x-y)K(y)\,dy,
\qquad x\in\R^d.
\end{equation}

Our first step is to
prove the following:

%
%le3.4 #&#
\begin{lemma}\label{le-1}
For any $\beta>0$
%
%
%e3.14 #&#
\begin{eqnarray}
\label{mo-12}
&& \lim_{m\to\infty}m^{-{(4-\alpha)/(2-\alpha)}}\nonumber
\\
&&\quad{}\times \log \E_0\exp \Biggl\{\beta \Biggl(t_m \int_0^{t_m}\!\!\int_{\R^d} \Biggl[
\sum_{j=1}^m K_\varepsilon
\bigl(x-B_j(s) \bigr) \Biggr]^2\,dx\,ds
\Biggr)^{1/2} \Biggr\}
\\
&&\qquad  =tM_\varepsilon(\beta),\nonumber
\end{eqnarray}
where
\begin{eqnarray*}
M_\varepsilon(\beta) &=& \sup_{g\in{\mathcal A}_d} \biggl\{\beta \biggl(
\int_0^1\!\!\int_{\R^d}
\biggl[\int_{\R^d}K_\varepsilon(y-x)g^2(s, y)\,dy
\biggr]^2 \,dx\,ds \biggr)^{1/2}
\\
&&\hspace*{110pt}{} - {1\over2}\int_0^1\!\!\int_{\R^d}\bigl
\llvert \nabla_x g(s,x)\bigr\rrvert ^2\,dx\,dy \biggr\}.
\end{eqnarray*}
\end{lemma}

\begin{pf} Indeed,
\begin{eqnarray*}
&& \Biggl(t_m \int_0^{t_m}\!\!\int
_{\R^d} \Biggl[\sum_{j=1}^m
K_\varepsilon \bigl(x-B_j(s) \bigr) \Biggr]^2\,dx\,ds
\Biggr)^{1/2}
\\[-1pt]
&&\qquad =t_m \Biggl( \int_0^{1}\!\!\int
_{\R^d} \Biggl[\sum_{j=1}^m
K_\varepsilon \bigl(x-B_j(t_ms) \bigr)
\Biggr]^2\,dx\,ds \Biggr)^{1/2}
\\[-1pt]
&&\qquad \ge t_m\int_0^{1}\!\!\int
_{\R^d}f(s,x) \Biggl[\sum_{j=1}^m
K_\varepsilon \bigl(x-B_j(t_ms) \bigr) \Biggr]\,dx\,ds
\\[-1pt]
&&\qquad =\sum_{j=1}^m\int_0^{t_m}
\bar{f} \biggl({s\over t_m}, B_j(s) \biggr)\,ds,
\end{eqnarray*}
where $f(s,x)\ge0$ is a compactly supported and continuous
function on $[0,1]\times\R^d$ with
\begin{eqnarray*}
\int_0^{1}\!\!\int_{\R^d}f^2(s,x)\,dx\,ds&=&1,
\\[-1pt]
\bar{f}(s, x)&=&\int_{\R^d}f(s,y)K_\varepsilon(y-x)\,dy,\qquad x\in\R^d,
\end{eqnarray*}
and the second step follows from the Cauchy--Schwarz inequality.
By independence,
\begin{eqnarray*}
&&\E_0\exp \Biggl\{\beta \Biggl(t_m \int_0^{t_m}\!\!\int_{\R^d} \Biggl[
\sum_{j=1}^m K_\varepsilon
\bigl(x-B_j(s) \bigr) \Biggr]^2\,dx\,ds
\Biggr)^{1/2} \Biggr\}
\\
&&\qquad \ge \biggl(\E_0\exp \biggl\{\beta\int_0^{t_m}
\bar{f} \biggl({s\over t_m}, B(s) \biggr)\,ds \biggr\}
\biggr)^m.
\end{eqnarray*}
Applying Proposition 3.1, \cite{CHSX} or (\ref{lo-14}) to the
right-hand side,
\begin{eqnarray*}
&&\lim_{m\to\infty}m^{-{(4-\alpha)/(2-\alpha)}}\log\E_0\exp \Biggl
\{\beta \Biggl(t_m \int_0^{t_m}\!\!\int
_{\R^d} \Biggl[\sum_{j=1}^m
K_\varepsilon \bigl(x-B_j(s) \bigr) \Biggr]^2\,dx\,ds
\Biggr)^{1/2} \Biggr\}
\\
&&\qquad \ge t\sup_{g\in{\mathcal A}_d} \biggl\{\beta \int_0^1\!\!\int_{\R^d} \bar{f}(s,x)g^2(s,x)\,dx\,ds -
{1\over2}\int_0^1\!\!\int
_{\R^d}\bigl\llvert \nabla_x g(s, x) \bigr\rrvert
^2\,dx\,ds \biggr\}
\\
&&\qquad =t\sup_{g\in{\mathcal A}_d} \biggl\{\beta \int_0^1\!\!\int_{\R^d} f(s,y) \biggl[\int_{\R^d}K_\varepsilon
(y-x)g^2(s,x)\,dx \biggr]\,dy\,ds
\\
&&\hspace*{155pt}{} -{1\over2}\int_0^1\!\!\int_{\R^d}\bigl
\llvert \nabla_x g(s, x) \bigr\rrvert ^2\,dx\,ds \biggr\}.
\end{eqnarray*}
Taking supremum over $f$ on the right-hand side leads to
the lower bound requested by (\ref{mo-12}).

The proof of the upper bound is harder. First, we perform the following
smooth truncation: let $l$: $\R^+\longrightarrow[0, 1]$ be a smooth function
satisfying the following properties:
$l(u)=1$ for $u\in[0,1]$, $l(u)=0$ for $u\ge3$ and $-1\le l'(u)\le0$
for all $u>0$. Let $M>0$ be a large number, and write
\[
Q(x)=K_\varepsilon(x) l \bigl(M^{-1}\llvert x\rrvert \bigr).
\]
One can easily see that $Q(x)$ is supported on
$B(0, 3M)=\{x\in\R^d; \llvert  x\rrvert  \le3M\}$ and that
\[
\int_{\R^d} \bigl[K_\varepsilon(x)-Q(x)
\bigr]^2\,dx\longrightarrow0\qquad (M\to\infty). %
\]
By the triangle inequality,
\begin{eqnarray*}
&& \Biggl(t_m \int_0^{t_m}\!\!\int
_{\R^d} \Biggl[\sum_{j=1}^m
\bigl(K_\varepsilon \bigl(x-B_j(s) \bigr)-Q
\bigl(x-B_j(s) \bigr) \bigr) \Biggr]^2\,dx\,ds
\Biggr)^{1/2}
\\
&&\qquad \le t_m^{1/2}\sum_{j=1}^m
\biggl(\int_0^{t_m}\!\!\int_{\R^d}
\bigl[K_\varepsilon \bigl(x-B_j(s) \bigr)-Q
\bigl(x-B_j(s) \bigr) \bigr]^2\,dx\,ds \biggr)^{1/2}
\\
&&\qquad =m t_m \biggl(\int_{\R^d} \bigl[K_\varepsilon(x)-Q(x)
\bigr]^2\,dx \biggr)^{1/2}
\\
&&\qquad =tm^{(4-\alpha)/(2-\alpha)} \biggl(\int
_{\R^d} \bigl[K_\varepsilon(x)-Q(x) \bigr]^2\,dx
\biggr)^{1/2}.
\end{eqnarray*}
This estimate shows that it suffices to establish the upper bound with
$K_\varepsilon(x)$ being replaced by $Q(x)$ for an arbitrarily large $M$.

Let $M>0$ be fixed. For $N>3M$, we have
\begin{eqnarray*}
&& \int_{\R^d} \Biggl[\sum_{j=1}^m
Q \bigl(x-B_j(s) \bigr) \Biggr]^2\,dx
\\
&&\qquad = \sum
_{z\in\Z^d}\int_{[-N, N]^d} \Biggl[\sum
_{j=1}^m Q \bigl(2zN+x-B_j(s) \bigr)
\Biggr]^2\,dx
\\
&&\qquad \le \int_{[-N, N]^d} \Biggl[\sum_{j=1}^m
Q_N \bigl(x-B_j(s) \bigr) \Biggr]^2\,dx
\\
&&\qquad =
m^2\int_{[-N, N]^d}\eta^2_m(s,
x)\,dx,
\end{eqnarray*}
where
%
%
%e3.15 #&#
\begin{eqnarray}\label{mo-12prime}
Q_N(x)&=&\sum_{z\in\Z^d}Q(2zN+x)
\quad\mbox{and}
\nonumber\\[-8pt]\\[-8pt]\nonumber
\eta_m(s, x)&=&{1\over m}
\sum_{j=1}^mQ_N
\bigl(x-B_j(s) \bigr).
\end{eqnarray}
Notice in the $z$-summation that defines $Q_N(\cdot)$, there is at most
one nonzero term for any $x\in\R^d$ by the assumption
that $N>3M$. Consequently, $Q_N(x)$ is a continuous periodic
extension (with the period $2N$) of $Q(x)$.

Further, by integration substitution
\[
\int_0^{t_m}\!\!\int_{[-N, N]^d}
\eta_m^2(s, x)\,dx\,ds=t_m\int_0^1\!\!\int_{[-N, N]^d}
\eta_m^2(t_ms, x)\,dx\,ds. %
\]

To establish the upper bound requested by (\ref{mo-12}), therefore, all we
need is to show that for any $M>0$,
%
%
%e3.16 #&#
\begin{eqnarray}
\label{mo-13}
&& \limsup_{N\to\infty}\limsup_{m\to\infty}
m^{-{(4-\alpha)/(2-\alpha)}}\nonumber
\nonumber\\[-8pt]\\[-8pt]\nonumber
&&\qquad{}\times \log\E_0\exp \biggl\{\beta mt_m \biggl(
\int_0^1\!\!\int_{\R^d}
\eta^2_m (t_ms,x)\,dx\,ds \biggr)^{1/2}
\biggr\}
\le tM_\varepsilon(\beta).\nonumber
\end{eqnarray}

We\vspace*{1pt} let $N>3M$ be fixed for a while and concentrate on the $m$-$\limsup$.
Unfortunately, $\eta_m(t_m(\cdot), \cdot)$ is not\vspace*{1pt} exponentially tight
when embedded into ${\mathcal L}^2([0,1]\times[-N, N]^d)$. In the following
we prove that with overwhelming
probability, for any $\delta>0$ there is a $C>0$ such that the range of
the ${\mathcal L}^2([0,1]\times[-N, N]^d)$-valued random variable
$\eta_m(t_m(\cdot), \cdot)$ is covered
by at most $\exp(Ct_m)$ $\delta$-balls in ${\mathcal L}^2([0,1]\times
[-N, N]^d)$.

Let $\upsilon>0$ be a small number, and define
$[s]_\upsilon=\upsilon[\upsilon^{-1}s]$.
By Jensen inequality,
\begin{eqnarray*}
&&\int_0^{t_m}\!\!\int_{[-N, N]^d}
\bigl[\eta_m(s,x)-\eta_m \bigl([s]_\upsilon, x
\bigr) \bigr]^2 \,dx\,ds
\\
&&\qquad \le{1\over m}\sum_{j=1}^m
\int_0^{t_m}\!\!\int_{[-N, N]^d}
\bigl[Q_N \bigl(x-B_j(s) \bigr)-Q_N
\bigl(x-B_j \bigl([s]_\upsilon \bigr) \bigr)
\bigr]^2\,dx\,ds.
\end{eqnarray*}
By independence,
\begin{eqnarray*}
\hspace*{-2pt}&&\E_0\exp \biggl\{\beta m\int_0^{t_m}\!\!\int_{[-N, N]^d} \bigl[\eta_m(s,x)- \eta_m
\bigl([s]_\upsilon, x \bigr) \bigr]^2 \,dx\,ds \biggr\}
\\
\hspace*{-2pt}&&\qquad \le \biggl(\E_0\exp \biggl\{\beta\int_0^{t_m}\!\!\int_{[-N, N]^d} \bigl[Q_N \bigl(x-B(s)
\bigr)-Q_N \bigl(x-B \bigl([s]_\upsilon \bigr) \bigr)
\bigr]^2\,dx\,ds \biggr\} \biggr)^m.
\end{eqnarray*}
Notice that
\begin{eqnarray*}
&&\int_0^{t_m}\!\!\int_{[-N, N]^d}
\bigl[Q_N \bigl(x-B(s) \bigr)-Q_N \bigl(x-B
\bigl([s]_\upsilon \bigr) \bigr) \bigr]^2\,dx\,ds
\\
&&\qquad \le\sum_{k}\int_{(k-1)\upsilon}^{k\upsilon}
\int_{[-N, N]^d} \bigl[Q_N \bigl(x-B(s) \bigr)-
Q_N \bigl(x-B \bigl((k-1)\upsilon \bigr) \bigr)
\bigr]^2\,dx\,ds
\\
&&\qquad =\sum_{k}\int_{(k-1)\upsilon}^{k\upsilon}
\int_{[-N, N]^d} \bigl[Q_N (x )-Q_N
\bigl(x+B(s)-B \bigl((k-1)\upsilon \bigr) \bigr) \bigr]^2\,dx\,ds,
\end{eqnarray*}
where the summation over $k$ runs from $k=1$ until $k=[\upsilon^{-1} t_m]+1$,
and the second step follows from the periodicity of the function
$Q_N(\cdot)$.
By increment-independence of the Brownian motion,
\begin{eqnarray*}
&&\E_0\exp \biggl\{\beta m\int_0^{t_m}\!\!\int_{[-N, N]^d} \bigl[\eta_m(s,x)- \eta_m
\bigl([s]_\upsilon, x \bigr) \bigr]^2 \,dx\,ds \biggr\}
\\
&&\qquad \le \biggl(\E_0\exp \biggl\{\beta\int_0^v\!\int_{[-N, N]^d} \bigl[Q_N (x )-Q_N
\bigl(x+B(s) \bigr) \bigr]^2\,dx\,ds \biggr\} \biggr)^{m([\upsilon^{-1} t_m]+1)}.
\end{eqnarray*}
By the continuity of $Q_N(\cdot)$ one can easily see that
%
%e3.17 #&#
\begin{eqnarray}
\E_0\exp \biggl\{\beta\int_0^v\!\int_{[-N, N]^d} \bigl[Q_N (x )-Q_N
\bigl(x+B(s) \bigr) \bigr]^2\,dx\,ds \biggr\} =\exp \bigl\{o(\upsilon)
\bigr\}\nonumber
\\
\eqntext{\bigl(\upsilon\to0^+ \bigr).}
\end{eqnarray}
Thus we have proved that for any $\beta>0$,
\begin{eqnarray*}
&& \lim_{\upsilon\to0^+}\limsup_{m\to\infty} m^{-{(4-\alpha)/(2-\alpha)}}
\\
&&\qquad{}\times \log\E_0\exp \biggl\{\beta m\int_0^{t_m}\!\!\int_{[-N, N]^d} \bigl[\eta_m(s,x)- \eta_m
\bigl([s]_\upsilon, x \bigr) \bigr]^2 \,dx\,ds \biggr\}=0.
\end{eqnarray*}
Write
\[
\Omega(\upsilon, \delta, m) = \biggl\{\int_0^{1}\!\!\int_{[-N, N]^d} \bigl[\eta_m(t_ms,x)-
\eta_m \bigl([t_ms]_\upsilon, x \bigr)
\bigr]^2 \,dx\,ds\le{\delta^2\over4} \biggr\}. %
\]
By variable substitution and
Chebyshev's inequality, for any $L>0$ one can take $\upsilon$
sufficiently small
so that
\[
\P_0 \bigl(\Omega(\upsilon, \delta, m)^c \bigr) \le
\exp \bigl\{-L m^{(4-\alpha)/(2-\alpha)} \bigr\} %
\]
for large $m$. Define
\[
\tau_j(H)=\inf \bigl\{s\ge0; \bigl\llvert
B_j(s) \bigr\rrvert \ge Hm^{(6-\alpha)/(2(2-\alpha))} \bigr\}
\quad\mbox{and}\quad \tau_{*}(H)=\min_{j\le m}
\tau_j(H),
\]
where $H>0$ is a large but fixed constant. By Gaussian tail,
\begin{eqnarray*}
\P_0 \bigl\{\tau_{*}(H)< t_m \bigr\}&\le& m
\P_0 \Bigl\{\max_{s\le t_m}\bigl\llvert B(s)\bigr
\rrvert \ge Hm^{(6-\alpha)/(2(2-\alpha))} \Bigr\}
\\
&\le& m\exp \bigl\{-CH^2
m^{(4-\alpha)/(2-\alpha)} \bigr\},
\end{eqnarray*}
where $C>0$ is a universal constant.

By the fact that $\eta_m (t_ms, x)$ is bounded by
a deterministic constant $C_N$ independent of $m$, for sufficiently small
$\upsilon$ and sufficiently large $H>0$,
%
%
%e3.18 #&#
\begin{eqnarray}
\label{mo-14} &&\E_0\exp \biggl\{\beta mt_m \biggl(
\int_0^1\!\!\int_{\R^d}
\eta^2_m (t_ms,x)\,dx\,ds \biggr)^{1/2}
\biggr\}\nonumber
\\
&&\qquad = \biggl[\E_0\exp \biggl\{\beta mt_m \biggl( \int_0^1\!\!\int_{\R^d}
\eta^2_m (t_ms,x)\,dx\,ds \biggr)^{1/2}
\biggr\} 1_{\Omega(\upsilon, \delta, m)};\tau_{*}\ge t_m
\biggr]
\nonumber\\[-8pt]\\[-8pt]\nonumber
&&\quad\qquad{} +\exp \bigl\{-(L-\beta C_N)m^{(4-\alpha)/(2-\alpha)} \bigr\}
\\
&&\quad\qquad{} + m\exp \bigl
\{- \bigl(H^2C-\beta C_N \bigr)m^{(4-\alpha)/(2-\alpha)} \bigr\}.
\nonumber
\end{eqnarray}
The second and the third terms on the right-hand side are negligible for
sufficiently large $L$ and
$H$.

We view $\eta_m(s, \cdot)$ ($s\ge0$) as a process (in $s$)
taking values in ${\mathcal L}^2([-N, N]^d)$.
Notice that the function $Q(x)$ is bounded and
Lipschitz continuous. These properties are inherited
by $Q_N(x)$ as a continuous periodic extension of $Q(x)$. Consequently,
$Q_N (\cdot-B_j(s) )$ is
bounded and Lipschitz continuous on $[-N, N]^d$ uniformly in $s\ge0$
and $j\ge1$
with a deterministic bound and a deterministic Lipschitz constant.
Hence there is a
deterministic and convex compact set ${\mathcal K}\subset{\mathcal
L}^2([-N, N]^d)$
such that $Q_N (\cdot-B_j(s) )\in{\mathcal K}$ a.s.
for every $s\ge0$ and $j=1,2,\ldots.$ As a convex linear combination
of $Q_N (\cdot-B_j(s) )$ ($j=1,\ldots, m$),
$\eta_m(s, \cdot)\in{\mathcal K}$ a.s. for any $s\ge0$ and
$m=1,2,\ldots.$ Let
$g_1,\ldots, g_l\in{\mathcal K}$ be a $(2^{-1}\delta)$-net of
${\mathcal K}$.
On the set $\Omega(\upsilon, \delta, m)$ the functions of the form
\[
g(s, x)=g_{i_k}(x)\qquad\mbox{as } s\in \biggl[{(k-1)\upsilon\over t_m},
{k\upsilon\over t_m} \biggr),  k=1, 2,\ldots, \bigl[
\upsilon^{-1}t_m \bigr]+1 %
\]
make a $\delta$-net (denoted as ${\mathcal N}_m^\delta$)
of the range of the ${\mathcal L}^2([0,1]\times[-N, N]^d)$-valued random
variable $\eta_m(t_m(\cdot), \cdot)$.
Indeed,
for any $k\ge1$ there is $g_{i_k}(x)=g_{i_k}(\omega, x)$ out of $\{
g_1,\ldots, g_l\}$ such that
\[
\int_{[-N, N]^d} \biggl\llvert \eta_m \biggl(
{(k-1)v\over t_m}, x \biggr) -g_{i_k}(x) \biggr\rrvert
^2 \,dx<{\delta^2\over4}. %
\]
Here the notation $g_{i_k}(\omega, x)$ indicates the randomness of picking
$g_{i_k}$.

Consequently,
\begin{eqnarray*}
&&\int_0^1\!\!\int_{[-N, N]^d}
\bigl\llvert \eta_m \bigl([t_ms]_\upsilon, x
\bigr)-g(s,x) \bigr\rrvert ^2 \,dx\,ds
\\
&&\qquad \le{\upsilon\over t_m}\sum_k \int
_{[-N, N]^d} \biggl\llvert \eta_m \biggl(
{(k-1)v\over t_m},  x \biggr) -g_{i_k}(x) \biggr
\rrvert ^2 \,dx< {\delta^2\over4}.
\end{eqnarray*}
So our assertion follows from the restriction by the set
$\Omega(\upsilon, \delta, m)$.

In addition, we can see that
$\#({\mathcal N}_m^\delta)\le l^{[\upsilon^{-1}t_m]+1}$.
Further, by
our construction of $g\in{\mathcal N}_m^\delta$,
%
%
%e3.19 #&#
\begin{eqnarray}
\label{mo-15} &&\int_0^1\!\!\int
_{[-N, N]^d} \bigl\llvert g(s,x)\bigr\rrvert ^2\,dx\,ds\nonumber
\\
&&\qquad  \le
{\upsilon\over t_m}\sum_k\int
_{[-N, N]^d}\bigl\llvert g_{i_k}(x)\bigr\rrvert
^2\,dx
\\
&&\qquad \le 2\sup_{h\in{\mathcal K}}\int_{[-N, N]^d}\bigl\llvert
h(x) \bigr\rrvert ^2\,dx<\infty, \qquad g\in{\mathcal
N}_m^\delta,\nonumber
\end{eqnarray}
and similarly, for any $0<u<1$,
%
%
%e3.20 #&#
\begin{eqnarray}
\label{mo-16}
\qquad\quad \int_0^u\!\!\int
_{[-N, N]^d} \bigl\llvert g(s,x)\bigr\rrvert ^2\,dx\,ds \le
u\sup_{h\in{\mathcal K}}\int_{[-N, N]^d}\bigl\llvert h(x)
\bigr\rrvert ^2\,dx,\qquad  g\in{\mathcal N}_m^\delta.
\end{eqnarray}
We emphasize the fact that the bounds in (\ref{mo-15}) and
(\ref{mo-16}) do not depends on $\delta$.

By the Hahn--Banach theorem, for each $g\in{\mathcal N}_m^\delta$
there is
$f\in{\mathcal L}^2([0,1]\times[-N, N]^d)$ such that
\begin{eqnarray*}
\int_0^1\!\!\int_{[-N, N]^d}
\bigl\llvert f(s,x)\bigr\rrvert ^2\,dx\,ds &=& 1,
\\
\int_0^1\!\!\int_{[-N, N]^d}
f(s,x)g(s,x)\,dx\,ds &=& \biggl(\int_0^1\!\!\int
_{[-N, N]^d} \bigl\llvert g(s,x)\bigr\rrvert ^2\,dx\,ds
\biggr)^{1/2}. %
\end{eqnarray*}
In view of the uniform bound (\ref{mo-15}) on $g\in{\mathcal
N}_m^\delta$,
for any given
$\sigma>0$ one can
take $\delta>0$ sufficiently small so that
\begin{eqnarray*}
&& \int_0^1\!\!\int_{[-N, N]^d}
f(s,x)h(s,x)\,dx\,ds
\\
&&\qquad >(1-\sigma) \biggl(\int_0^1\!\!\int_{[-N, N]^d} \bigl\llvert h(s,x)\bigr\rrvert ^2\,dx\,ds
\biggr)^{1/2} %
\end{eqnarray*}
for every $h\in B(g,\delta)$ and $g\in{\mathcal N}_m^\delta$. By bound
(\ref{mo-16}), we may make $u>0$ sufficiently small
(but independent of $f$) so that $\llvert   f(s,x)\rrvert  \le1$ for $0\le
s\le u$
and $x\in[-N, N]^d$, due to the fact that one can change the definition
of $f(s,x)$ on $[0,u]\times[-N, N]^d$ without drastically changing
the value of the integral on $[0,1]\times[-N, N]^d$.
Finally, we may make each $f(s,x)$ continuous and
bounded on $[0,1]\times[-N, N]^d$ by (\ref{mo-15})
and the fact that these kinds of functions are dense in
${\mathcal L}^2([0,1]\times[-N, N]^d)$.
Denote
the collection of such $f$ by $ ({\mathcal N}_m^\delta )^*$.
Our way\vspace*{1pt} of using the Hahn--Banach theorem defines a surjective map from
${\mathcal N}_m^\delta$ to $ ({\mathcal N}_m^\delta )^*$. Consequently,
$\# ( ({\mathcal N}_m^\delta )^* )\le\# ({\mathcal
N}_m^\delta )
=l^{[v^{-1}t_m]+1}\le\exp(Ct_m)$, where the constant $C>0$ is
independent of $m$
(but dependent on $l$ and $v$).

On the set $\Omega(\upsilon, \delta, m)$, in particular,
\begin{eqnarray*}
&& \biggl(\int_0^1\!\!\int_{[-N, N]^d}
\eta_m^2(t_ms, x)\,dx\,ds \biggr)^{1/2}
\\
&&\qquad \le(1-\sigma)^{-1} \max_{f\in ({\mathcal N}_m^\delta )^*} \int_0^1\!\!\int_{[-N, N]^d}f(s,x)
\eta_m(t_ms, x)\,dx\,ds.
\end{eqnarray*}
Therefore,
\begin{eqnarray*}
\hspace*{-3pt}&&\E_0 \biggl[\exp \biggl\{\beta mt_m \biggl(\int_0^1\!\!\int_{[-N, N]^d}
\eta_m^2(t_ms, x)\,dx\,ds \biggr)^{1/2}
\biggr\}1_{\Omega(\upsilon, \delta, m)};\tau_{*}(H)\ge t_m
\biggr]
\\
\hspace*{-3pt}&&\qquad \le\exp(Ct_m)\max_{f\in ({\mathcal N}_m^\delta )^*} \E_0
\biggl[\exp \biggl\{{\beta mt_m\over1-\sigma} \int_0^1\!\!\int_{[-N, N]^d}f(s,x)\eta_m(t_ms, x)\,dx\,ds
\biggr\};
\\
\hspace*{-3pt}&&\hspace*{299pt} \tau_{*}(H)\ge t_m \biggr].
\end{eqnarray*}
Notice that
%
%
%e3.21 #&#
\begin{eqnarray}
\label{mo-17} &&\int_0^1\!\!\int
_{[-N, N]^d}f(s,x)\eta_m(t_ms, x)\,dx\,ds\nonumber
\\
&&\qquad ={1\over m}\sum_{j=1}^m
\int_0^1 \biggl[\int_{[-N, N]^d}
f(s, x)Q_N \bigl(x-B_j(t_ms) \bigr) \,dx
\biggr]\,ds
\\
&&\qquad ={1\over mt_m}\sum_{j=1}^m
\int_0^{t_m}\tilde{f} \biggl(
{s\over t_m}, B_j(s) \biggr) \,ds,\nonumber
\end{eqnarray}
where
\[
\tilde{f}(s,x)=\int_{[-N, N]^d}f(s, y)Q_N(y-x)\,dy.
\]

Summarizing our argument since (\ref{mo-14}), we conclude that
%
%
%e3.22 #&#
\begin{eqnarray}
\label{mo-18} \qquad&&\limsup_{m\to\infty} m^{-{(4-\alpha)/(2-\alpha)}}\log
\E_0\exp \biggl\{\beta mt_m \biggl( \int_0^1\!\!\int_{\R^d}
\eta^2_m (t_ms,x)\,dx\,ds \biggr)^{1/2}
\biggr\}\nonumber
\\
&&\qquad \le t\limsup_{m\to\infty}{1\over t_m}\log\max
_{f\in ({\mathcal
N}_m^\delta )^*} \E_0 \biggl[\exp \biggl\{
{\beta\over1-\sigma} \int_0^{t_m}\tilde{f}
\biggl({s\over t_m}, B(s) \biggr)\,ds \biggr\};
\\
&&\hspace*{261pt} \tau(H)\ge
t_m \biggr]. \nonumber
\end{eqnarray}
Here we recall our notation
\[
\tau(H)=\inf \bigl\{s\ge0; \bigl\llvert B(s)\bigr\rrvert \ge
Hm^{(6-\alpha)/(2(2-\alpha))} \bigr\}. %
\]

Let $f\in ({\mathcal N}_m^\delta )^*$. For large $m$
\begin{eqnarray*}
\hspace*{-3pt}&&\E_0 \biggl[\exp \biggl\{{\beta\over1-\sigma} \int
_0^{t_m}\tilde{f} \biggl({s\over t_m},
B(s) \biggr)\,ds \biggr\}; \tau(H)\ge t_m \biggr]
\\
\hspace*{-3pt}&&\qquad \le\exp \biggl\{{\beta\over1-\sigma} \biggr\} \E_0 \biggl[
\exp \biggl\{{\beta\over1-\sigma} \int_1^{t_m}
\tilde{f} \biggl({s\over t_m}, B(s) \biggr)\,ds \biggr\};
\tau(H)\ge t_m \biggr]
\\
\hspace*{-3pt}&&\qquad = \exp \biggl\{{\beta\over1-\sigma} \biggr\}
\\
\hspace*{-3pt}&&\quad\qquad{}\times  \int_{B(0, Hm^{(6-\alpha)/(2(2-\alpha))})}
p_1(x)\E_x \biggl[ \exp \biggl\{ {\beta\over1-\sigma}
\\
\hspace*{-3pt}&&\hspace*{196pt}{}\times \int_0^{t_m-1} \tilde{f} \biggl(
{s+1\over t_m}, B(s) \biggr)\,ds \biggr\};
\\
\hspace*{-3pt}&&\hspace*{287pt} \tau(H)\ge t_m \biggr]\,dx,
\end{eqnarray*}
where $p_1(x)$ is the density function of $B(1)$ and the second step
follows from Markov's property. Here and elsewhere, we adopt the notation
$B(0, r)$ for the \mbox{$d$-}dimensional ball with the center 0 and the radius $r>0$.

By the bound $p_1(x)\le(2\pi)^{-d/2}$,
the right-hand side is bounded
by a constant multiple of
\begin{eqnarray*}
&&\int_{B(0, Hm^{(6-\alpha)/(2(2-\alpha))})}\E_x \biggl[\exp \biggl\{
{\beta\over1-\sigma} \int_0^{t_m-1}\tilde{f}
\biggl({s+1\over t_m}, B(s) \biggr)\,ds \biggr\};  \tau(H)\ge
t_m \biggr]
\\
&&\qquad \le \bigl\llvert B \bigl(0, Hm^{(6-\alpha)/(2(2-\alpha))} \bigr) \bigr\rrvert
\\
&&\quad\qquad{}\times  \exp \biggl\{
\int_0^{t_m-1}\sup_{g\in{\mathcal F}_d} \biggl(
{\beta\over
1-\sigma} \int_{\R^d} \tilde{f} \biggl(
{s+1\over t_m}, x \biggr)g^2(x)\,dx
\\
&&\hspace*{126pt}\quad\qquad{}  -{1\over2} \int_{\R^d}\bigl\llvert
\nabla g(x)\bigr\rrvert ^2\,dx \biggr)\,ds \biggr\},
\end{eqnarray*}
where the inequality follows from (\ref{v-5}) in Lemma \ref{v-3} in the
\hyperref[appendix]{Appendix}.
By variable substitution,
\begin{eqnarray*}
&&\int_0^{t_m-1}\sup_{g\in{\mathcal F}_d}
\biggl({\beta\over1-\sigma} \int_{\R^d} \tilde{f} \biggl(
{s+1\over t_m}, x \biggr)g^2(x)\,dx -
{1\over2} \int_{\R^d}\bigl\llvert \nabla g(x)
\bigr\rrvert ^2\,dx \biggr)\,ds
\\
&&\qquad =t_m\int_{t_m^{-1}}^1\sup
_{g\in{\mathcal F}_d} \biggl({\beta\over
1-\sigma} \int
_{\R^d} \tilde{f}(s, x)g^2(x)\,dx
\\
&&\hspace*{117pt}{} -
{1\over2} \int_{\R^d}\bigl\llvert \nabla g(x)
\bigr\rrvert ^2\,dx \biggr)\,ds
\\
&&\qquad \le t_m\int_0^1\sup
_{g\in{\mathcal F}_d} \biggl({\beta\over1-\sigma} \int
_{\R^d} \tilde{f}(s, x)g^2(x)\,dx -
{1\over2} \int_{\R^d}\bigl\llvert \nabla g(x)
\bigr\rrvert ^2\,dx \biggr)\,ds
\\
&&\qquad =t_m\sup_{g\in{\mathcal A}_d} \biggl({\beta\over1-\sigma}
\int_0^1\!\!\int_{\R^d}
\tilde{f}(s, x)g^2(s, x)\,dx\,ds
\\
&&\hspace*{94pt}{}  -{1\over2} \int_0^1\!\!\int_{\R^d}\bigl
\llvert \nabla_x g(s, x)\bigr\rrvert ^2\,dx\,ds \biggr).
\end{eqnarray*}
Further,
\begin{eqnarray*}
&&\int_0^1\!\!\int_{\R^d}
\tilde{f}(s, x)g^2(s, x)\,dx\,ds
\\
&&\qquad =\int_0^1\!\!\int_{[-N, N]^d}f(s,x) \biggl[\int_{\R^d}Q_N(y-x)g^2(s,y)\,dy
\biggr]\,dx
\\
&&\qquad \le \biggl(\int_0^1\!\!\int
_{[-N, N]^d} \biggl[\int_{\R^d}Q_N(y-x)g^2(s,y)\,dy
\biggr]^2\,dx \biggr)^{1/2}.
\end{eqnarray*}

Summarizing our estimate,
\begin{eqnarray*}
&&\max_{f\in ({\mathcal N}_m^\delta )^*} \E_0 \biggl[\exp \biggl\{
{\beta\over1-\sigma} \int_0^{t_m}\tilde{f}
\biggl({s\over t_m}, B(s) \biggr)\,ds \biggr\}; \tau(H)\ge
t_m \biggr]
\\
&&\qquad \le Cm^{(6-\alpha)d/(2(2-\alpha))}\exp \biggl\{M_{\varepsilon, N} \biggl(
{\beta\over1-\sigma} \biggr)t_m \biggr\}.
\end{eqnarray*}
Here we introduce the notation
\begin{eqnarray*}
M_{\varepsilon, N}(\beta) &=& \sup_{g\in{\mathcal A}_d} \biggl\{ \biggl( \int_0^1\!\!\int_{[-N, N]^d} \biggl[
\int_{\R^d} Q_N(y-x)g^2(s,y)\,dy
\biggr]^2 \,dx \biggr)^{1/2}
\\
&&\hspace*{119pt}{}  -{1\over2} \int_0^1\!\!\int_{\R^d} \bigl
\llvert \nabla_x g(s, x)\bigr\rrvert ^2\,dx\,ds \biggr\}.
\end{eqnarray*}
By (\ref{mo-18}), therefore,
%
%
%e3.23 #&#
\begin{eqnarray}
\label{mo-19}
\qquad&& \limsup_{m\to\infty} m^{-{(4-\alpha)/(2-\alpha)}}\log
\E_0\exp \biggl\{\beta mt_m \biggl( \int_0^1\!\!\int_{\R^d}
\eta^2_m (t_ms,x)\,dx\,ds \biggr)^{1/2}
\biggr\}
\nonumber\\[-8pt]\\[-8pt]\nonumber
&&\qquad  \le tM_{\varepsilon, N} \biggl({\beta\over1-\sigma} \biggr).
\end{eqnarray}
By Lemma \ref{va-1} in the \hyperref[appendix]{Appendix},
\[
\limsup_{N\to\infty}M_{\varepsilon, N} \biggl({\beta\over1-\sigma}
\biggr)\le M_\varepsilon \biggl({\beta\over1-\sigma
} \biggr). %
\]
Finally, the requested (\ref{mo-13}) follows from the obvious fact
that the right-hand side of the above inequality tends to $M_\varepsilon
(\beta)$
as $\sigma\to0^+$.
\end{pf}

By (\ref{mo-11}),
\begin{eqnarray*}
&& \int_0^{t_m}\!\!\int_{\R^d}
\Biggl[\sum_{j=1}^m K_\varepsilon
\bigl(x-B_j(s) \bigr) \Biggr]^2\,dx\,ds
\\
&&\qquad  =mt_m
\gamma_\varepsilon(0)+2\sum_{1\le j<k\le m}\int
_0^{t_m} \gamma_\varepsilon
\bigl(B_j(s)-B_k(s) \bigr)\,ds. %
\end{eqnarray*}
The first term on the right-hand side is deterministic and negligible. Thus
Lemma~\ref{le-1} (with $\beta$ replaced by $\beta/\sqrt{2}$)
can be restated as
%
%
%e3.24 #&#
\begin{eqnarray}
\label{mo-20}
&& \lim_{m\to\infty}m^{-{(4-\alpha)/(2-\alpha)}}\nonumber
\\
&&\quad{}\times \log \E_0\exp \biggl\{\beta \biggl(t_m\sum
_{1\le j<k\le m}\int_0^{t_m}
\gamma_\varepsilon \bigl(B_j(s)-B_k(s) \bigr)\,ds
\biggr)^{1/2} \biggr\}
\\
&&\qquad =tM_\varepsilon \biggl({\beta\over\sqrt{2}}
\biggr).\nonumber
\end{eqnarray}
The next step is to squash $\varepsilon$ to zero.

%
%le3.5 #&#
\begin{lemma}\label{le-2} For any integer $n\ge1$,
\begin{eqnarray*}
&& \E_0 \biggl[\sum_{1\le j<k\le m}\int
_0^{t_m} \gamma_\varepsilon
\bigl(B_j(s)-B_k(s) \bigr)\,ds \biggr]^n
\\
&&\qquad  \le
\E_0 \biggl[\sum_{1\le j<k\le m}\int
_0^{t_m} \gamma \bigl(B_j(s)-B_k(s)
\bigr)\,ds \biggr]^n. %
\end{eqnarray*}
\end{lemma}

\begin{pf} By Fourier transform
\begin{eqnarray*}
\hspace*{-5pt}&&\sum_{1\le j<k\le m}\int_0^{t_m}
\gamma_\varepsilon \bigl(B_j(s)-B_k(s) \bigr)\,ds
\\
\hspace*{-5pt}&&\qquad =(2\pi)^{-d}
\\
\hspace*{-5pt}&&\quad\qquad{}\!\times  \int_{\R^d}\exp \bigl\{-\varepsilon
\llvert \lambda\rrvert ^2 \bigr\}\hat{\gamma }(\lambda) \biggl[\int _0^{t_m}\!\!\sum_{1\le j<k\le m}
\exp \bigl\{-i\lambda\cdot \bigl(B_j(s)-B_k(s) \bigr)
\bigr\}\,ds \biggr]\,d\lambda,
\end{eqnarray*}
where $\hat{\gamma}(\lambda)$ is the Fourier transform of $\gamma(x)$;
see (\ref{intro-33}). Here we shall use the fact that
$\hat{\gamma}(\lambda)>0$ in our setting.
Hence
\begin{eqnarray*}
&&\E_0 \biggl[\sum_{1\le j<k\le m}\int
_0^{t_m} \gamma_\varepsilon
\bigl(B_j(s)-B_k(s) \bigr)\,ds \biggr]^n
\\
&&\qquad =(2\pi)^{-nd}\int_{(\R^d)^n}\,d\lambda_1
\cdots\, d\lambda_n \exp \Biggl\{-\varepsilon\sum
_{l=1}^n\llvert \lambda_l \rrvert
^2 \Biggr\} \Biggl(\prod_{l=1}^n
\hat{\gamma}(\lambda_l) \Biggr)
\\
&&\quad\qquad{}\times \int_{[0,t_m]^n}\E_0\prod
_{l=1}^n\sum_{1\le j<k\le m}\exp
\bigl\{-i\lambda _l\cdot \bigl(B_j(s)-B_k(s)
\bigr) \bigr\}\,ds_1\cdots\, ds_n.
\end{eqnarray*}
Notice that
\[
\E_0\prod_{l=1}^n\sum
_{1\le j<k\le m}\exp \bigl\{-i\lambda_l\cdot
\bigl(B_j(s)-B_k(s) \bigr) \bigr\}>0. %
\]
The right-hand side is less than or equal to
\begin{eqnarray*}
&&(2\pi)^{-nd}\int_{(\R^d)^n}\,d\lambda_1
\cdots\, d\lambda_n \Biggl(\prod_{l=1}^n
\hat{\gamma}(\lambda_l) \Biggr)
\\
&&\quad\times \int_{[0,t_m]^n}\E_0\prod
_{l=1}^n\sum_{1\le j<k\le m}\exp
\bigl\{-i\lambda _l\cdot \bigl(B_j(s)-B_k(s)
\bigr) \bigr\}\,ds_1\cdots\, ds_n
\\
&&\qquad = \E_0 \biggl[\sum_{1\le j<k\le m}\int
_0^{t_m} \gamma \bigl(B_j(s)-B_k(s)
\bigr)\,ds \biggr]^n.
\end{eqnarray*}\upqed
\end{pf}

By using \cite{Chen}, Lemma 1.2.6, page 13, twice with $p=2$,
Lemma \ref{le-2} and (\ref{mo-20}) lead to
\begin{eqnarray*}
&& \lim_{m\to\infty}m^{-{(4-\alpha)/(2-\alpha)}}\log\E_0\exp \biggl\{
\beta \biggl(t_m\sum_{1\le j<k\le m}\int
_0^{t_m} \gamma \bigl(B_j(s)-B_k(s)
\bigr)\,ds \biggr)^{1/2} \biggr\}
\\[-3pt]
&&\qquad  \ge tM_\varepsilon \biggl(
{\beta\over\sqrt{2}} \biggr) %
\end{eqnarray*}
for every $\varepsilon>0$. Notice that
\[
\liminf_{\varepsilon\to0^+}M_\varepsilon \biggl({\beta\over\sqrt{2}}
\biggr) \ge M \biggl({\beta\over\sqrt{2}} \biggr)= \biggl(
{\beta\over\sqrt{2}} \biggr)^{4/(4-\alpha)}M(1),
\]
where
\begin{eqnarray*}
M(\beta) &=& \sup_{g\in{\mathcal A}_d} \biggl\{\beta \biggl(\int_0^1\!\!\int_{\R^d} \biggl[
\int_{\R^d}K(y-x)g^2(s, y)\,dy \biggr]^2\,dx
\biggr)^{1/2}
\\
&&\hspace*{95pt}{}  -{1\over2}\int_0^1\!\!\int_{\R^d}\bigl\llvert \nabla_x g(s,x) \bigr
\rrvert ^2\,dx\,dy \biggr\} %
\end{eqnarray*}
and the second step comes from
the fact that $M(\beta)=\beta^{4/(4-\alpha)}M(1)$ resulted from
replacing $g(s,x)$ by
$ \beta^{d/(4-\alpha)}g(s, \beta^{2/(4-\alpha)}x)$ in
the variation $M(\beta)$.

Hence we reach the lower bound
%
%
%e3.25 #&#
\begin{eqnarray}
\label{mo-21}
&&\liminf_{m\to\infty}m^{-{(4-\alpha)/(2-\alpha)}}\nonumber
\\
&&\quad {}\times \log
\E_0\exp \biggl\{ \beta \biggl(t_m\sum
_{1\le j<k\le m}\int_0^{t_m} \gamma
\bigl(B_j(s)-B_k(s) \bigr)\,ds \biggr)^{1/2}
\biggr\}
\\
&&\qquad \ge t \biggl({\beta\over\sqrt{2}} \biggr)^{4/(4-\alpha)}M(1).\nonumber
\end{eqnarray}

Write $\zeta_\varepsilon(x)=\gamma(x)-\gamma_\varepsilon(x)$.
To have the correspondent upper bound, we prove the following:

%
%le3.6 #&#
\begin{lemma}\label{le-3}
For every $\beta>0$,
\begin{eqnarray*}
&& \lim_{\varepsilon\to0^+}\limsup_{m\to\infty} m^{-{(4-\alpha)/(2-\alpha)}}
\log\E_0\exp \biggl\{{\beta\over m} \sum
_{1\le j<k\le m}\int_0^{t_m}
\zeta_\varepsilon \bigl(B_j(s)-B_k(s) \bigr)\,ds
\biggr\}
\\
&&\qquad  =0. %
\end{eqnarray*}
\end{lemma}

\begin{pf} By Jensen's inequality
\begin{eqnarray*}
&&\E_0\exp \biggl\{{\beta\over m} \sum
_{1\le j<k\le m}\int_0^{t_m}
\zeta_\varepsilon \bigl(B_j(s)-B_k(s) \bigr)\,ds
\biggr\}
\\
&&\qquad \ge\exp \biggl\{{\beta\over m} \E_0\sum
_{1\le j<k\le m}\int_0^{t_m}
\zeta_\varepsilon \bigl(B_j(s)-B_k(s) \bigr)\,ds
\biggr\} \ge1,
\end{eqnarray*}
where the second inequality follows from Lemma \ref{le-2} with $n=1$.
Thus we only need to prove the upper bound estimate.

Write
\[
\sum_{1\le j<k\le m}\int_0^{t_m}
\zeta_\varepsilon \bigl(B_j(s)-B_k(s) \bigr)\,ds =
{1\over2}\sum_{j=1}^m\sum
_{k\dvtx  k\neq j} \int_0^{t_m}
\zeta_\varepsilon \bigl(B_j(s)-B_k(s) \bigr)\,ds.
\]
By H\"older's inequality
\begin{eqnarray*}
&&\E_0\exp \biggl\{{\beta\over m} \sum
_{1\le j<k\le m}\int_0^{t_m}
\zeta_\varepsilon \bigl(B_j(s)-B_k(s) \bigr)\,ds
\biggr\}
\\
&&\qquad \le\prod_{j=1}^m \biggl(
\E_0\exp \biggl\{{\beta\over2} \sum
_{k\dvtx  k\neq j} \int_0^{t_m}
\zeta_\varepsilon \bigl(B_j(s)-B_k(s) \bigr)\,ds
\biggr\} \biggr)^{1/m}
\\
&&\qquad =\E_0\exp \Biggl\{{\beta\over2}\sum
_{k=2}^m\int_0^{t_m}
\zeta_\varepsilon \bigl(B_1(s)-B_k(s) \bigr)\,ds
\Biggr\}.
\end{eqnarray*}

We now make use of Fourier transform again. Notice that
$\hat{\zeta}_\varepsilon(\lambda)=(1-e^{-\varepsilon\llvert  \lambda \rrvert  ^2})
\hat{\gamma}(\lambda)\ge0$. By Fourier inversion
\begin{eqnarray*}
&&\sum_{k=2}^m\int_0^{t_m}
\zeta_\varepsilon \bigl(B_1(s)-B_k(s) \bigr)\,ds
\\
&&\qquad =(2\pi)^{-d}\int_{\R^d} \hat{\zeta}_\varepsilon(
\lambda) \Biggl[\int_0^{t_m}\sum
_{k=2}^m\exp \bigl\{-i \lambda\cdot
\bigl(B_1(s)-B_k(s) \bigr) \bigr\}\,ds \Biggr]\,d\lambda.
\end{eqnarray*}
For any integer $n\ge1$, by the independence between $B_1(s)$ and
$\{B_2(s),\ldots, B_m(s)\}$,
\begin{eqnarray*}
&&\E_0 \Biggl[\sum_{k=2}^m
\int_0^{t_m} \zeta_\varepsilon
\bigl(B_1(s)-B_k(s) \bigr)\,ds \Biggr]^n
\\
&&\qquad =(2\pi)^{-nd}\int_{(\R^d)^n}\,d\lambda_1
\cdots\, d\lambda_m \Biggl(\prod_{l=1}^n
\hat{\zeta}_\varepsilon(\lambda_l) \Biggr)
\\
&&\quad\qquad{}\times \int
_{[0,t_m]^n}\E_0\exp \Biggl\{-i\sum
_{l=1}^n\lambda_l\cdot
B_1(s_l) \Biggr\}
\\
&&\hspace*{75pt}{}\times\E_0 \Biggl(\prod_{l=1}^n
\sum_{k=2}^m\exp \bigl\{i
\lambda_l\cdot B_k(s_l) \bigr\}
\Biggr)\,ds_1\cdots\, ds_n.
\end{eqnarray*}
By the fact that
\[
0<\E_0\exp \Biggl\{-i\sum_{l=1}^n
\lambda_l\cdot B_1(s_l) \Biggr\} \le1
\quad \mbox{and}\quad \E_0 \Biggl(\prod
_{l=1}^n\sum_{k=2}^m
\exp \bigl\{i \lambda_l\cdot B_k(s_l)
\bigr\} \Biggr)>0, %
\]
the right-hand side is less than or equal to
\begin{eqnarray*}
&&(2\pi)^{-nd}\int_{(\R^d)^n}\,d\lambda_1
\cdots\, d\lambda_m \Biggl(\prod_{l=1}^n
\hat{\zeta}_\varepsilon(\lambda_l) \Biggr)
\\[-1pt]
&&\quad{}\times  \int
_{[0,t_m]^n}\E_0 \Biggl(\prod
_{l=1}^n\sum_{k=2}^m
\exp \bigl\{i \lambda_l\cdot B_k(s_l)
\bigr\} \Biggr)\,ds_1\cdots\, ds_n
\\[-1pt]
&&\qquad =\E_0 \Biggl[\sum_{k=2}^m
\int_0^{t_m} \zeta_\varepsilon
\bigl(B_k(s) \bigr)\,ds \Biggr]^n.
\end{eqnarray*}
Therefore, for any $n=1,2,\ldots$
\[
\E_0 \Biggl[\sum_{k=2}^m
\int_0^{t_m} \zeta_\varepsilon
\bigl(B_1(s)-B_k(s) \bigr)\,ds \Biggr]^n \le
\E_0 \Biggl[\sum_{k=2}^m\int
_0^{t_m} \zeta_\varepsilon
\bigl(B_k(s) \bigr)\,ds \Biggr]^n. %
\]
By Taylor expansion we conclude that
\begin{eqnarray*}
&&\E_0\exp \Biggl\{{\beta\over2}\sum
_{k=2}^m\int_0^{t_m}
\zeta_\varepsilon \bigl(B_1(s)-B_k(s) \bigr)\,ds
\Biggr\}
\\
&&\qquad \le \E_0\exp \Biggl\{{\beta\over2}\sum
_{k=2}^m\int_0^{t_m}
\zeta_\varepsilon \bigl(B_k(s) \bigr)\,ds \Biggr\}
= \biggl(\E_0\exp \biggl\{{\beta\over2}\int
_0^{t_m} \zeta_\varepsilon \bigl(B(s)
\bigr)\,ds \biggr\} \biggr)^{m-1}.
\end{eqnarray*}

Summarizing our argument, we have reduced the problem to the proof of
%
%
%e3.26 #&#
\begin{eqnarray}
\label{mo-22} \limsup_{\varepsilon\to0^+}\lim_{m\to\infty}
{1\over t_m} \log\E_0\exp \biggl\{\beta\int
_0^{t_m} \zeta_\varepsilon \bigl(B(s)
\bigr)\,ds \biggr\}\le0
\end{eqnarray}
for any $\beta>0$.

For the sake of simplicity we consider the case when $t_m$ goes to
infinity along
the integer times. Notice that $\hat{\zeta}_\varepsilon(\lambda)>0$ for all
$\lambda\in\R^d$.
Using the same argument as that used in the proof of
Lemma \ref{le-2}, one can show that for any $x\in\R^d$,
\[
\E_x \biggl[\int_0^1
\zeta_\varepsilon \bigl(B(s) \bigr)\,ds \biggr]^n \le
\E_0 \biggl[\int_0^1
\zeta_\varepsilon \bigl(B(s) \bigr)\,ds \biggr]^n, \qquad n=1,2,
\ldots. %
\]
By Taylor expansion
\[
\E_x\exp \biggl\{\beta\int_0^1
\zeta_\varepsilon \bigl(B(s) \bigr)\,ds \biggr\} \le\E_0\exp
\biggl\{\beta\int_0^1 \zeta_\varepsilon
\bigl(B(s) \bigr)\,ds \biggr\},\qquad x\in\R^d. %
\]
By Markov's property,
\begin{eqnarray*}
&&\E_0\exp \biggl\{\beta\int_0^{t_m}
\zeta_\varepsilon \bigl(B(s) \bigr)\,ds \biggr\}
\\
&&\qquad =\E_0 \biggl[\exp \biggl\{\beta\int_0^{t_m-1}
\zeta_\varepsilon \bigl(B(s) \bigr)\,ds \biggr\}\E_{B(t_m-1)}\exp \biggl\{
\beta\int_0^1 \zeta_\varepsilon
\bigl(B(s) \bigr)\,ds \biggr\} \biggr]
\\
&&\qquad \le\E_0\exp \biggl\{\beta\int_0^{t_m-1}
\zeta_\varepsilon \bigl(B(s) \bigr)\,ds \biggr\}\E_0\exp \biggl\{
\beta\int_0^1 \zeta_\varepsilon
\bigl(B(s) \bigr)\,ds \biggr\}.
\end{eqnarray*}
Continuing this procedure we have
\[
\E_0\exp \biggl\{\beta\int_0^{t_m}
\zeta_\varepsilon \bigl(B(s) \bigr)\,ds \biggr\} \le \biggl(\E_0
\exp \biggl\{\beta\int_0^1
\zeta_\varepsilon \bigl(B(s) \bigr)\,ds \biggr\} \biggr)^{t_m}.
\]
Finally, the requested (\ref{mo-22}) follows from the obvious fact that
\[
\lim_{\varepsilon\to0^+}\E_0\exp \biggl\{\beta\int
_0^1 \zeta_\varepsilon \bigl(B(s) \bigr)\,ds
\biggr\}=1.
\]\upqed
\end{pf}

Write
\[
Z_m=\sum_{1\le j<k\le m} \int
_0^{t_m}\zeta_\varepsilon
\bigl(B_j(s)-B_k(s) \bigr)\,ds %
\]
and $Z_m^+=\max\{0, Z_m\}$. Given $\delta>0$ and $\beta>0$
\begin{eqnarray*}
\E_0\exp \Bigl\{\beta\sqrt{t_m
Z_m^+} \Bigr\} &\le&\exp \{\beta mt_m\delta \}+
\E_0\exp \bigl[ \{\beta\sqrt{t_m Z_m} \}; Z_m\ge\delta^2m^2t_m
\bigr]
\\
&\le&\exp \{\beta mt_m\delta \}+ \E_0\exp \biggl\{
{\beta\over\delta m}Z_m \biggr\}.
\end{eqnarray*}
By Lemma \ref{le-3},
\[
\limsup_{\varepsilon\to0^+}\limsup_{m\to\infty}m^{-{(4-\alpha)/(2-\alpha) }}
\log \E_0\exp \Bigl\{\beta\sqrt{t_m
Z_m^+} \Bigr\} \le\beta t\delta. %
\]
Since $\delta>0$ can be arbitrarily small,
%
%
%e3.27 #&#
\begin{eqnarray}
\label{mo-23} \lim_{\varepsilon\to0^+}\limsup_{m\to\infty}m^{-{(4-\alpha)/(2-\alpha) }}
\log \E_0\exp \Bigl\{\beta\sqrt{t_m
Z_m^+} \Bigr\} =0.
\end{eqnarray}

We now return to the variation $M_\varepsilon(\beta)$
introduced at the beginning of this subsection. By Jensen's
inequality for any $g\in{\mathcal A}_d$,
\[
\int_{\R^d} \biggl[\int_{\R^d}K_\varepsilon(y-x)g^2(s,
y)\,dy \biggr]^2\,dx \le\int_{\R^d} \biggl[\int
_{\R^d}K(y-x)g^2(s, y)\,dy \biggr]^2\,dx.
\]
Consequently, $M_\varepsilon(\beta)\le M(\beta)$ for any $\beta>0$.

In view of (\ref{mo-20}) and (\ref{mo-23}), a standard argument
of exponential approximation via H\"older inequality leads to
%
%
%e3.28 #&#
\begin{eqnarray}
\label{mo-24} &&\limsup_{m\to\infty}m^{-{(4-\alpha)/(2-\alpha)}}\nonumber
\\
&&\quad{}\times \log
\E_0\exp \biggl\{\beta \biggl(t_m\sum
_{1\le j<k\le m}\int_0^{t_m} \gamma
\bigl(B_j(s)-B_k(s) \bigr)\,ds \biggr)^{1/2}
\biggr\}
\\
&&\qquad \le tM \biggl({\beta\over\sqrt{2}} \biggr) =t \biggl({\beta\over\sqrt{2}}
\biggr)^{4/(4-\alpha)}M(1).
\nonumber
\end{eqnarray}
Combining (\ref{mo-21}) and (\ref{mo-24})
\begin{eqnarray*}
&& \lim_{m\to\infty}m^{-{(4-\alpha)/(2-\alpha)}}\log \E_0\exp \biggl
\{\beta \biggl(t_m\sum_{1\le j<k\le m}\int
_0^{t_m} \gamma \bigl(B_j(s)-B_k(s)
\bigr)\,ds \biggr)^{1/2} \biggr\}
\\
&&\qquad =t \biggl({\beta\over\sqrt{2}}
\biggr)^{4/(4-\alpha)}M(1).
\end{eqnarray*}
Replacing $\beta$ by $t^{-1/2}\beta$ leads to
%
%
%e3.29 #&#
\begin{eqnarray}\label{mo-25}
&&\lim_{m\to\infty}m^{-{(4-\alpha)/(2-\alpha)}}\nonumber
\\
&&\quad{} \times \log \E_0\exp \biggl\{\beta m^{(4-\alpha)/(2(2-\alpha))}
\nonumber\\[-8pt]\\[-8pt]\nonumber
&&\hspace*{71pt}{}\times  \biggl({1\over m}
\sum_{1\le j<k\le m}\int_0^{t_m}
\gamma \bigl(B_j(s)-B_k(s) \bigr)\,ds
\biggr)^{1/2} \biggr\}
\\
&&\qquad =t^{(2-\alpha)/(4-\alpha)} \biggl({\beta\over\sqrt{2}} \biggr)^{4/(4-\alpha)}M(1).\nonumber
\end{eqnarray}

In addition, notice that for any $\varepsilon>0$ there is a constant
$C_\varepsilon>0$
such that $\gamma_\varepsilon(x)\le C_\varepsilon$ for all $x\in\R^d$. Thus
\[
{1\over m}\sum_{1\le j<k\le m}\int
_0^{t_m}\gamma_\varepsilon
\bigl(B_j(s)-B_k(s) \bigr)\,ds \le C_\varepsilon
mt_m=C_\varepsilon tm^{(4-\alpha)/(2-\alpha)}. %
\]
Together with Lemma \ref{le-3}, this implies that for every $\beta>0$
\[
\limsup_{m\to\infty} m^{-{(4-\alpha)/(2-\alpha)}}\log\E_0\exp
\biggl\{{\beta\over m} \sum_{1\le j<k\le m}\int
_0^{t_m}\gamma \bigl(B_j(s)-B_k(s)
\bigr)\,ds \biggr\} <\infty. %
\]
Using the second half of Lemma \ref{A-1} in the \hyperref[appendix]{Appendix} with
\[
b_m=m^{(4-\alpha)/(2-\alpha)},\qquad p={2\over2-\alpha},
\qquad C_0=t{2-\alpha\over4} \biggl(
{4M(1)\over4-\alpha} \biggr)^{(4-\alpha)/(2-\alpha)} %
\]
and
\[
X_m={1\over m}\sum_{1\le j<k\le m}
\int_0^{t_m}\gamma \bigl(B_j(s)-B_k(s)
\bigr)\,ds, %
\]
we obtain
\begin{eqnarray*}
&&\lim_{m\to\infty} m^{-{(4-\alpha)/(2-\alpha)}}\log\E_0\exp \biggl
\{{\beta\over m} \sum_{1\le j<k\le m}\int
_0^{t_m}\gamma \bigl(B_j(s)-B_k(s)
\bigr)\,ds \biggr\}
\\
&&\qquad =t{2-\alpha\over4} \biggl({4M(1)\over4-\alpha}
\biggr)^{(4-\alpha)/(2-\alpha)} \beta^{2/(2-\alpha)} =t{\mathcal E}(d, \gamma) \biggl(
{\beta\over2} \biggr)^{2/(2-\alpha)},
\end{eqnarray*}
where the last equality follows from Lemma \ref{va-2} in the \hyperref[appendix]{Appendix}.
By the identity in the law
\[
{1\over m} \sum_{1\le j<k\le m}\int
_0^{t_m}\gamma \bigl(B_j(s)-B_k(s)
\bigr)\,ds \buildrel d\over=\sum_{1\le j<k\le m}\int
_0^{t}\gamma \bigl(B_j(s)-B_k(s)
\bigr)\,ds, %
\]
we have
%
%
%e3.30 #&#
\begin{eqnarray}\label{mo-26}
&& \lim_{m\to\infty} m^{-{(4-\alpha)/(2-\alpha)}}\log
\E_0\exp \biggl\{\beta \sum_{1\le j<k\le m}\int
_0^{t}\gamma \bigl(B_j(s)-B_k(s)
\bigr)\,ds \biggr\}
\nonumber\\[-8pt]\\[-8pt]\nonumber
&&\qquad  =t{\mathcal E}(d, \gamma) \biggl({\beta\over2}
\biggr)^{2/(2-\alpha)}.
\end{eqnarray}

Proposition \ref{prop-3} follows from representations
(\ref{mo-9}) and (\ref{mo-26}) (with $\beta=\theta^2$).
%\end{pf}

%
%re3.7 #&#
\begin{remark}\label{remark-1}
Bertini and Cancrini (Theorem 2.6,
\cite{BC})
claimed
a precise formula for $\E u(t, 0)^m$ in the setting of Theorem \ref{th-4}.
Unfortunately,
their result is false due to incorrectly using the Skorokhod lemma.
On the other hand, (\ref{mo-8}) shows that the relation in
Bertini--Cancrini's formulation is asymptotically sound.
\end{remark}

%
%re3.8 #&#
\begin{remark}\label{remark}
By Theorem 6.1 in \cite{CHSX},
under the assumptions of Theorem \ref{th-2},
%
%
%e3.31 #&#
\begin{eqnarray}
\label{mo-27}
&& \lim_{t\to\infty}t^{-{(4-\alpha-2\alpha_0)/(2-\alpha)}}\log\E
u(t,0)^m
\nonumber\\[-8pt]\\[-8pt]\nonumber
&&\qquad =m^{(4-\alpha)/(2-\alpha)} \biggl({\theta^2\over2}
\biggr)^{2/(2-\alpha)} {\mathcal E}(\alpha_0, d, \gamma)
\end{eqnarray}
for every integer $m\ge1$.
Comparing this to (\ref{mo-3}), we find the $m$-limit and the $t$-limit are
completely consistent. The same can be claimed in the context of
Theorem \ref{th-2prime}. The situation
is slightly different when it
comes to the cases labeled (2) in Table~\ref{table1}
where $V(t,x)$ is white in time. Take the setting
of Theorem \ref{th-4}, for example. Let $m=2$ in (\ref{mo-9}),
%
%
%e3.32 #&#
\begin{eqnarray}
\label{mo-29} \E u(t,0)^2&=&\E_0\exp \biggl\{
\theta^2\int_0^t
\delta_0 \bigl(B_1(s)-B_2(s) \bigr) \,ds
\biggr\}
\nonumber\\[-8pt]\\[-8pt]
&=&\E_0\exp \biggl\{{\theta^2\over\sqrt{2}} \int
_0^t\delta_0 \bigl(B(s) \bigr)\,ds
\biggr\} =\E_0\exp \biggl\{{\theta^2\over\sqrt{2}}\bigl\llvert
B(t)\bigr\rrvert \biggr\},
\nonumber
\end{eqnarray}
where the last equation follows the well-known identity in law between
the Brownian local time and the reflected Brownian motion. Thus
%
%
%e3.33 #&#
\begin{equation}
\label{mo-28} \lim_{t\to\infty}{1\over t}\log\E u(t,
0)^2={\theta^4\over4}=2^3 {\theta^4\over32}.
\end{equation}
In comparison to (\ref{mo-8}) [keep in mind that $m=2$ in (\ref{mo-28})],
we witness a small but interesting gene mutation occurring
during the course $m\to\infty$.
\end{remark}

%s4 #&#
\section{Modulus continuity}\label{con}

The main goal of this section is to measure the degree of the
continuity of $u(t,x)$ in the space variable $x$ by estimating
the difference $u(t,x)-u(t,y)$. In the settings of Theorems \ref{th-3}~and~\ref{th-4}
we use the bound
\begin{eqnarray*}
&& \bigl(\E\bigl\llvert u (t,x)-u(t,y)\bigr\rrvert ^{2m}
\bigr)^{1/(2m)}
\\
&&\qquad \le C\llvert x-y\rrvert + \biggl(8m\int_0^t
\bigl(\E u (s, 0)^{2m} \bigr)^{1/m} {\mathcal
I}_s\,ds \biggr)^{1/2}
\\
&&\qquad \le C\llvert x-y\rrvert + \bigl(\E u(t, 0)^{2m}
\bigr)^{1/(2m)} \biggl(8m\int_0^t{
\mathcal I}_{t-s}\,ds \biggr)^{1/2}
\end{eqnarray*}
established by Conus et al. ((9.49), \cite{CJKS}),
where
\begin{eqnarray*}
{\mathcal I}_s&=&\int_{\R^d\times\R^d}\gamma(z_1-z_2)h_s(z_1)h_s(z_2)\,dz_1\,dz_2,
\\
h_s(z)&=&\bigl\llvert p_s(z-x)-p_s(z-y)
\bigr\rrvert, \qquad z\in\R^d %
\end{eqnarray*}
and $C>0$ represents, here and else where in this section, a constant
independent of $m$ and $x, y$ that takes possibly different values
when appearing in different places.

By (9.51) in \cite{CJKS},
\[
{\mathcal I}_s\le C(t-s)^{-{\alpha/2}}\cdot \biggl(
{\llvert   x-y\rrvert  \over
(t-s)^{1/2}}\wedge 1 \biggr). %
\]
Thus, for $\llvert   x-y\rrvert  \le \sqrt{t}$,
\[
\int_0^t{\mathcal I}_{t-s}\,ds\le C
\biggl\{\int_0^{\llvert   x-y\rrvert  ^2} s^{-{\alpha/2}}\,ds+\llvert
x-y\rrvert \int_{\llvert   x-y\rrvert  ^2}^ts^{-{(\alpha +1)/2}} \,ds
\biggr\} \le C\llvert x-y\rrvert. %
\]
This estimate gives the bound
\[
\E\bigl\llvert u(t,x)-u(t,y)\bigr\rrvert ^{2m} \le C^mm!
\llvert x-y \rrvert ^m\E u(t, 0)^{2m} %
\]
or
%
%
%e4.1 #&#
\begin{equation}
\label{con-1} \E\bigl\llvert u(t,x)-u(t,y)\bigr\rrvert ^{m} \le
C^m(m!)^{1/2}\llvert x-y\rrvert ^{m/2}\E u (t,
0)^{m}.
\end{equation}
By the classic theory on chaining method (see, e.g., Lemma 9, \cite{CLR}),
(\ref{con-1}) leads to:

%
%le4.1 #&#
\begin{lemma}\label{le-4} In the settings of Theorems \ref{th-3}~and~\ref{th-4}, $u(t,x)$
yields a continuous modification. Moreover, for any $0<\delta<1$ and
bounded domain $D\subset\R^d$ there is
a $C_\delta(D)>0$ such that for $m\delta>2d$,
%
%
%e4.2 #&#
\begin{eqnarray}
\label{con-2} \E\mathop{\sup_{x\neq y}}_{x, y
\in D} \biggl\llvert
{u(t,x)-u(t,y)\over\llvert   x-y\rrvert  ^{\delta/2}} \biggr\rrvert ^m \le C_\delta(D)
(m!)^{1/2}\E u(t,0)^m.
\end{eqnarray}
\end{lemma}

We now consider the setting of Theorems \ref{th-1}, \ref{th-2}
and \ref{th-2prime}
where the solution $u(t,x)$ yields the
Feynman--Kac representation (\ref{intro-13})
[with $u_0(x)\equiv1$ according to our agreement]. For any $p>1$, write
\[
u_p(t, x)=\E_x\exp \biggl\{p\theta\int
_0^tV \bigl(t-s, B(s) \bigr)\,ds \biggr\}.
\]

%
%le4.2 #&#
\begin{lemma}\label{le-5} In the setting of Theorems \ref{th-1}, \ref
{th-2} and \ref{th-2prime}, $u(t,x)$
yields a continuous modification. Moreover, for any $p>1$ such that
$q\equiv p(p-1)^{-1}$ is an even number and for
any bounded
domain $D\subset\R^d$, there is
a $C_p(D)>0$ and $\delta>0$ such that for $m\delta>2d$:

\begin{longlist}[(2)]
\item[(1)] in the setting of Theorem \ref{th-1},
%
%
%e4.3 #&#
\begin{eqnarray}
\label{con-2prime} \E\mathop{\sup_{x\neq y}}_{x, y
\in D} \biggl\llvert
{u(t,x)-u(t,y)\over\llvert   x-y\rrvert  ^{\delta/2}} \biggr\rrvert ^m \le(m!)^{1/2}C_p(D)^m
\bigl\{\E u_p(t,0)^m \bigr\}^{1/p};
\end{eqnarray}

\item[(2)] in the settings of Theorems \ref{th-2} and \ref{th-2prime},
%
%
%e4.4 #&#
\begin{eqnarray}
\label{con-3} \E\mathop{\sup_{x\neq y}}_{x, y
\in D} \biggl\llvert
{u(t,x)-u(t,y)\over\llvert   x-y\rrvert  ^{\delta/2}} \biggr\rrvert ^m \le m!C_p(D)^m
\bigl\{\E u_p(t,0)^m \bigr\}^{1/p}.
\end{eqnarray}
\end{longlist}
\end{lemma}

\begin{pf} The main part of the proof is to establish a bound similar to
(\ref{con-1}).
By the mean-value theorem,
\[
\bigl\llvert e^\xi-e^\eta\bigr\rrvert \le\llvert \xi-\eta
\rrvert \max \bigl\{e^\xi, e^\eta \bigr\}. %
\]
By the Feynman--Kac representation (\ref{intro-13}) and H\"older's inequality,
for any $y\in\R^d$
\begin{eqnarray*}
&& \bigl\llvert u(t,0)- u(t,y) \bigr\rrvert
\\
&&\qquad  \le\E_0 \biggl\llvert
\exp \biggl\{\theta\int_0^tV \bigl(t-s, B(s)
\bigr)\,ds \biggr\}- \exp \biggl\{\theta\int_0^tV
\bigl(t-s, y+B(s) \bigr)\,ds \biggr\} \biggr\rrvert
\\
&&\qquad \le\theta
\E_0 \biggl( \biggl\llvert \int_0^tV
\bigl(t-s, B(s) \bigr)\,ds -\int_0^tV
\bigl(t-s, y+B(s) \bigr)\,ds \biggr\rrvert
\\
&&\hspace*{55pt}{} \times\max \biggl\{\exp \biggl\{
\theta\int_0^tV \bigl(t-s, B(s) \bigr)\,ds
\biggr\},
\\
&&\hspace*{34pt}\hspace*{93pt} \exp \biggl\{\theta\int_0^tV
\bigl(t-s, y+B(s) \bigr)\,ds \biggr\} \biggr\} \biggr)
\\
&&\qquad \le 2\theta \biggl(
\E_0 \biggl\llvert \int_0^tV
\bigl(t-s, B(s) \bigr)\,ds -\int_0^tV \bigl(s,
y+B(s) \bigr)\,ds \biggr\rrvert ^q \biggr)^{1/q}
\\
&&\quad\qquad{}\times  \bigl
\{u_p(t, 0)+u_p(t, y) \bigr\}^{1/p}.
\end{eqnarray*}
By H\"older's inequality again,
\begin{eqnarray*}
&&\E\bigl\llvert u(t,0)-u(t,y)\bigr\rrvert ^m
\\
&&\qquad \le (2\theta)^m \biggl\{\E \biggl(\E_0 \biggl\llvert
\int_0^tV \bigl(t-s, B(s) \bigr)\,ds -\int
_0^tV \bigl(t-s, y+B(s) \bigr)\,ds \biggr\rrvert
^q \biggr)^m \biggr\}^{1/q}
\\
&&\quad\qquad{}\times \bigl\{\E \bigl(u_p(t, 0)+u_p(t, y)
\bigr)^m \bigr\}^{1/p}.
\end{eqnarray*}
Notice that
\begin{eqnarray*}
&&\E \biggl(\E_0 \biggl\llvert \int_0^tV
\bigl(t-s, B(s) \bigr)\,ds -\int_0^tV
\bigl(t-s, y+B(s) \bigr)\,ds \biggr\rrvert ^q \biggr)^m
\\
&&\qquad \le \E\otimes\E_0 \biggl\llvert \int_0^tV
\bigl(t-s, B(s) \bigr)\,ds -\int_0^tV
\bigl(t-s, y+B(s) \bigr)\,ds \biggr\rrvert ^{qm}.
\end{eqnarray*}
By the triangle inequality and by the stationarity of $u_p(t,x)$ in $x$,
\[
\bigl\{\E \bigl(u_p(t, 0)+u_p(t, y)
\bigr)^m \bigr\}^{1/p}\le2^{{(p+1)/p}m} \bigl(\E
u_p(t,0)^m \bigr)^{1/p}. %
\]
Set
\[
S_t(y)= \biggl\{\int_0^t\!\!\int
_0^t\gamma_0(r-s) \bigl(\gamma
\bigl(B(r)-B(s) \bigr) -\gamma \bigl(y+B(r)-B(s) \bigr) \bigr)\,dr\,ds \biggr
\}^{1/2}. %
\]

Notice the fact that the difference
\[
\int_0^tV \bigl(t-s, B(s) \bigr)\,ds -\int
_0^tV \bigl(t-s, y+B(s) \bigr)\,ds %
\]
is a Gaussian conditioning on the Brownian motion with conditional
variance $2S_t(y)^2$.
By the (conditional) Gaussian property,
\begin{eqnarray*}
&& \E_0\otimes\E \biggl[\int_0^tV(t-s,
B_s)\,ds -\int_0^tV(t-s,
y+B_s)\,ds \biggr]^{qm}
\\
&&\qquad =(qm-1)!!(\sqrt{2})^{qm}
\E_0S_t(y)^{qm}. %
\end{eqnarray*}
So we have
%
%
%e4.5 #&#
\begin{eqnarray}\label{con-4}
&& \E\bigl\llvert u(t,0)-u(t,y)\bigr\rrvert ^m
\nonumber\\[-8pt]\\[-8pt]\nonumber
&&\qquad \le2^{{(p+1)/p}m}( \sqrt{2}\theta)^m \bigl((qm-1)!!
\E_0S_t(y)^{qm} \bigr)^{1/q} \bigl(
\E u_p(t,0)^m \bigr)^{1/p}.
\end{eqnarray}

Let $\hat{\gamma}(\lambda)$ be the Fourier transform [see (\ref
{intro-33})] of
$\gamma(x)$. By Fourier inversion
\begin{eqnarray*}
S_t(y)^2 &=&(2\pi)^{-d}\int
_{\R^d}\hat{\gamma}(\lambda) \bigl[1-e^{-i\lambda\cdot
y} \bigr]
\\
&&\hspace*{50pt}{}\times
\biggl[\int_0^t\!\!\int_0^t
\gamma_0(r-s) \exp \bigl\{i\lambda\cdot \bigl(B(s)-B(r) \bigr)
\bigr\}\,dr\,ds \biggr]\,d\lambda.
\end{eqnarray*}
In the setting of Theorem \ref{th-1}
\begin{eqnarray*}
S_t^2&\le& C_\delta\llvert y\rrvert
^\delta \int_{\R^d}\llvert \lambda \rrvert
^\delta \bigl\llvert \hat{ \gamma}(\lambda)\bigr\rrvert \biggl\llvert
\int_0^t\!\!\int_0^t
\gamma_0(r-s) \exp \bigl\{i \lambda\cdot \bigl(B(s)-B(r) \bigr)
\bigr\}\,dr\,ds \biggr\rrvert \,d\lambda
\\
&\le& C_\delta\llvert y\rrvert ^\delta \biggl(\int
_{\R^d}\llvert \lambda \rrvert ^\delta \bigl\llvert
\hat{\gamma}( \lambda)\bigr\rrvert \,d\lambda \biggr) \int_0^t\!\!\int_0^t\gamma_0(r-s)\,dr\,ds,
\end{eqnarray*}
where $\delta>0$ is chosen according to the assumption (\ref{intro-18}).
By
(\ref{con-4}), by Stirling's formula and by the stationary of $u(t,x)$
in $x$, we reach the bound
\[
\E\bigl\llvert u(t,0)-u(t,y)\bigr\rrvert ^m\le C(D)^m(m!)^{1/2}
\llvert y\rrvert ^{m\delta/2} \bigl(\E u_p(t,0)^m
\bigr)^{1/p}, %
\]
which leads to (\ref{con-2prime}) with possibly smaller $\delta$ and larger
$C(D)$.

We now come to the setting of Theorems \ref{th-2}~and~\ref{th-2prime}.
Notice that $\hat{\gamma}(\lambda)$ is equal to a positive constant
multiple of $\llvert  \lambda \rrvert  ^{-(d-\alpha)}$,
$ \prod_{i=1}^d\llvert  \lambda_i\rrvert  ^{-(1-\alpha_i)}$
[in the notation of $\lambda=(\lambda_1,\ldots,\lambda_d)$]
and $1$ in connection to $(1)\times\mathrm{(I)}$, $(1)\times\mathrm{(II)}$
and $(1)\times\mathrm{(III)}$ (labeled in Table~\ref{table1}), respectively.
By the first relation in (\ref{lo-12}) and the representation
of $S_t(y)$ given above,
\begin{eqnarray*}
S_t^2(y)&=&C\int_{\R\times\R^d}
\bigl[1-e^{-i\lambda\cdot y} \bigr]\hat {\gamma}(\lambda) \biggr\llvert \int
_0^t| u-s|
^{-{(1+\alpha_0)/2}} \exp \bigl\{i\lambda\cdot B(s) \bigr
\}\,ds \biggr\rrvert ^2\,du\,d\lambda
\\
&\le& C_\delta| y| ^\delta \int
_{\R\times\R^d} \Gamma_\delta(\lambda) \biggl\llvert \int
_0^t | u-s|
^{-{(1+\alpha_0)/2}} \exp \bigl\{i\lambda\cdot B(s)  \bigr\}\,ds \biggr\rrvert
^2\,du\,d\lambda,
\end{eqnarray*}
where $\Gamma_\delta(\lambda)=\llvert  \lambda \rrvert  ^\delta\hat{\gamma
}(\lambda)$ and
$\delta>0$ is a small number. We claim that for sufficiently small
$\delta$
the process
%
%e4.6 #&#
\begin{eqnarray}
Z_T= \biggl(\int_{\R\times\R^d}\Gamma_\delta(
\lambda) \biggl\llvert \int_0^T|
u-s | ^{-{(1+\alpha_0)/2}} \exp \bigl\{i\lambda\cdot B(s)
\bigr\}\,ds \biggr\rrvert ^2\,du\,d\lambda
\biggr)^{1/2},\nonumber
\\
\eqntext{T\ge0}
\end{eqnarray}
takes finite values almost surely. For the sake of simplicity
we show this by controlling
$\E Z_T^2$ through a ``usual'' computation without justification, which
is easy to be installed.

By the first relation in (\ref{lo-12}) $Z_T^2$ is a constant multiple of
\[
\int_{\R^d}\Gamma_\delta(\lambda) \biggl[\int _0^T\!\!\int_0^T
\llvert r-s\rrvert ^{-\alpha_0}\exp \bigl\{i\lambda\cdot \bigl(B(r)-B(s)
\bigr) \bigr\}\,dr\,ds \biggr]\,d\lambda %
\]
whose expectation is equal to
\begin{eqnarray*}
&&\int_{\R^d}\Gamma_\delta(\lambda) \biggl[\int _0^T\!\!\int_0^T
\llvert r-s\rrvert ^{-\alpha_0}\E_0 \exp \bigl\{i\lambda\cdot
\bigl(B(r)-B(s) \bigr) \bigr\}\,dr\,ds \biggr]\,d\lambda
\\
&&\qquad =\int_{\R^d}\Gamma_\delta(\lambda) \biggl[\int _0^T\!\!\int_0^T
\llvert r-s\rrvert ^{-\alpha_0}\exp \biggl\{-{\llvert  r-s\rrvert  \over2}\llvert
\lambda \rrvert ^2 \biggr\}\,dr\,ds \biggr]\,d\lambda
\\
&&\qquad =\int _0^T\!\!\int_0^T
\llvert r-s\rrvert ^{-d/2}\llvert r-s\rrvert ^{-\alpha_0}
\\
&&\hspace*{63pt}{}\times  \biggl[\int
_{\R^d} \Gamma_\delta \biggl({\lambda\over\sqrt{\llvert  r-s\rrvert  }}
\biggr) \exp \biggl\{-{1\over2}\llvert \lambda \rrvert
^2 \biggr\}\,d\lambda \biggr]\,dr\,ds
\\
&&\qquad = \biggl(\int_{\R^d} \Gamma_\delta(\lambda) \exp
\biggl\{-{1\over2}\llvert \lambda \rrvert ^2 \biggr\}\,d
\lambda \biggr)\int _0^T\!\!\int
_0^T \llvert r-s\rrvert ^{-\alpha_0-({1/2})(\alpha+\delta)}\,dr\,ds,
\end{eqnarray*}
where the second equality follows from the Fubini theorem and
integration substitution, and the third equality follows from the fact
that $\Gamma_\delta(C\lambda)=C^{-(d-\alpha-\delta)}\Gamma_\delta
(\lambda)$
for $C>0$ and $\lambda\in\R^d$. It is
easy to see that the $\lambda$-integral
on the right-hand side
is finite. We mention the fact that $\alpha_0+{1\over2}\alpha<1$
under our assumptions. Consequently, one
can make the time-integral finite by making $\delta>0$
sufficiently small so $ \alpha_0+{1\over2}(\alpha+\delta)<1$.

Clearly, $Z_T$ is a continuous and nonnegative process in this case.
By the triangle inequality, for any $S, T>0$
$Z_{S+T}\le Z_S+Z_T'$ where
\begin{eqnarray*}
Z_T'&=& \biggl(\int_{\R\times\R^d}
\Gamma_\delta(\lambda) \biggl\llvert \int_S^{S+T}
| u-s| ^{-{(1+\alpha_0)/2}} \exp \bigl\{i\lambda\cdot B(s)
\bigr\}\,ds \biggr\rrvert ^2\,du\,d
\lambda \biggr)^{1/2}
\\
&=& \biggl(\int_{\R\times\R^d}\Gamma_\delta(\lambda) \biggl
\llvert \int_0^T| u-s|
^{-{(1+\alpha_0)/2}}
\\[-10pt]
&&\hspace*{82pt}{}\times  \exp \bigl\{i\lambda\cdot \bigl(B(T+s)-B(T) \bigr) \bigr\}\,ds
\biggr\rrvert ^2\,du\,d \lambda \biggr)^{1/2}.
\end{eqnarray*}
Consequently, $Z_T'$ is independent of $\{B(s); 0\le s\le T\}$
and $Z_T'\buildrel d\over=Z_T$.
By \cite{Chen}, Theorem 1.3.5, page 21,
$Z_T$ is exponential integrable, and the limit
\[
L\equiv\lim_{T\to\infty}{1\over T}\log
\E_0\exp \{Z_T \} %
\]
exists.

By Brownian scaling, on the other hand,
\[
Z_T= \biggl({T\over t} \biggr)^{(4-\alpha-\delta-2\alpha_0)/4}Z_t.
\]
Applying a suitable variable substitution, we conclude
\begin{eqnarray*}
&& \lim_{a\to\infty}a^{-{4/(\alpha+\delta+2\alpha_0)}}\log\E_0\exp \bigl\{
\beta a^{(4-\alpha-\delta-2\alpha_0)/(\alpha+\delta+2\alpha_0)}Z_t \bigr\}
\\
&&\qquad  =Lt\beta^{4/(4-\alpha-\delta-2\alpha_0)}
\end{eqnarray*}
for any $\beta>0$. By G\"artner--Ellis theorem (Theorem 1.2.4, page~11,
\cite{Chen})
\[
\lim_{a\to\infty}a^{-{4/(\alpha+\delta+2\alpha_0)}}\log\P_0
\{Z_t\ge a\} =-C %
\]
for some $C>0$. Consequently,
%
%e4.7 #&#
\begin{eqnarray}
\E_0 Z_t^{qm}=qm\int_0^\infty
a^{qm-1}\P_0\{Z_t\ge a\}\,da \le C^m
(m!)^{q(\alpha+\delta+2\alpha_0)/4},\nonumber
\\
\eqntext{m=1,2,\ldots.}
\end{eqnarray}
In view of (\ref{con-4}), by the stationary of $u(t,x)$ in $x$, we
obtain the bound
%
%e4.8 #&#
\begin{eqnarray}
\E_0\bigl\llvert u(t,x)-u(t,y)\bigr\rrvert ^m\le
C^m\llvert x-y\rrvert ^{m\delta/2} (m!)^{(2+\alpha+\delta+2\alpha_0)/4} \bigl(\E
u_p(t,0)^m \bigr)^{1/p},\nonumber
\\
\eqntext{m=1,2, \ldots}
\end{eqnarray}
uniformly
for all $x, y\in\R^d$. Notice that
$ {2+\alpha+\delta+2\alpha_0\over4}<1$.
This bound leads to (\ref{con-3}).
\end{pf}

%s5 #&#
\section{Tail probability and proof of the upper bounds}\label{tail}

A central piece of our approach relies on the precise
large deviations for $u(t,x)$.
These kinds of results certainly have their independent values. We list
them as part of the major theorems. Recall our assumption that
$u_0(x)=1$.

%
%th5.1 #&#
\begin{theorem}\label{th-4prime} Under the assumption of Theorem \ref{th-1},
%
%
%e5.1 #&#
\begin{eqnarray}
&& \lim_{a\to\infty}a^{-2}\log\P \bigl\{\log
u(t,0)\ge\lambda a \bigr\}
\nonumber\\[-8pt]\label{up-1}  \\[-8pt]\nonumber
&&\qquad  =-{\lambda^2\over2\theta^2} \biggl(\gamma(0)\int_0^t\!\!\int_0^t
\gamma_0(r-s)\,dr\,ds \biggr)^{-1},
\\
&& \lim_{a\to\infty}a^{-2}\log\P \Bigl\{\log
\max_{x\in D}u(t,x)\ge\lambda a \Bigr\}
\nonumber\\[-8pt]\label{up-2}  \\[-8pt]\nonumber
&&\qquad  =-{\lambda^2\over2\theta^2}
\biggl(\gamma(0)\int_0^t\!\!\int
_0^t\gamma_0(r-s)\,dr\,ds
\biggr)^{-1}\nonumber
\end{eqnarray}
for any $t>0$, $\lambda>0$ and bounded domain $D\subset\R^d$.
\end{theorem}

%
%th5.2 #&#
\begin{theorem}\label{th-5} Under the assumption of Theorem \ref{th-2},
%
%
%e5.2 #&#
\begin{eqnarray}
&&\lim_{a\to\infty}a^{-{(4-\alpha)/2}}\log\P \bigl\{
\log u(t,0)\ge \lambda a \bigr\}
\nonumber\\[-8pt]\label{up-3} \\[-8pt]\nonumber
&&\qquad =-{4\over\theta^2} \biggl({2-\alpha\over{\mathcal E}(\alpha_0, d, \gamma)}
\biggr)^{(2-\alpha)/2} \biggl({\lambda\over4-\alpha} \biggr)^{(4-\alpha)/2}
t^{-{(4-\alpha-2\alpha_0)/2}},
\nonumber
\\
&&\lim_{a\to\infty} a^{-{(4-\alpha)/2}} \log\P \Bigl\{
\log\max_{x\in D}u(t,x)\ge\lambda a \Bigr\}
\nonumber\\[-8pt] \label{up-4} \\[-8pt]
&&\qquad = -{4\over\theta^2} \biggl({2-\alpha\over{\mathcal E}(\alpha_0, d, \gamma)}
\biggr)^{(2-\alpha)/2} \biggl({\lambda\over4-\alpha} \biggr)^{(4-\alpha)/2}
t^{-{(4-\alpha-2\alpha_0)/2}}
\nonumber
\end{eqnarray}
for any $t>0$, $\lambda>0$ and bounded domain $D\subset\R^d$,
where ${\mathcal E}(\alpha_0, d, \gamma)$ is the variation given
in (\ref{intro-22}).
\end{theorem}

%
%th5.3 #&#
\begin{theorem}\label{th-5prime} Under the assumption of Theorem \ref{th-2prime},
%
%
%e5.3 #&#
%e5.4 #&#
\begin{eqnarray}
&& \lim_{a\to\infty}a^{-3/2}\log\P \bigl\{\log
u(t,0)\ge\lambda a \bigr\}
\nonumber\\[-8pt]\label{up-3prime}  \\[-8pt]\nonumber
&&\qquad  =-{4\over\theta^2} \sqrt{1\over{\mathcal E}(\alpha_0, 1, \delta_0)}
\biggl({\lambda\over3} \biggr)^{3/2} t^{-{(3-2\alpha_0)/2}},
\\
&& \lim_{a\to\infty} a^{-3/2} \log\P \Bigl\{
\log \max_{x\in D}u(t,x)\ge\lambda a \Bigr\}
\nonumber\\[-8pt]\label{up-4prime}  \\[-8pt]
&&\qquad =-
{4\over\theta^2} \sqrt{1\over{\mathcal E}(\alpha_0, 1, \delta_0)} \biggl(
{\lambda\over3} \biggr)^{3/2} t^{-{(3-2\alpha_0)/2}}\nonumber
\end{eqnarray}
for any $t>0$, $\lambda>0$ and bounded domain $D\subset\R^d$.
\end{theorem}

%
%th5.4 #&#
\begin{theorem}\label{th-6} Under the assumption of Theorem \ref{th-3},
%
%
%e5.5 #&#
\begin{eqnarray}
&& \lim_{a\to\infty}a^{-{(4-\alpha)/2}}\log\P \bigl\{\log
u(t,0)\ge \lambda a \bigr\}
\nonumber\\[-8pt]\label{up-5}  \\[-8pt]\nonumber
&&\qquad  =-{4\over\theta^2} \biggl(
{2-\alpha\over t{\mathcal E}(d, \gamma)} \biggr)^{(2-\alpha)/2} \biggl({\lambda\over4-\alpha}
\biggr)^{(4-\alpha)/2},
\\
&& \lim_{a\to\infty} a^{-{(4-\alpha)/2}} \log\P \Bigl\{
\log\max_{x\in D}u(t,x)\ge\lambda a \Bigr\}
\nonumber\\[-8pt]\label{up-6}  \\[-8pt]\nonumber
&&\qquad  =-
{4\over\theta^2} \biggl({2-\alpha\over t{\mathcal E}(d, \gamma)} \biggr)^{(2-\alpha)/2}
\biggl({\lambda\over4-\alpha} \biggr)^{(4-\alpha)/2}\nonumber
\end{eqnarray}
for any $t>0$, $\lambda>0$ and bounded domain $D\subset\R^d$,
where ${\mathcal E}(d, \gamma)$ is the variation given
in (\ref{intro-27}).
\end{theorem}

%
%th5.5 #&#
\begin{theorem}\label{th-7} When $\gamma_0(\cdot)=\delta_0(\cdot)$,
$\gamma(\cdot)=\delta_0(\cdot)$ and $\alpha=d=1$,
%
%
%e5.6 #&#
%e5.7 #&#
\begin{eqnarray}
\label{up-7} \lim_{a\to\infty} a^{-{3/2}} \log\P \bigl\{
\log u(t,0)\ge\lambda a \bigr\} &=&-{4\over\theta^2} \biggl(
{6\over t} \biggr)^{1/2} \biggl({\lambda\over3}
\biggr)^{3/2},
\\
\label{up-8} \lim_{a\to\infty} a^{-{3/2}} \log\P \Bigl\{
\log\max_{x\in D}u(t,x)\ge\lambda a \Bigr\} &=&-
{4\over\theta^2} \biggl({6\over t} \biggr)^{1/2}
\biggl({\lambda\over3} \biggr)^{3/2}.
\end{eqnarray}
\end{theorem}

Due to similarity we only prove Theorem \ref{th-5}. By H\"older's inequality,
for any $b>1$,
\[
\bigl(\E u(t,0)^{[b]} \bigr)^{1/[b]}\le \bigl(\E
u(t,0)^b \bigr)^{1/b}\le \bigl(\E u(t,0)^{[b]+1}
\bigr)^{1/([b]+1)}. %
\]
Thus the limit in (\ref{mo-3}) (Proposition \ref{prop-2})
can be extended to noninteger $m$. So (\ref{mo-3})
can be re-written as
\begin{eqnarray*}
&& \lim_{a\to\infty}a^{-{(4-\alpha)/2}}\log\E\exp \bigl\{\beta
a^{(2-\alpha)/2} \log u(t,0) \bigr\}
\\
&&\qquad  =\beta^{(4-\alpha)/(2-\alpha)} \biggl(
{\theta^2\over2} \biggr)^{2/(2-\alpha)} t^{(4-\alpha-2\alpha_0)/(2-\alpha)} {\mathcal E}(
\alpha_0, d,\gamma) %
\end{eqnarray*}
for every $\beta>0$.

We now face a problem in using the G\"artner--Ellis theorem: the
exponential moment asymptotics is established only for $\beta>0$,
and the random variable $\log u(t,0)$ takes negative values with
positive probability. To resolve this problem, we
notice that $\P\{u(t,0)\ge1\}>0$ and
\begin{eqnarray*}
&&\E\exp \bigl\{\beta a^{(2-\alpha)/2} \log u(t,0) \bigr\}
\\
&&\qquad =\E \bigl[\exp \bigl\{\beta a^{(2-\alpha)/2} \log u(t,0) \bigr\};
u(t,0)< 1 \bigr]
\\
&&\quad\qquad{} +\E \bigl[\exp \bigl\{\beta a^{(2-\alpha)/2} \log u(t,0) \bigr
\}; u(t,0)\ge1 \bigr]
\\
&&\qquad \le 1+\E \bigl[\exp \bigl\{\beta a^{(2-\alpha)/2} \log u(t,0) \bigr\};
 u(t,0)\ge1 \bigr].
\end{eqnarray*}
We have that for any $\beta>0$
\begin{eqnarray*}
&& \liminf_{a\to\infty}a^{-{(4-\alpha)/2}}\log\E \bigl[ \exp \bigl\{\beta
a^{(2-\alpha)/2} \log u(t,0) \bigr\} \mid u(t,0)\ge1 \bigr]
\\
&&\qquad  \ge
\beta^{(4-\alpha)/(2-\alpha)} \biggl({\theta^2\over2} \biggr)^{2/(2-\alpha)}
t^{(4-\alpha-2\alpha_0)/(2-\alpha)} {\mathcal E}(\alpha_0, d,\gamma). %
\end{eqnarray*}
On the other hand, by the bound
\begin{eqnarray*}
&& \E \bigl[ \exp \bigl\{\beta a^{(2-\alpha)/2} \log u(t,0) \bigr\} \mid u(t,0)
\ge1 \bigr]
\\
&&\qquad \le \bigl(\P \bigl\{u(t,0)\ge1 \bigr\} \bigr)^{-1}\E \exp
\bigl\{\beta a^{(2-\alpha)/2} \log u(t,0) \bigr\}, %
\end{eqnarray*}
we have that for any $\beta>0$
\begin{eqnarray*}
&& \limsup_{a\to\infty}a^{-{(4-\alpha)/2}}\log\E \bigl[ \exp \bigl\{\beta
a^{(2-\alpha)/2} \log u(t,0) \bigr\} \mid u(t,0)\ge1 \bigr]
\\
&&\qquad  \le
\beta^{(4-\alpha)/(2-\alpha)} \biggl({\theta^2\over2} \biggr)^{2/(2-\alpha)}
t^{(4-\alpha-2\alpha_0)/(2-\alpha)} {\mathcal E}(\alpha_0, d,\gamma). %
\end{eqnarray*}
Thus
\begin{eqnarray*}
&& \lim_{a\to\infty}a^{-{(4-\alpha)/2}}\log\E \bigl[ \exp \bigl\{\beta
a^{(2-\alpha)/2} \log u(t,0) \bigr\} \mid u(t,0)\ge1 \bigr]
\\
&&\qquad  =
\beta^{(4-\alpha)/(2-\alpha)} \biggl({\theta^2\over2} \biggr)^{2/(2-\alpha)}
t^{(4-\alpha-2\alpha_0)/(2-\alpha)} {\mathcal E}(\alpha_0, d,\gamma). %
\end{eqnarray*}
By the G\"artner--Ellis
theorem for nonnegative random variables (Theorem 1.2.4, page 11, \cite{Chen}),
\begin{eqnarray*}
&&\lim_{a\to\infty}a^{-{(4-\alpha)/2}}\log \P \bigl\{\log u(t,0)\ge
\lambda a \mid u(t,0)\ge1 \bigr\}
\\
&&\qquad = -\sup_{\beta>0} \biggl\{\beta\lambda- \beta^{(4-\alpha)/(2-\alpha)}
\biggl({\theta^2\over2} \biggr)^{2/(2-\alpha)} t^{(4-\alpha-2\alpha_0)/(2-\alpha)} {\mathcal
E}(\alpha_0, d,\gamma) \biggr\}
\\
&&\qquad =-{4\over\theta^2} \biggl({2-\alpha\over{\mathcal E}(\alpha_0, d, \gamma)}
\biggr)^{(2-\alpha)/2} \biggl({\lambda\over4-\alpha} \biggr)^{(4-\alpha)/2}
t^{-{(4-\alpha-2\alpha_0)/2}}.
\end{eqnarray*}
Therefore, (\ref{up-3}) follows from the fact that
\[
\P \bigl\{\log u(t,0)\ge\lambda a \bigr\}= \P \bigl\{u(t,0)\ge1 \bigr\}\cdot \P
\bigl\{\log u(t,0)\ge\lambda a \mid u(t,0)\ge1 \bigr\}. %
\]

It remains to prove (\ref{up-4}). By (\ref{up-3}) and
the stationary of $u(t,x)$ in $x$, we only
need to prove the upper bound. Without loss of generality,
we may assume that $0\in D$. Notice that
\[
\sup_{x\in D}\bigl\llvert u(t,x)-u(t, 0)\bigr\rrvert
^m\le \operatorname{diam}(D)^{m\delta/2} \sup_{x\in D}
\biggl\llvert {u(t,x)-u(t, 0)\over\llvert  x\rrvert  ^{\delta
/2}} \biggr\rrvert ^m,
\]
where $\delta>0$ is determined by (\ref{con-3}) in Lemma \ref{le-5}.
By (\ref{con-3})
\begin{eqnarray*}
&& \limsup_{m\to\infty}m^{-{(4-\alpha)/(2-\alpha)}}\log\E\sup_{x\in D}
\bigl\llvert u(t,x)-u(t, 0)\bigr\rrvert ^m
\\
&&\qquad \le{1\over p}
\limsup_{m\to\infty}m^{-{(4-\alpha)/(2-\alpha)}}\log\E u_p(t,0)^m
\end{eqnarray*}
for any $p>1$ with $q\equiv p(p-1)^{-1}$ being
an even number. Here we point out that $u_p(t,x)$ is
the solution of the parabolic
Anderson equation (\ref{intro-1}) satisfying the assumption given in
Theorem \ref{th-2} with $\theta$ being replaced by $\theta p$.
By Proposition \ref{prop-2}, the $\limsup$ on the right-hand side is
equal to
\[
\biggl({(p\theta)^2\over2} \biggr)^{2/(2-\alpha)} t^{(4-\alpha-2\alpha_0)/(2-\alpha)} {
\mathcal E}(\alpha_0, d,\gamma). %
\]
Since $p>1$ can be made arbitrarily close to 1, we conclude that
\begin{eqnarray*}
&& \limsup_{m\to\infty}m^{-{(4-\alpha)/(2-\alpha)}}\log\E\sup_{x\in D}
\bigl\llvert u(t,x)-u(t, 0)\bigr\rrvert ^m
\\
&&\qquad  \le \biggl(
{\theta^2\over2} \biggr)^{2/(2-\alpha)} t^{(4-\alpha-2\alpha_0)/(2-\alpha)} {\mathcal E}(
\alpha_0, d, \gamma). %
\end{eqnarray*}
Using Chebyshev's inequality instead of the G\"artner--Ellis theorem,
%
%
%e5.8 #&#
\begin{eqnarray}
\label{up-9}
\qquad && \limsup_{a\to\infty}a^{-{(4-\alpha)/2}}\log \P \Bigl
\{ \log\sup_{x\in D} \bigl\llvert u(t,x)-u(t, 0)\bigr\rrvert \ge
\lambda a \Bigr\}\nonumber
\\
&&\qquad \le -\sup_{\beta>0} \biggl\{\beta\lambda- \beta^{(4-\alpha)/(2-\alpha)}
\biggl({\theta^2\over2} \biggr)^{2/(2-\alpha)} t^{(4-\alpha-2\alpha_0)/(2-\alpha)} {\mathcal
E}(\alpha_0, d,\gamma) \biggr\}
\\
&&\qquad =-{4\over\theta^2} \biggl({2-\alpha\over{\mathcal E}(\alpha_0, d, \gamma)}
\biggr)^{(2-\alpha)/2} \biggl({\lambda\over4-\alpha} \biggr)^{(4-\alpha)/2}
t^{-{(4-\alpha-2\alpha_0)/2}}.
\nonumber
\end{eqnarray}

By the triangle inequality,
\[
\sup_{x\in D}u(t,x)\le u(t,0)+\sup_{x\in D}
\bigl\llvert u(t,x)-u(t, 0)\bigr\rrvert. %
\]
Hence,
\begin{eqnarray*}
&&\log\sup_{x\in D}u(t,x)\le\log \Bigl(u(t,0)+\sup
_{x\in D} \bigl\llvert u(t,x)-u(t, 0)\bigr\rrvert \Bigr)
\\[-1pt]
&&\qquad \le\log2+\max \Bigl\{\log u(t,0), \log\sup_{x\in D}
\bigl\llvert u(t,x)-u(t, 0)\bigr\rrvert \Bigr\}.
\end{eqnarray*}
For any $0<\lambda'<\lambda$, therefore,
\begin{eqnarray*}
&& \P \Bigl\{\log\max_{x\in D}u(t,x)\ge\lambda a \Bigr\}
\\[-1pt]
&&\qquad \le\P
\bigl\{\log u(t,0)\ge\lambda' a \bigr\}+ \P \Bigl\{\log\max
_{x\in D}\bigl\llvert u(t,x)-u(t,0)\bigr\rrvert \ge
\lambda' a \Bigr\} %
\end{eqnarray*}
for large $a$. Thus
\begin{eqnarray*}
&&\limsup_{a\to\infty}a^{-{(4-\alpha)/2}}\log\P \Bigl\{\log \max
_{x\in D}u(t,x)\ge\lambda a \Bigr\}
\\[-1pt]
&&\qquad \le \max \Bigl\{\limsup_{a\to\infty}a^{-{(4-\alpha)/2}}\log\P \bigl\{
\log u(t, 0) \ge\lambda' a \bigr\},
\\[-1pt]
&&\hspace*{58pt} \limsup_{a\to\infty}a^{-{(4-\alpha)/2}}\log\P \Bigl\{\log
\max_{x\in D}\bigl\llvert u(t,x)-u(t,0)\bigr\rrvert \ge
\lambda' a \Bigr\} \Bigr\}
\\[-1pt]
&&\qquad \le-{4\over\theta^2} \biggl({2-\alpha\over{\mathcal E}(\alpha_0, d,
\gamma)}
\biggr)^{(2-\alpha)/2} \biggl({\lambda'\over4-\alpha} \biggr)^{(4-\alpha)/2}
t^{-{(4-\alpha-2\alpha_0)/2}},
\end{eqnarray*}
where the last step follows from (\ref{up-5}) and (\ref{up-9}). Since
$\lambda'$ can be arbitrarily close to $\lambda$, we have finally established
the upper bound requested by
(\ref{up-4}).
%\end{pf}

Having Theorems \ref{th-4prime}--\ref{th-7}
installed, we are ready to prove the upper bounds in Theorems \ref
{th-1}, \ref{th-2}, \ref{th-2prime}, \ref{th-3}
and \ref{th-4}. Again, due to similarity we only
prove the upper bound requested by Theorem \ref{th-3}.
That is,
%
%
%e5.9 #&#
\begin{eqnarray}
\label{up-10}
&& \limsup_{R\to\infty}(\log R)^{-{2/(4-\alpha)}}\log\max
_{\llvert  x\rrvert
\le R}u(t,x)
\nonumber\\[-9pt]\\[-9pt]\nonumber
&&\qquad  \le{4-\alpha\over4} \biggl(
{4t{\mathcal E}(d, \gamma)\over2-\alpha} \biggr)^{(2-\alpha)/(4-\alpha)} \theta^{4/(4-\alpha)}d^{2/(4-\alpha)}
\qquad\mbox{a.s.}
\end{eqnarray}

To this end, we set ${\mathcal N}_R=\Z^d\cap B(0,R)$ and write $Q=[-1,1]^d$.
Notice that
\[
\max_{\llvert  x\rrvert  \le R}u(t,x)\le\max_{z\in{\mathcal N}_R}\max
_{x\in
z+Q}u(t,x). %
\]
For any $\lambda>0$ satisfying
\[
\lambda>{4-\alpha\over4} \biggl({4t{\mathcal E}(d, \gamma)\over2-\alpha}
\biggr)^{(2-\alpha)/(4-\alpha)} \theta^{4/(4-\alpha)}d^{2/(4-\alpha)}, %
\]
by stationarity of $u(t,x)$ in $x$
\begin{eqnarray*}
&& \P \Bigl\{\log\max_{\llvert  x\rrvert  \le R}u(t,x)\ge \lambda(\log
R)^{2/(4-\alpha)} \Bigr\}
\\
&&\qquad  \le\#({\mathcal N}_R)\P \Bigl\{\log\max
_{x\in Q}u(t,x)\ge \lambda(\log R)^{2/(4-\alpha)} \Bigr\}.
\end{eqnarray*}
By (\ref{up-6}) in Theorem \ref{th-6} there is a $\delta>0$ such that
\[
\P \Bigl\{\log\max_{x\in Q}u(t,x)\ge \lambda(\log
R)^{2/(4-\alpha)} \Bigr\}\le\exp \bigl\{-(d+\delta)\log R \bigr\} %
\]
as $R$ is sufficiently large. Consequently,
\[
\P \Bigl\{\log\max_{\llvert  x\rrvert  \le R}u(t,x)\ge \lambda(\log
R)^{2/(4-\alpha)} \Bigr\}\le CR^{-\delta} %
\]
with the constant $C>0$ independent of $R$. With this bound
\[
\sum_{n=1}^\infty\P \Bigl\{\log\max
_{\llvert  x\rrvert  \le2^n}u(t,x)\ge \lambda \bigl(\log2^n
\bigr)^{2/(4-\alpha)} \Bigr\}<\infty. %
\]
By the Borel--Cantelli lemma
\[
\limsup_{n\to\infty} \bigl(\log2^n \bigr)^{-{2/(4-\alpha)}}
\log\max_{\llvert  x\rrvert  \le2^n}u(t,x) \le\lambda\qquad\mbox{a.s.} %
\]
The $\limsup$ can be extended from the sequence $2^n$ to $R$ due to
the monotonicity of the quantity
$ \log\max_{\llvert  x\rrvert  \le R}u(t,x)$ in $R$.
Finally, (\ref{up-10}) follows from the fact that $\lambda$ can be arbitrarily
close to the limit value appearing on the right-hand side of~(\ref{up-10}).
%\end{pf}

%s6 #&#
\section{Link to the long-term asymptotics: The case of time independence}\label{time}

A~classic quenched law (Theorem 5.1, \cite{CM-1}) by Carmona and Molchanov
stated that
for a homogeneous and time-independent Gaussian potential
$V(x)$ whose covariance function
$\gamma(x)$ satisfies the conditions comparable to the ones assumed in
Theorem \ref{th-1},
%
%
%e6.1 #&#
\begin{eqnarray}
\label{time-1} \lim_{t\to\infty}{1\over t\sqrt{\log t}}\log
\E_0\exp \biggl\{\theta\int_0^tV
\bigl(B(s) \bigr)\,ds \biggr\}=\theta\sqrt {2d\gamma(0)} \qquad\mbox{a.s.}
\end{eqnarray}

In his recent work, Chen \cite{C-1} considers the case of the time
independent Gaussian field $V(x)$ with the covariance
function $\gamma(\cdot)$ in the forms given
in Table~\ref{table1}. More specifically, under the assumption $0<\alpha<2\wedge d$,
and for the $\gamma(\cdot)$ of types (I) and (II) (labeled in Table~\ref{table1})
(Corollary 1.2 and Theorem 1.3, \cite{C-1}),
\begin{eqnarray}
\label{time-2} &&\lim_{t\to\infty}t^{-1}(\log
t)^{-{2/(4-\alpha)}}\log \E_0\exp \biggl\{\theta\int
_0^tV \bigl(B(s) \bigr)\,ds \biggr\}
\nonumber\\[-8pt]\\[-8pt]\nonumber
&&\qquad = {4-\alpha\over4} \biggl({4{\mathcal E}(d, \gamma)\over2-\alpha}
\biggr)^{(2-\alpha)/(4-\alpha)}\theta^{4/(4-\alpha)}d^{2/(4-\alpha)} \qquad\mbox{a.s.}
\nonumber
\end{eqnarray}
When $d=1$ and $\gamma(\cdot)=\delta_0(\cdot)$ (Theorem 1.4, \cite{C-1}),
%
%
%e6.2 #&#
\begin{eqnarray}
\label{time-2prime} \lim_{t\to\infty}t^{-1}(\log
t)^{-2/3} \log \E_0\exp \biggl\{\theta\int
_0^tV \bigl(B(s) \bigr)\,ds \biggr\} =
{3\over4}\root3\of{2\over3}\qquad\mbox{a.s.}
\end{eqnarray}
We mention that the right-hand side of (\ref{time-2}) was initially
given in terms of the best constant $\kappa(\gamma, d)$ of the Soblev-type
inequality
\[
\int_{\R^d\times\R^d}\gamma(x-y)f^2(x)f^2(y)\,dx\,dy
\le C\llVert f\rrVert _2^{4-\alpha} \llVert \nabla f\rrVert
_2^\alpha \qquad f\in W^{1,2} \bigl(
\R^d \bigr) %
\]
and can be switched into the current form, thanks to the identity
\[
{\mathcal E}(d, \gamma)={2-\alpha\over2}\alpha^{\alpha/(2-\alpha)} \kappa(
\gamma, d)^{2/(2-\alpha)} %
\]
which can be derived in the same way as (7.3) in \cite{CHSX}.

The striking resemblance of the pairs (\ref{intro-29})
versus (\ref{time-1}),
(\ref{intro-30}) versus (\ref{time-2}) and
(\ref{intro-31}) versus (\ref{time-2prime})
suggests a possible link
between the time asymptotics and the spatial asymptotics. In this section
we explore this problem by providing an alternative treatment to the long-term
asymptotics. For similarity, we only consider~(\ref{time-2}).

For the sake of simplicity we assume that $t$ goes $\infty$ along the integer
points.
Given $R>0$, define $\tau(R)=\inf\{s\ge0; \llvert  B(s)\rrvert  \ge
R\}$.
For any function $R(t)\uparrow\infty$ ($t\to\infty$), by Markov's property
%
%
%e6.3 #&#
\begin{eqnarray}
\label{time-3}
&& \E_0 \biggl[\exp \biggl\{\theta\int
_0^tV \bigl(B(s) \bigr)\,ds \biggr\};
\tau \bigl(R(t) \bigr)\ge t \biggr]
\nonumber\\[-9pt]\\[-9pt]\nonumber
&&\qquad  \le \biggl(\max_{\llvert  x\rrvert  \le R(t)}
\E_x\exp \biggl\{\theta\int_0^1
V \bigl(B(s) \bigr)\,ds \biggr\} \biggr)^t.
\end{eqnarray}
Applying (\ref{intro-30}) we have
\begin{eqnarray}
\label{time-4} \qquad&&\limsup_{t\to\infty}t^{-1} \bigl(\log
R(t) \bigr)^{-{2/(4-\alpha)}}\log \E_0 \biggl[\exp \biggl\{\theta\int
_0^tV \bigl(B(s) \bigr)\,ds \biggr\}; \tau \bigl(R(t) \bigr)\ge t \biggr]\nonumber
\nonumber\\[-9pt]\\[-9pt]\nonumber
&&\qquad \le {4-\alpha\over4} \biggl({4{\mathcal E}(d, \gamma)\over2-\alpha}
\biggr)^{(2-\alpha)/(4-\alpha)}
\theta^{4/(4-\alpha)}d^{2/(4-\alpha)} \qquad\mbox{a.s.}\nonumber
\end{eqnarray}

Further, let $R_k(t)=t(\log t)^{k+1}$ ($k=0,1,2,\ldots$).
\begin{eqnarray*}
&&\E_0\exp \biggl\{\theta\int_0^tV
\bigl(B(s) \bigr)\,ds \biggr\}
\\[-1pt]
&&\qquad = \E_0 \biggl[\exp \biggl\{\theta\int_0^tV
\bigl(B(s) \bigr)\,ds \biggr\}; \tau \bigl(R_0(t) \bigr)\ge t \biggr]
\\[-1pt]
&&\quad\qquad{} +\sum_{k=1}^\infty\E_0
\biggl[\exp \biggl\{\theta\int_0^tV \bigl(B(s)
\bigr)\,ds \biggr\}; \tau \bigl(R_{k-1}(t) \bigr)< t\le\tau
\bigl(R_k(t) \bigr) \biggr]
\\[-1pt]
&&\qquad \le\E_0 \biggl[\exp \biggl\{\theta\int_0^tV
\bigl(B(s) \bigr)\,ds \biggr\}; \tau \bigl(R_0(t) \bigr)\ge t \biggr]
\\[-1pt]
&&\quad\qquad{} +\sum_{k=1}^\infty \biggl(
\E_0 \biggl[\exp \biggl\{2\theta\int_0^t
V \bigl(B(s) \bigr)\,ds \biggr\}; \tau \bigl(R_k(t) \bigr)\ge t
\biggr] \biggr)^{1/2}
\\[-1pt]
&&\hspace*{10pt}\quad\qquad{}\times  \bigl(\P_0 \bigl\{\tau
\bigl(R_{k-1}(t) \bigr) < t \bigr\} \bigr)^{1/2}.
\end{eqnarray*}
By Gaussian tail
%
%e6.4 #&#
\begin{eqnarray}
\P_0 \bigl\{\tau \bigl(R_{k-1}(t) \bigr)< t \bigr\} \le
\exp \biggl\{-C{R_{k-1}(t)^2\over t} \biggr\}=\exp \bigl\{-Ct(\log
t)^{2k} \bigr\},\nonumber
\\
\eqntext{k=1,2,\ldots.}
\end{eqnarray}
Together with (\ref{time-4}), this shows that the infinite series on
the right-hand side of the decomposition is negligible. Applying
(\ref{time-4})
to the first term [with $R(t)=R_0(t)$] on the right-hand
side of the decomposition leads to the upper bound requested by~(\ref{time-2}).

Relation (\ref{time-3}) is reversible with some nonsubstantial but
technically involved modification, so (\ref{intro-30}) also applies to
the lower bound for (\ref{time-2}).
We skip this part of the argument.

%
%re6.1 #&#
\begin{remark}\label{remark-2}
An asymptotic bound similar to (\ref{time-4})
can be extended to the setting of time-dependence with some obvious
modification. However, it is unlikely to
be sharp in
the settings given in Table~\ref{table1} (with $\alpha_0>0$, of course). Compared
with the case of time independence, much less
is known about the quenched long-term asymptotics in the setting
of time-dependence.
\end{remark}

%sA #&#
\begin{appendix}
\section*{Appendix}\label{appendix}

%sA.1 #&#
\subsection{Feynman--Kac bounds}

\setcounter{equation}{0}
\setcounter{theorem}{0}

For any open domain $D\in\R^d$, define
${\mathcal F}_d(D)$ as the class of the functions $g$ supported
in $D$ such that $\llVert  g\rrVert  _2=1$ and $\llVert  \nabla g\rrVert  _2<\infty$.
Write
%
%
%eA.1 #&#
\begin{eqnarray}
\label{v-1} \tau_D=\inf \bigl\{s\ge0; B(s)\notin D
\bigr\}.
\end{eqnarray}
For a function
$f$ defined on $D$, set
%
%
%eA.2 #&#
\begin{eqnarray}
\label{v-2} \lambda_D(f)=\sup_{g\in{\mathcal F}_d(D)} \biggl
\{ \int_D f(x)g^2(x)\,dx-{1\over2}
\int_D\bigl\llvert \nabla g(x)\bigr\rrvert
^2\,dx \biggr\}.
\end{eqnarray}

%
%leA.1 #&#
\begin{lemma}\label{v-3} Let $t>0$, and let the function
$f(s,x)$ be continuous and bounded on $[0, t]\times \operatorname{cl}(D)$.
Then for any $t>0$,
%
%
%eA.3 #&#
%eA.4 #&#
\begin{eqnarray}
&& \int_D\E_x \biggl[\exp
\biggl\{\int_0^tf \bigl(t-s, B(s) \bigr)\,ds
\biggr\}; \tau_D\ge t \biggr]\,dx
\nonumber\\[-8pt]\label{v-4}  \\[-8pt]\nonumber
&&\qquad \le\llvert D\rrvert \exp
\biggl\{\int_0^t\lambda_D
\bigl(f(s,\cdot) \bigr)\,ds \biggr\},
\\
&& \int_D\E_x \biggl[\exp
\biggl\{\int_0^tf \bigl(s, B(s) \bigr)\,ds
\biggr\}; \tau_D\ge t \biggr]\,dx
\nonumber\\[-8pt]\label{v-5}  \\[-8pt]
&&\qquad \le\llvert D\rrvert \exp
\biggl\{\int_0^t\lambda_D
\bigl(f(s, \cdot) \bigr)\,ds \biggr\}.\nonumber
\end{eqnarray}
\end{lemma}

\begin{pf} By the Feynman--Kac formula (e.g., Theorem 2.3, page 133,
\cite{Freidlin}
with $g(s,x)=0$), the function
\[
u(s,x)=\E_x \biggl[\exp \biggl\{\int_0^sf
\bigl(s-u, B(u) \bigr)\,du \biggr\}; \tau_D\ge s \biggr],
\qquad x\in D %
\]
solves the initial-boundary problem
\[
\qquad \cases{ \partial_su(s,x) ={1\over2}\Delta u(s,x)+f(s,x)u(s,x), &\quad $(s,x)\in(0,t]\times D$,
\vspace*{3pt}\cr
u(s,x)=0, &\quad $x\in\partial D$,
\vspace*{3pt}\cr
u(0, x)=1, &\quad $x\in D$.}
\]
Hence
\begin{eqnarray*}
{d\over ds}\int_Du^2(s,x)\,dx&=&2
\int_D u(s,x)\partial_su(s,x) \,dx
\\
&=&2 \biggl\{\int_Df(s,x)u^2(s,x)\,dx -
{1\over2}\int_D\bigl\llvert
\nabla_x u(s,x)\bigr\rrvert ^2\,dx \biggr\}
\\
&\le&2\lambda_D \bigl(f(s,\cdot) \bigr)\int_Du^2(s,x)\,dx.
\end{eqnarray*}
Notice that the function
\[
U(s)=\int_Du^2(s,x)\,dx %
\]
has the initial value $U(0)=\llvert   D\rrvert  $.
Thus by Gronwall's inequality
\[
\int_Du^2(t,x)\,dx\le \llvert D\rrvert \exp
\biggl\{2\int_0^t\lambda_D
\bigl(f(s, \cdot) \bigr)\,ds \biggr\}. %
\]
Therefore, (\ref{v-4}) follows from the Cauchy--Schwarz inequality:
\[
\int_Du(t,x)\,dx\le\sqrt{\llvert D\rrvert } \biggl\{\int
_Du^2(t,x)\,dx \biggr\}^{1/2}.
\]

Replacing $f(s, x)$ by $f_t(s, x)=f(t-s, x)$ in (\ref{v-4})
leads to (\ref{v-5}).
\end{pf}

%sA.2 #&#
\subsection{A lemma on the large deviations}

Let $\{X_m\}$ be a sequence of nonnegative random variables and
$b_m$ be a sequence of positive numbers such that $b_m\to\infty$
as $m\to\infty$.

%
%leA.2 #&#
\begin{lemma}\label{A-1} Assume that there is $p>1$ and $C_0>0$
such that for any $\beta>0$,
%
%
%eA.5 #&#
%eA.6 #&#
\begin{eqnarray}
\label{A-2}
\limsup_{m\to\infty}{1\over b_m}\log\E
\exp \{ \beta X_m \} &\le& C_0 \beta^p,
\\
\label{A-3} \qquad \liminf_{m\to\infty}{1\over b_m}\log\E
\exp \bigl\{ \beta b_m^{1/2} X_m^{1/2}
\bigr\} &\ge& {p+1\over p}(pC_0)^{1/(p+1)} \biggl(
{\beta\over2} \biggr)^{(2p)/(p+1)}.
\end{eqnarray}
Then we have
%
%
%eA.7 #&#
\begin{eqnarray}
\label{A-4} \lim_{m\to\infty}{1\over b_m}\log\E\exp
\{ \beta X_m \} = C_0 \beta^p\qquad
\forall\beta>0.
\end{eqnarray}

The same claim holds if we weaken the first assumption (\ref{A-2}) to
%
%
%eA.8 #&#
\begin{eqnarray}
\label{A-5} \limsup_{m\to\infty}{1\over b_m}\log\E
\exp \{ \beta X_m \} <\infty\qquad \forall\beta>0
\end{eqnarray}
and strengthen the second assumption (\ref{A-3}) into
%
%
%eA.9 #&#
\begin{eqnarray}
\label{A-6}
&& \lim_{m\to\infty}{1\over b_m}\log\E\exp
\bigl\{ \beta b_m^{1/2} X_m^{1/2} \bigr
\}
\nonumber\\[-8pt]\\[-8pt]\nonumber
&&\qquad  = {p+1\over p}(pC_0)^{1/(p+1)} \biggl(
{\beta\over2} \biggr)^{(2p)/(p+1)} \qquad \forall\beta>0.
\end{eqnarray}
\end{lemma}

\begin{pf} Due to similarity, we only prove the first claim.
By (\ref{A-2}) and by a standard way of using Chebyshev's
inequality, for any $\lambda>0$,
\begin{eqnarray*}
\limsup_{m\to\infty}{1\over b_m}\log\P\{X_m
\ge\lambda b_m\} &\le&-\sup_{\beta>0} \bigl\{\lambda
\beta-C_0 \beta^p \bigr\}
\\
& =&-{p-1\over p}(C_0p)^{-{1/(p-1)}}
\lambda^{p/(p-1)},
\end{eqnarray*}
and for any $\beta>0$,
\[
\limsup_{m\to\infty}{1\over b_m}\log\E\exp \bigl\{ \beta
b_m^{1/2} X_m^{1/2} \bigr\} <\infty.
\]
By Varadhan's integral lemma (Lemma 4.3.6, \cite{DZ}),
\begin{eqnarray*}
&&\limsup_{m\to\infty}{1\over b_m} \log\E\exp \bigl\{
\beta b_m^{1/2} X_m^{1/2} \bigr\}
\\
&&\qquad \le\sup_{\lambda>0} \biggl\{\beta\lambda^{1/2}-
{p-1\over p}(C_0p)^{-{1/(p-1)}}\lambda
^{p/(p-1)} \biggr\}
\\
&&\qquad ={p+1\over p}(pC_0)^{1/(p+1)} \biggl(
{\beta\over2} \biggr)^{(2p)/(p+1)}.
\end{eqnarray*}
Together with (\ref{A-3}) and the G\"artner--Ellis theorem (Theorem
1.2.4, page 11,
\cite{Chen}), we have
\begin{eqnarray*}
&&\lim_{m\to\infty}{1\over b_m}\log\P \{X_m
\ge\lambda b_m \}
\\
&&\qquad =-\sup_{\beta>0} \biggl\{\beta\sqrt{\lambda} -
{p+1\over p}(pC_0)^{1/(p+1)} \biggl(
{\beta\over2} \biggr)^{(2p)/(p+1)} \biggr\}
\\
&&\qquad =-{p-1\over p}(C_0p)^{-{1/(p-1)}}
\lambda^{p/(p-1)}.
\end{eqnarray*}
Finally, by Varadhan's integral lemma (Lemma 4.3.6, \cite{DZ})
\begin{eqnarray*}
\lim_{m\to\infty}{1\over b_m}\log\E\exp \{ \beta
X_m \} &=&\sup_{\lambda>0} \biggl\{\beta\lambda -
{p-1\over p}(C_0p)^{-{1/(p-1)}}\lambda^{p/(p-1)}
\biggr\}
\\
&=& C_0 \beta^p.
\end{eqnarray*}\upqed
\end{pf}

%sA.3 #&#
\subsection{Variations}

Recall that for any $\varepsilon>0$ and $\beta>0$,
\begin{eqnarray*}
M_\varepsilon(\beta)&=&\sup_{g\in{\mathcal A}_d} \biggl\{\beta \biggl(
\int_0^1\!\!\int_{\R^d}
\biggl[\int_{\R^d}K_\varepsilon(y-x)g^2(s, y)\,dy
\biggr]^2 \,dx\,ds \biggr)^{1/2}
\\
&&\hspace*{111pt}{} - {1\over2}\int_0^1\!\!\int_{\R^d}\bigl
\llvert \nabla_x g(s,x)\bigr\rrvert ^2\,dx\,dy \biggr\},
\\
M_{\varepsilon, N}(\beta) &=& \sup_{g\in{\mathcal A}_d} \biggl\{ \beta \biggl(
\int_0^1\!\!\int_{[-N, N]^d}
\biggl[\int_{\R^d} Q_N(y-x)g^2(s,y)\,dy
\biggr]^2 \,dx \biggr)^{1/2}
\\
&&\hspace*{126pt}{} -{1\over2} \int_0^1\!\!\int_{\R^d}\bigl
\llvert \nabla_x g(s, x) \bigr\rrvert ^2\,dx\,ds \biggr\},
\end{eqnarray*}
where $K_\varepsilon(x)$ and $Q_N(x)$ are defined in (\ref{mo-11prime}) and
(\ref{mo-12prime}), respectively.

%
%leA.3 #&#
\begin{lemma}\label{va-1} In the settings marked (2) in Table~\ref{table1}, for
any $\varepsilon>0$ and $\beta>0$,
\[
\limsup_{N\to\infty}M_{\varepsilon, N} (\beta)\le M_\varepsilon(
\beta). %
\]
\end{lemma}

\begin{pf} Notice that for any $0\le s\le1$,
\[
\int_{\R^d} Q_N(y-x)g^2(s,y)\,dy =
\int_{\R^d}Q(y-x)\tilde{g}^2(s,y)\,dy,
\]
where
\[
\tilde{g}(x)=\sqrt{\sum_{{\mathbf k}\in\Z^d}g^2(2{\mathbf
k}N+x)},\qquad x\in \R^d. %
\]
Here we recall
\[
Q(x)=K_\varepsilon(x) l \bigl(M^{-1}\llvert x\rrvert \bigr),
\]
where $l$: $\R^+\longrightarrow[0, 1]$ is a smooth function
satisfying the following properties:
$l(u)=1$ for $u\in[0,1]$, $l(u)=0$ for $u\ge3$ and $-1\le l'(u)\le0$
for all $u>0$. By the fact that $Q(\cdot)$ is
supported on $[-3M, 3M]^d$,
\begin{eqnarray*}
M_{\varepsilon, N}(\beta) &=& \sup_{g\in{\mathcal A}_d} \biggl\{ \biggl( \int_0^1\!\!\int_{[-N, N]^d} \biggl[
\int_{[-\widetilde{N}, \widetilde{N}]^d} Q(y-x)\tilde{g}^2(s,y)\,dy
\biggr]^2 \,dx \biggr)^{1/2}
\\
&&\hspace*{134pt}{} -{1\over2} \int_0^1\!\!\int_{\R^d} \bigl
\llvert \nabla_x g(s, x)\bigr\rrvert ^2\,dx\,ds \biggr\},
\end{eqnarray*}
where $\widetilde{N}=N+3M$.
We omit the rest of the proof as it follows from the constructive
argument used in \cite{Chen-1}, Lemma A.1, with some minor
modification.
\end{pf}

We also use the following notation:
\begin{eqnarray*}
M(\beta) &=&\sup_{g\in{\mathcal A}_d} \biggl\{\beta \biggl( \int_0^1\!\!\int_{\R^d} \biggl[
\int_{\R^d}K(y-x)g^2(s, y)\,dy \biggr]^2
\,dx\,ds \biggr)^{1/2}\\
&&\hspace*{106pt}{}  -{1\over2} \int_0^1\!\!\int_{\R^d} \bigl\llvert \nabla_x g(s,x)\bigr
\rrvert ^2\,dx\,dy \biggr\}, %
\\
{\mathcal E}(d, \gamma)&\equiv&\sup_{g\in{\mathcal F}_d} \biggl\{ \int
_{\R^d\times\R^d}\gamma(x-y) g^2(x)g^2(y)\,dx\,dy
\\
&&\hspace*{87pt}{} -
{1\over2}\int_{\R^d}\bigl\llvert \nabla g(x)
\bigr\rrvert ^2\,dx \biggr\}. %
\end{eqnarray*}

%
%leA.4 #&#
\begin{lemma}\label{va-2} Under the assumptions of Theorems \ref{th-3}~and~\ref{th-4},
\[
{\mathcal E}(d, \gamma)={2-\alpha\over2}2^{\alpha/(2-\alpha)} \biggl(
{4M(1)\over4-\alpha} \biggr)^{(4-\alpha)/(2-\alpha)}. %
\]
\end{lemma}

\begin{pf} By (\ref{lo-12}), $M(1)$ can be rewritten as
\begin{eqnarray*}
M(1) &=& \sup_{g\in{\mathcal A}_d} \biggl\{ \biggl(\int_0^1\!\!\int_{\R^d\times\R^d}\gamma (x-y)g^2(s,x)g^2(s,y)\,dx\,dy
\biggr)^{1/2}
\\
&&\hspace*{106pt}{} -{1\over2}\int_0^1\!\!\int_{\R^d}\bigl\llvert \nabla_x g(s,x)\bigr
\rrvert ^2\,dx\,dy \biggr\}. %
\end{eqnarray*}
Define
\begin{eqnarray*}
{\mathcal E}'(d,\gamma)&=&\sup_{g\in{\mathcal A}_d} \biggl\{
\int_0^1\!\!\int_{\R^d\times\R^d}
\gamma(x-y)g^2(s,x)g^2(s,y)\,dx\,dy\,ds
\\
&&\hspace*{93pt}{} - {1\over2}
\int_0^1\!\!\int_{\R^d}\bigl
\llvert \nabla_x g(s,x)\bigr\rrvert ^2\,dx\,ds \biggr\}.
\end{eqnarray*}
Replacing $\gamma_0(s)=\llvert   s\rrvert  ^{-\alpha_0}$ by $\gamma_0(s)=\delta_0(s)$
in (7.4), \cite{CHSX}, we have
\[
{\mathcal E}'(d, \gamma)={2-\alpha\over2}2^{\alpha/(2-\alpha)}
\biggl({4M(1)\over4-\alpha} \biggr)^{(4-\alpha)/(2-\alpha)}. %
\]
Therefore, it remains to show that
\[
{\mathcal E}'(d,\gamma)={\mathcal E}(d,\gamma). %
\]

Indeed, taking $g(s,\cdot)=g(\cdot)\in{\mathcal F}_d$
leads to
${\mathcal E}'(d,\gamma)\ge{\mathcal E}(d,\gamma)$. On the other hand,
by the relation ${\mathcal A}_d= \{g(\cdot, \cdot);
g(s,\cdot)
\in{\mathcal F}_d$ $\forall0\le s\le1 \}$,
\begin{eqnarray*}
{\mathcal E}'(d,\gamma)&\le&\int_0^1
\sup_{g\in{\mathcal A}_d} \biggl\{ \int_{\R^d\times\R^d}
\gamma(x-y)g^2(s,x)g^2(s,y)\,dx\,dy
\\
&&\hspace*{96pt}{}  -{1\over2}
\int_{\R^d}\bigl\llvert \nabla_x g(s,x) \bigr
\rrvert ^2\,dx\,dy \biggr\}\,ds
\\
&=&{\mathcal E}(d,\gamma).
\end{eqnarray*}\upqed
\end{pf}

In connection to the variation ${\mathcal E}(\alpha_0, d, \gamma)$
given in
(\ref{intro-22}), write
\begin{eqnarray*}
{\mathcal E}(0, d,\gamma)&=&\sup_{g\in{\mathcal A}_d} \biggl\{\int_0^1\!\!\int_0^1\!\!\int_{\R^d\times\R^d}\gamma(x-y) g^2(s, x)g^2(r,
y)\,dx\,dy\,dr\,ds
\\
&&\hspace*{118pt}{}-{1\over2}\int_0^1\!\!\int
_{\R^d}\bigl\llvert \nabla_x g(s, x)\bigr\rrvert
^2\,dx\,ds \biggr\}.
\end{eqnarray*}

%
%leA.5 #&#
\begin{lemma}\label{va-3} Under the assumptions of Corollary \ref{co-2},
${\mathcal E}(0, d,\gamma)={\mathcal E}(d,\gamma)$.
\end{lemma}

\begin{pf} The direction of $\ge$ is obvious. We now consider
opposite direction.
Let $g\in{\mathcal A}_d$, and write
\[
\tilde{g}(x)= \biggl(\int_0^1g^2(s,x)\,ds
\biggr)^{1/2}. %
\]
Then $\tilde{g}\in{\mathcal F}_d$ and
\begin{eqnarray*}
&& \int_0^1\!\!\int_0^1\!\int_{\R^d\times\R^d}\gamma(x-y) g^2(s, x)g^2(r,
y)\,dx\,dy\,dr\,ds
\\
&&\qquad  =\int_{\R^d\times\R^d}\gamma(x-y) \tilde{g}^2(x)
\tilde{g}^2(y)\,dx\,dy %
\end{eqnarray*}
and
\[
\nabla\tilde{g}(x)= \biggl(\int_0^1g^2(s,x)\,ds
\biggr)^{-1/2} \int_0^1g(s,x)
\nabla_xg(s,x)\,ds. %
\]
Hence
\begin{eqnarray*}
\bigl\llvert \nabla\tilde{g}(x)\bigr\rrvert &\le& \biggl(\int_0^1g^2(s,x)\,ds
\biggr)^{-1/2}\int_0^1\bigl\llvert
g(s,x) \bigr\rrvert \cdot \bigl\llvert \nabla_xg(s,x)\bigr\rrvert \,ds
\\
&\le&
\biggl(\int_0^1\bigl\llvert
\nabla_xg(s,x) \bigr\rrvert ^2 \,ds \biggr)^{1/2}.
\end{eqnarray*}
Thus
\[
\int_{\R^d}\bigl\llvert \nabla\tilde{g}(x)\bigr\rrvert
^2\,dx \le \int_0^1\!\!\int
_{\R^d} \bigl\llvert \nabla_xg(s,x)\bigr\rrvert
^2 \,dx\,ds. %
\]

Summarizing our estimate,
\begin{eqnarray*}
&& \int_0^1\!\!\int_0^1\!\int_{\R^d\times\R^d}\gamma(x-y) g^2(s, x)g^2(r,
y)\,dx\,dy\,dr\,ds
\\
&&\quad{} -{1\over2} \int_0^1\!\!\int_{\R^d}\bigl\llvert \nabla_xg(s,x)\bigr
\rrvert ^2 \,dx\,ds
\\
&&\qquad \le \int_{\R^d\times\R^d}\gamma(x-y)\tilde{g}^2(x)
\tilde{g}^2(y)\,dx\,dy -{1\over2}\int_{\R^d}
\bigl\llvert \nabla\tilde{g}(x)\bigr\rrvert ^2\,dx\le{\mathcal E}(d,
\gamma).
\end{eqnarray*}
Taking supremum over $g\in{\mathcal A}_d$ on the left-hand side
completes the proof.
\end{pf}
\end{appendix}

\section*{Acknowledgment} The author is grateful to two anonymous referees
for their careful reading of the manuscript and for
making numerous corrections and suggestions.

%\begin{supplement}[id=suppA]
%\slink[doi]{10.1214/00-AOPXXXXSUPP} %[doi,text={...}] - jei reikia
%suskaldyti doi
%\sfilename{aopXXXX\_supp.pdf}
%\end{supplement}

% imsref loaded by linak, 2015-03-23 15:09:29
%

\printaddresses
\end{document}